\documentclass{article}
\usepackage{amsmath}
\usepackage{amsthm}
\usepackage{amssymb}

\theoremstyle{definition}
\newtheorem{defi}{Definition}

\theoremstyle{plain}
\newtheorem{thm}[defi]{Theorem}
\newtheorem{defthm}[defi]{Definition and Theorem}
\newtheorem{cor}[defi]{Corollary}
\newtheorem{lem}[defi]{Lemma}

\newtheorem{prop}[defi]{Proposition}

\newcounter{enuroman}
\renewcommand{\theenuroman}{\roman{enuroman}}

\newcommand{\forces}{\Vdash}
\newcommand{\re}{{\upharpoonright}}

\newcommand{\Bool}{[\![}
\newcommand{\Boor}{]\!]}

\newcommand{\embed}{{<\!\!\circ\;}}

\newcommand{\A}{{\cal A}}
\newcommand{\B}{{\cal B}}

\newcommand{\D}{{\cal D}}

\newcommand{\F}{{\cal F}}

\newcommand{\I}{{\cal I}}

\newcommand{\M}{{\cal M}}
\newcommand{\N}{{\cal N}}

\renewcommand{\P}{{\cal P}}

\newcommand{\SSS}{{\cal S}}
\newcommand{\U}{{\cal U}}

\newcommand{\BB}{{\mathbb B}}

\newcommand{\DD}{{\mathbb D}}
\newcommand{\EE}{{\mathbb E}}

\newcommand{\LL}{{\mathbb L}}

\newcommand{\NN}{{\mathbb N}}
\newcommand{\PP}{{\mathbb P}}
\newcommand{\QQ}{{\mathbb Q}}

\renewcommand{\aa}{{\mathfrak a}}

\newcommand{\bb}{{\mathfrak b}}

\newcommand{\cc}{{\mathfrak c}}
\newcommand{\dd}{{\mathfrak d}}

\newcommand{\ee}{{\mathfrak e}}
\renewcommand{\gg}{{\mathfrak g}}

\newcommand{\hh}{{\mathfrak h}}

\newcommand{\ii}{{\mathfrak i}}

\newcommand{\mm}{{\mathfrak m}}

\newcommand{\pp}{{\mathfrak p}}
\newcommand{\rr}{{\mathfrak r}}

\renewcommand{\ss}{{\mathfrak s}}

\renewcommand{\tt}{{\mathfrak t}}
\newcommand{\uu}{{\mathfrak u}}

\newcommand{\add}{{\mathsf{add}}}
\newcommand{\cov}{{\mathsf{cov}}}
\newcommand{\non}{{\mathsf{non}}}
\newcommand{\cof}{{\mathsf{cof}}}
\newcommand{\Add}[1]{{{\mathsf{add}}({\cal #1})}}
\newcommand{\Cov}[1]{{{\mathsf{cov}}({\cal #1})}}
\newcommand{\Non}[1]{{{\mathsf{non}}({\cal #1})}}
\newcommand{\Cof}[1]{{{\mathsf{cof}}({\cal #1})}}

\newcommand{\Fn}{{\mathsf{Fn}}}

\newcommand{\UBD}{{\mathsf{EUB}}}
\newcommand{\COF}{{\mathsf{COB}}}
\newcommand{\COB}{{\mathsf{COB}}}
\newcommand{\EUB}{{\mathsf{EUB}}}

\newcommand{\ran}{{\mathrm{ran}}}

\newcommand{\dom}{{\mathrm{dom}}}

\newcommand{\supp}{{\mathrm{supp}}}
\newcommand{\Dp}{{\mathrm{Dp}}}
\newcommand{\On}{{\mathrm{On}}}

\newcommand{\dir}{{\mathrm{dir}}}

\newcommand{\Sym}{{\mathrm{Sym}}}

\newcommand{\Spec}{{\mathrm{Spec}}}

\newcommand{\cf}{{\mathrm{cf}}}

\newcommand{\Ult}{{\mathrm{Ult}}}
\newcommand{\comp}{{\mathsf{comp}}}
\newcommand{\abs}{{\mathrm{abs}}}

\newcommand{\zero}{{\mathbf{0}}}
\newcommand{\one}{{\mathbf{1}}}
\newcommand{\PI}{{\mathbf{\Pi}}}
\newcommand{\SIGMA}{\mathbf{\Sigma}}

\newcommand{\noi}{\noindent}
\newcommand{\sub}{\subseteq}
\newcommand{\sem}{\setminus}
\newcommand{\twoom}{2^\omega}

\newcommand{\omlom}{\omega^{<\omega}}
\newcommand{\omom}{\omega^\omega}
\newcommand{\omloms}{[\omega]^{<\omega}}

\newcommand{\omoms}{[\omega]^\omega}

\newcommand{\ha}{\,{}\hat{}\,}
\newcommand{\la}{\langle}
\newcommand{\ra}{\rangle}

\newcommand{\Loleriar}{\mbox{$\Longleftrightarrow$}}

\title{Modern Forcing Techniques related to Finite Support Iteration: \\ Ultrapowers, templates, and submodels}

\author{J\"org Brendle\thanks{Partially supported by Grant-in-Aid for Scientific Research
   (C) 18K03398, Japan Society for the Promotion of Science.}  \\
   Graduate School of System Informatics \\
   Kobe University \\
   Rokko-dai 1-1, Nada-ku \\
   Kobe 657-8501, Japan \\
   email: {\sf brendle@kobe-u.ac.jp}}

\begin{document}
\maketitle

\begin{abstract}
\noindent This is an expository paper about several sophisticated forcing techniques closely related
to standard finite support iterations of ccc partial orders. We focus on the following four topics:
ultrapowers of forcing notions, iterations along templates, Boolean ultrapowers of forcing notions,
and restrictions of forcing notions to elementary submodels.
\end{abstract}

%%%%%%%%%%%%%%%

\section*{Introduction}

The method of finite support iteration (fsi) of ccc forcing, originally developed by Solovay and Tennenbaum to prove the
consistency of Suslin's hypothesis~\cite{ST71}, has since been used for a plethora of independence proofs,
both in set theory proper and in other areas of mathematics like topology or algebra. One recurring theme has
been its use for independence proofs about cardinal invariants of the continuum, that is, cardinal numbers
describing the combinatorial structure of the Baire space $\omom$ or Cantor space $\twoom$ and typically taking
values between the first uncountable cardinal $\aleph_1$ and the size $\cc$ of the continuum. Examples of such
cardinal invariants are the {\em unbounding number} $\bb$ and the {\em dominating number} $\dd$, the
least size of a family $F$ of functions in $\omom$ such that no single function dominates modulo finite all functions in $F$
(such that all functions are dominated modulo finite by a member of $F$, respectively).
The fact that the continuum can be made arbitrarily large adds to the versatility of the fsi method. 

However, there are situations when an fsi cannot be used (or when it is not known whether it can be used).
In such cases, a countable support iteration of proper forcing~\cite{Sh98} may be appropriate. This method,
though, makes the continuum have size at most $\aleph_2$ and therefore is of no use for distinguishing three or
more cardinal invariants. Therefore, a number of intricate and sophisticated methods, which are to some extent
modifications of fsi, have been developed for solving specific problems, and later used for further  results
about cardinal invariants. The purpose of the present survey paper is to introduce four such methods with the hope
of making them more accessible to researchers in the field. Specifically, we shall discuss
\begin{itemize}
\item ultrapowers of partial orders (Section~\ref{ultrapowers}),
\item iterations along templates (Section~\ref{templates}),
\item Boolean ultrapowers of partial orders (Section~\ref{boolultrapowers}), and
\item restrictions of partial orders to elementary submodels (Section~\ref{submodels}).
\end{itemize}
The first two methods were introduced by Shelah~\cite{Sh04} to prove the consistency of $\dd < \aa$,
first using a measurable cardinal and then from Con(ZFC). Here $\aa$ is the {\em almost disjointness number}, that is,
the least size of an infinite maximal almost disjoint (mad) family of infinite sets of natural numbers, one of the
most important cardinal invariants of the continuum because of its numerous applications in general topology. The third was introduced in two papers,
one, still unpublished, by Raghavan and Shelah~\cite{RSta} dealing with the consistency of $\dd (\lambda) < \aa (\lambda)$
where $\lambda$ is regular uncountable, and another, by Goldstern, Kellner, and Shelah~\cite{GKS19}
showing the consistency of Cicho\'n's maximum, the statement that all cardinal invariants in Cicho\'n's diagram
can be simultaneously distinct. Both use large cardinal assumptions. The fourth, then, was first used
by Goldstern, Kellner, Mej\'ia, and Shelah~\cite{GKMSta1} to prove the consistency of Cicho\'n's maximum from Con(ZFC).

For each of the four topics, we first describe the technique and prove a number of basic results in one
or two subsections. In the next subsection, we present one proof obtained by the respective method in full detail
and in the final subsection, we provide an overview of results obtained by the same method.

It should be noted that while only the template technique inherently is an fsi-style method, apart from the mentioned work by
Raghavan and Shelah, the others so far have been used in an fsi context. For ultrapowers this seems to
have to do with the fact that the rather complicated limit construction in the iteration (a limit strictly larger than the
direct limit is taken but the ccc is preserved, see Lemmata~\ref{iu-limit} and~\ref{iu-ccc}) only works in this case,
but for the other two methods, Boolean ultrapowers and submodels, we expect more applications to higher
cardinal invariants and, thus, e.g., to forcing notions which are $< \lambda$-closed and $\lambda^+$-cc for some
uncountable regular $\lambda$.

%

%

%

%%%

\subsection{Preliminaries}
\label{prelim}

We assume basic knowledge about forcing theory (see~\cite{Je03} and~\cite{Ku13}), as well as some practice with 
cardinal invariants and their interplay with forcing (see~\cite{BJ95}, \cite{Bl10}, and~\cite{Ha17}).

We often freely switch between the partial order (p.o.) language and the complete Boolean algebra (cBa) language
when dealing with forcing. We use $\PP \embed \QQ$ to denote that a p.o. $\PP$ completely embeds into a p.o. $\QQ$, and
$\PP \star \dot \QQ$ is the two-step iteration of $\PP$ with $\dot \QQ$. ``$\lim \dir$" denotes the direct limit of a
directed system of forcing notions. 

Let $\PP$ and $\QQ$ be complete Boolean algebras (cBas) with $\PP \sub \QQ$. Then the
{\em projection mapping} $h^\QQ_\PP : \QQ \to \PP$ is defined by $h^\QQ_\PP (q) = \prod \{ p \in \PP : q \leq_\QQ p \}$ for $q \in \QQ$.
Notice that in this context, $\PP \embed \QQ$ is equivalent to saying that for all  $q \in \QQ$ with $ q > \zero$, we also have
$h^\QQ_\PP (q) > \zero$.

A forcing notion $( \PP , \leq )$ is said to be {\em Suslin ccc} if $\PP$ is ccc and $\PP \sub \omom$,
$\leq \; \sub (\omom)^2$, as well as $\bot \sub (\omom)^2$ are $\SIGMA^1_1$ sets.  See~\cite[Section 3.6]{BJ95} for 
basic properties about Suslin ccc forcing. The following two basic lemmas are crucial and we therefore include the short proofs.

\begin{lem}   \label{Suslinccc-absoluteness}
If $\PP$ is Suslin ccc and $A \sub \PP$ is countable, ``$A $ is a maximal antichain" is a $\PI^1_1$ statement
and therefore absolute between models of ZFC. 
\end{lem}

\begin{proof}
Let $A = \{ x_n : n \in \omega \} \sub \PP$. ``$A$ is a maximal antichain" iff
\begin{itemize}
\item %$\forall n \neq m \; (x_n \bot x_m)$ iff 
   $\forall n \neq m \; \forall y \; ( y \notin \PP \lor y \not\leq x_n \lor y \not\leq x_m )$
\item $\forall y \; ( y \notin \PP \lor \exists n \; \neg (y \bot x_n) )$
\end{itemize}
Both formulas are $\PI^1_1$, and therefore $\SIGMA^1_1$-absoluteness applies.
\end{proof}

\begin{lem}   \label{Suslinccc-completeness}
Assume $\PP , \PP'$ are partial orders with $\PP \embed \PP '$. Also let $\QQ$ be a Suslin ccc forcing.
Then $\PP \star \dot \QQ \embed \PP ' \star \dot \QQ$ where the first $\dot \QQ$ is the $\PP$-name for 
(the interpretation of the code of) $\QQ$ in $V^{\PP}$ and the second the corresponding $\PP '$-name.
\end{lem}

Note here that $\dot \QQ$ is still (forced to be) ccc in $\PP$- and $\PP'$-generic extensions~\cite[Theorem 3.6.6]{BJ95}.

\begin{proof}
Assume $B = \{ (p_\alpha, \dot q_\alpha) : \alpha < \kappa \}$ is a maximal antichain in $\PP \star \dot \QQ$. We need to show
$B$ is still maximal in $\PP ' \star \dot \QQ$. 

Assume $(p' , \dot q') \in \PP ' \star \dot \QQ$. Let $G '$ be a $\PP '$-generic filter over $V$ with $p' \in G'$.
Then $G: = G' \cap \PP$ is $\PP$-generic over $V$ and, in $V[G]$, $\{ \dot q_\alpha [G] : p_\alpha \in G \}$
is a maximal antichain in $\dot \QQ [G]$. (This means in particular that $\{ \alpha : p_\alpha \in G \}$
is at most countable.) By the previous lemma, this is still a maximal antichain in $\dot \QQ [G ']$ in $V[G ']$.
Therefore, in $V [G']$, there is $\alpha$ with $p_\alpha \in G$ such that $\dot q_\alpha [G]$ and $\dot q ' [G']$ are compatible.
Hence there are $p'' \leq p_\alpha, p '$ in $\PP '$ and a $\PP '$-name $\dot q ''$ for an element of $\dot \QQ$ such that
$(p'' , \dot q '') \leq (p_\alpha, \dot q_\alpha), (p', \dot q ')$, as required.
\end{proof}

In Sections~\ref{ultrapowers} and~\ref{templates}, we need the following basic notion from~\cite{Br05} (see also~\cite[Definition 1]{Brta}). 
Let $\PP_{0 \land 1} \embed \PP_j \embed \PP_{0 \lor 1}$, $j \in \{ 0,1\}$, be cBa's. We say {\em projections in the diagram 
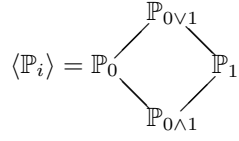
\[ 
\begin{picture}(80,50)(0,0)
\put(48,5){\makebox(0,0){$\PP_{0\land 1}$}}
\put(7,25){\makebox(0,0){$\la \PP_i \ra  = \PP_{0}$}}
\put(68,25){\makebox(0,0){$\PP_{1}$}}
\put(48,45){\makebox(0,0){$\PP_{0\lor 1}$}}
\put(26,28){\line(1,1){13}}
\put(28,20){\line(1,-1){11}}
\put(51,8){\line(1,1){13}}
\put(53,40){\line(1,-1){11}}
\end{picture}
\]
are correct} if either of the following three equivalent conditions holds:
\begin{itemize}
\item $h^{0 \lor 1}_1 (p_0) = h^0_{0 \land 1} (p_0)$ for all $p_0 \in \PP_0$,
\item $h^{0 \lor 1}_0 (p_1) = h^1_{0 \land 1} (p_1)$ for all $p_1 \in \PP_1$,
\item whenever $h^0_{0 \land 1} (p_0) = h^1_{0 \land 1} (p_1)$ then $p_0$ and $p_1$
are compatible in $\PP_{0 \lor 1}$. 
\end{itemize}
Notice this implies (but is not equivalent to)  $\PP_{0\land 1} = \PP_0 \cap \PP_1$. 
A typical example for a diagram with correct projections is given by letting $\PP_{0\land 1} = \{ \zero, \one \}$ and
$\PP_{0 \lor 1}$ the usual product forcing, that is, the completion of $(\PP_0 \sem \{ \zero \} ) \times (\PP_1 \sem \{ \zero \})$.
Another important example is obtained by letting $\PP_{0 \land 1} \embed \PP_0$ be arbitrary forcing notions and putting
$\PP_1 := \PP_{0\land 1} \star \dot \QQ$ and $\PP_{0 \lor 1} := \PP_0 \star \dot \QQ$, where $\QQ$ is a Suslin ccc forcing notion.
In both cases, correctness is straightforward. For an example of a non-correct diagram with $\PP_{0\land 1} = \PP_0 \cap \PP_1$ 
see~\cite[Counterexample 4]{Brta}.  More on correctness can be found in Section 1 of the latter work. We mention the following because
we will use it later in Lemma~\ref{los-forcing}.

\begin{lem}[{see~\cite[Observation 1 (ii)]{Brta}}]   \label{correctness-quotient}
Projections in the diagram 
\[\begin{picture}(80,50)(0,0)
\put(48,5){\makebox(0,0){$\PP_{0\land 1}$}}
\put(7,25){\makebox(0,0){$\la \PP_i \ra  = \PP_{0}$}}
\put(68,25){\makebox(0,0){$\PP_{1}$}}
\put(48,45){\makebox(0,0){$\PP_{0\lor 1}$}}
\put(26,28){\line(1,1){13}}
\put(28,20){\line(1,-1){11}}
\put(51,8){\line(1,1){13}}
\put(53,40){\line(1,-1){11}}
\end{picture}\]
are correct iff projections in the diagram 
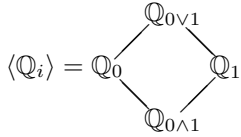
\[\begin{picture}(80,50)(0,0)
\put(48,5){\makebox(0,0){$\QQ_{0\land 1}$}}
\put(7,25){\makebox(0,0){$\la \QQ_i \ra  = \QQ_{0}$}}
\put(68,25){\makebox(0,0){$\QQ_{1}$}}
\put(48,45){\makebox(0,0){$\QQ_{0\lor 1}$}}
\put(26,28){\line(1,1){13}}
\put(28,20){\line(1,-1){11}}
\put(51,8){\line(1,1){13}}
\put(53,40){\line(1,-1){11}}
\end{picture}\]
are correct for every $\PP_{0 \land 1}$-generic filter $G_{0 \land 1}$ where $\QQ_i = \PP_i / G_{0 \land 1}$ (so $\QQ_{0 \land 1}$ is the trivial forcing).
\end{lem}

\begin{proof}
If $p = h^0_{0 \land 1} (p_0) = h^1_{0 \land 1} (p_1)$ and $p_0 \in \PP_0$ and $p_1 \in \PP_1$ are incompatible, take a generic $G_{0 \land 1}$ containing $p$ to see that
$p_0$ and $p_1$ are still incompatible in $\QQ_{0 \lor 1}$. On the other hand, if there are a generic $G_{0 \land 1}$ and $q_0 \in \QQ_0$ and $q_1 \in \QQ_1$ incompatible
in $\QQ_{0 \lor 1}$, then this incompatibility is forced by some $p \in G_{0 \land 1}$. Letting $p_i = (p, \dot q_i)$, we see that $h^i_{0 \land 1} (p_i) =p $, for $i = 0,1$,
and $p_0$ and $p_1$ are incompatible in $\PP_{0 \lor 1}$.
\end{proof}

Correctness can be used to show complete embeddability between direct limits.

\begin{lem}   \label{correctness-limdir}
Let $K$ be a directed index set. Assume $( \PP_k : k \in K )$ and $( \QQ_k : k \in K )$ are systems of cBa's such that $\PP_k \embed \PP_\ell$,
$\QQ_k \embed \QQ_\ell$, and $\PP_k \embed \QQ_k$ for any $k \leq \ell$. Assume further projections in all diagrams of the form
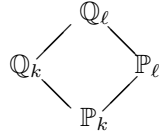
\[ 
\begin{picture}(80,50)(0,0)
\put(48,5){\makebox(0,0){$\PP_k$}}
\put(22,25){\makebox(0,0){$\QQ_k$}}
\put(68,25){\makebox(0,0){$\PP_\ell$}}
\put(48,45){\makebox(0,0){$\QQ_\ell$}}
\put(26,28){\line(1,1){13}}
\put(28,20){\line(1,-1){11}}
\put(51,8){\line(1,1){13}}
\put(53,40){\line(1,-1){11}}
\end{picture}
\]
are correct for $k \leq \ell$. Then $\PP : = \lim \dir_{k \in K} \PP_k$ completely embeds into $\QQ := \lim \dir_{k \in K} \QQ_k$. 
Furthermore, correctness is preserved in the sense that projections in all diagrams of the form
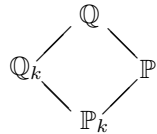
\[ 
\begin{picture}(80,50)(0,0)
\put(48,5){\makebox(0,0){$\PP_k$}}
\put(22,25){\makebox(0,0){$\QQ_k$}}
\put(68,25){\makebox(0,0){$\PP$}}
\put(46,45){\makebox(0,0){$\QQ$}}
\put(26,28){\line(1,1){13}}
\put(28,20){\line(1,-1){11}}
\put(51,8){\line(1,1){13}}
\put(53,40){\line(1,-1){11}}
\end{picture}
\]
for $k \in K$ are correct.
\end{lem}

\begin{proof}
Let $A \sub \bigcup_{k \in K} \PP_k$ be a maximal antichain in $\PP$. We have to show $A$ is still maximal in $\QQ$.
Choose $q \in \QQ$. Then $q \in \QQ_k$ for some $k \in K$. By maximality of $A$ there is $p \in A$ such that
$h^{\QQ_k}_{\PP_k} (q) $ is compatible with $p$ in $\PP$ and thus with $h^{\PP}_{ \PP_k} (p) $ in $\PP_k$. Find $\ell \geq k$
such that $p \in \PP_\ell$. By correctness, $p$ and $q$ are compatible in $\QQ_\ell$ and thus in $\QQ$, as required.
Preservation of correctness is straightforward.
\end{proof}

%%%%%%%%%%%%%%%%%%%%%%%

%

%

%

%

%

%

%%%%%%%%%%%%%%%

\section{Ultrapowers}
\label{ultrapowers}

Assume $\PP$ is a ccc partial order and $\kappa$ is a measurable cardinal as witnessed by the $\kappa$-complete ultrafilter $\D$. 
Then the ultrapower  $\PP^\kappa / \D$ is again a ccc partial order and $\PP$ completely embeds into $\PP^\kappa / \D$ so that we may
view $\PP^\kappa / \D$ as a two-step iteration of $\PP$ and some remainder forcing (see Subsection~\ref{ultrapowers-po} for details). 
$\PP^\kappa / \D$ shares many properties with
$\PP$ and some objects added by $\PP$ will actually be preserved by the ultrapower, while on the other hand, if $\PP$ forces that $\aa$ has size at least 
$\kappa$, then $\PP^\kappa / \D$ destroys all mad families of the intermediate extension $V^\PP$. This simple and ingenious observation,
due to Shelah, forms the basis of his consistency proofs of $\dd < \aa$ and $\uu < \aa$~\cite{Sh04}. 

Typically, the ultrapower operation is applied to iterations $\P = \la \PP_\gamma : \gamma \leq \mu \ra$, and ultrapowers of iterations
are again iterations. When iterating the process of taking such ultrapowers the question arises what to do in limit steps.
One option is a direct limit (this has been used e.g. in Theorem~\ref{b<as}), but often embedding the iterations obtained by
taking ultrapowers into a larger iteration is necessary (e.g. for $\dd < \aa$ and $\uu < \aa$). Technical details of this are discussed in Subsection~\ref{ultrapowers-iterations}.

In Subsection~\ref{a>d-measurable} we present a complete proof of Shelah's consistency of $\dd < \aa$ from a measurable (Theorem~\ref{a>d-meas}),
and in Subsection~\ref{ultrapowers-further} we discuss further results obtained by the ultrapower method. Our exposition follows to some extent our earlier~\cite{Br07}.

%\bigskip

\subsection{Ultrapowers of partial orders}
\label{ultrapowers-po}

Let $\kappa$ be a measurable cardinal and let $\D$ be a $\kappa$-complete ultrafilter on $\kappa$.
For a p.o. $\PP$ and $f \in \PP^\kappa$, 
\[ [f] = f / \D = \{ g \in \PP^\kappa : \{ \alpha < \kappa : f(\alpha) = g(\alpha) \} \in \D \} \]
is the equivalence class of $f$ modulo $\D$. The {\em ultrapower} $\PP^\kappa / \D$ of $\PP$ consists of all such equivalence classes.
It is partially ordered by
\[ [f] \leq [g]\; \mbox{ iff } \; \{ \alpha < \kappa : f(\alpha) \leq g(\alpha) \} \in \D \]
As usual, we identify $p \in \PP$ with the class $[f]$ of the constant function $f (\alpha) = p$, $\alpha < \kappa$, and
thus construe $\PP$ as a subset of $\PP^\kappa /  \D$.

If we consider the elementary embedding $j_\D : V \to M$ derived from $\D$, $j_\D (\PP) = \PP^\kappa / \D$ is a p.o. in $M \sub V$,
and results like the next two lemmata have alternative shorter proofs as absoluteness arguments, based on the fact that $M$ is $\kappa$-closed.
However, when considering iterations of length some $\mu$ later on, we want to think of their ultrapowers as being of length $\mu$ and not of length $j_\D (\mu)$
(though they technically are), and when iterating the process of taking ultrapowers we will always use the same ultrafilter $\D$ and not its iterates.
Since these issues may make the $j_\D (\PP)$ notation confusing we prefer the $\PP^\kappa / \D$ approach.

\begin{lem} \label{ultra-embed}
If $\PP$ is $\kappa$-cc then $\PP \embed \PP^\kappa / \D$.
\end{lem}

\begin{proof}
Fix $\nu < \kappa$, and let $A = \{ p_\gamma : \gamma < \nu \}$ be a maximal antichain in $\PP$. Given arbitrary $f \in \PP^\kappa$,
for all $\alpha < \kappa$ there is $\gamma < \nu$ such that $f(\alpha)$ and $p_\gamma$ are compatible. By $\kappa$-completeness
of $\D$, there is $\gamma < \nu$ such that $\{ \alpha < \kappa : f(\alpha)$ and $p_\gamma$ are compatible$\}$ belongs to $\D$.
Hence $[f]$ is compatible with $p_\gamma$ in $\PP^\kappa / \D$, and $A$ is still a maximal antichain in $\PP^\kappa / \D$.
\end{proof}

Notice that the converse holds as well. If $\PP$ is not $\kappa$-cc, then there is a maximal antichain $\{ p_\gamma : \gamma < \mu \} \sub \PP$
with $\mu \geq \kappa$, and $[f]$ given by $f(\alpha) = p_\alpha$ for $\alpha < \kappa$ is an element of $\PP^\kappa / \D$ incompatible with
all $p_\gamma$. Hence $\PP$ does not completely embed into $\PP^\kappa/ \D$.

\begin{lem} \label{ultra-chain}
If $\PP$ is $\nu$-cc for some $\nu < \kappa$ then so is $\PP^\kappa / \D$.
\end{lem}

\begin{proof}
Take arbitrary elements $f_\gamma \in \PP^\kappa$, $\gamma < \nu$. By the $\nu$-cc of $\PP$, for each $\alpha < \kappa$, there are $\gamma < \delta 
< \nu$ such that $f_\gamma(\alpha)$ and $f_\delta (\alpha)$ are compatible. By $\kappa$-completeness of $\D$, there are $\gamma < \delta 
< \nu$ such that $\{ \alpha < \kappa : f_\gamma (\alpha)$ and $f_\delta (\alpha)$ are compatible$\}$ belongs to $\D$. Thus $[f_\gamma]$ and
$[f_\delta]$ are compatible. Therefore every antichain of $\PP^\kappa / \D$ has size less than $\nu$.

(The absoluteness argument for this proof runs as follows: Since $V \models ``\PP$ is $\nu$-cc", $M \models `` j_\D (\PP)$ is $j_\D (\nu)$-cc". But $j_\D (\nu ) = \nu$,
and if $A \in V$ were an antichain in $j_D (\PP)$ of size $\nu$, $A$ would belong to $M$ by $\kappa$-closure of $M$, a contradiction. Hence $V  \models ``j_\D (\PP)$ is $\nu$-cc".)
\end{proof}

On the other hand, if $\PP$ is not $\nu$-cc for any $\nu < \kappa$, then $\PP^\kappa / \D$ is not $\kappa$-cc. Indeed, let
$A_\alpha = \{ p_{\alpha,\gamma} : \gamma < \nu_\alpha \}$ be a maximal antichain in $\PP$ of size $\nu_\alpha \geq |\alpha|$ for each $\alpha < \kappa$, and fix $p \in \PP$ arbitrarily.
Define $f_\gamma \in \PP^\kappa$ for $\gamma < \kappa$ by \[f_\gamma (\alpha ) = \left\{ \begin{array}{ll} p_{\alpha,\gamma} & \mbox{ if } \nu_\alpha > \gamma \\ 
p & \mbox{ otherwise} \end{array} \right. \]
It is easy to see that $\{ [f_\gamma] : \gamma < \kappa \}$ is an antichain in $\PP^\kappa / \D$.
For more general versions of the two lemmata, see Lemmas~\ref{Boolultra-embed} and~\ref{Boolultra-chain} in Subsection~\ref{boolultrapowers-po}.

For the remainder of this section, we assume $\PP$ is ccc. Therefore $\PP \embed \PP^\kappa /\D$ and $\PP^\kappa / \D$ is ccc by the two previous lemmata.

We next describe the relationship between $\PP$-names and $\PP^\kappa / \D$-names for real numbers. First notice that given $\kappa$ many maximal antichains
$\{ p_n^\alpha : n \in \omega \}$, $\alpha < \kappa$, in $\PP$, letting $f_n (\alpha) = p_n^\alpha$ for $\alpha < \kappa$ we obtain a maximal antichain 
$\{ [f_n] : n \in \omega \}$ in $\PP^\kappa /\D$. Furthermore all maximal antichains of $\PP^\kappa / \D$ are of this form. Now recall that a $\PP$-name $\dot x$ for
a real in $\omom$ is given by maximal antichains $\{ p_{n,i} : n \in \omega \}$ and numbers $\{ k_{n,i} : n \in \omega \}$, $i \in \omega$, such that
\[ p_{n,i} \forces_\PP \dot x (i) = k_{n,i} . \]
Therefore, a $\PP^\kappa / \D$-name $\dot y$ for a real is given by maximal antichains $\{ p_{n,i}^\alpha : n \in \omega \}$ and numbers $\{ k_{n,i} : n \in \omega \}$, 
$i \in \omega$ and $\alpha < \kappa$, such that, letting $f_{n,i} (\alpha) = p_{n,i}^\alpha$,
\[ [f_{n,i}] \forces_{\PP^\kappa / \D}  \dot y (i) = k_{n,i} . \]
Since $\{ p_{n,i}^\alpha : n \in \omega \}$ and $\{ k_{n,i} : n \in \omega \}$, $i \in \omega$, determine a $\PP$-name $\dot x^\alpha$ for a real, 
we may think of $\dot y$ as the {\em average} or {\em mean} of the $\dot x^\alpha$ and write $\dot y = \la \dot x^\alpha : \alpha < \kappa \ra / \D$.
Notice that every $\PP^\kappa / \D$-name for a real is of this form.

\begin{lem} \label{ultra-mad}
Let $\PP$ be ccc. Assume $\dot \A$ is a $\PP$-name for a mad family of size at least $\kappa$. Then $\PP^\kappa /\D$ forces that $\dot \A$ is not maximal.
In particular, if $\PP$ forces $\aa \geq \kappa$, then no a.d. family of $V^\PP$ is maximal in $V^{\PP^\kappa / \D}$.
\end{lem}

\begin{proof}
Assume $\dot \A = \{ \dot A^\gamma : \gamma < \nu \}$ where $\nu \geq \kappa$ and all $\dot A^\gamma$ are $\PP$-names for infinite subsets of $\omega$.
Then $\dot A = \la \dot A^\alpha : \alpha < \kappa \ra / \D$ is a $\PP^\kappa / \D$-name for an infinite subset of $\omega$ by the preceding discussion.

Fix $\gamma < \nu$. Since for all $\alpha < \kappa$ with $\alpha \neq \gamma$,
\[ \forces_\PP | \dot A^\alpha \cap \dot A^\gamma | < \aleph_0, \]
we see that $\{ \alpha < \kappa : \forces_\PP | \dot A^\alpha \cap \dot A^\gamma | < \aleph_0 \}$ belongs to $\D$. Therefore
\[ \forces_{\PP^\kappa / \D}  | \dot A \cap \dot A^\gamma | < \aleph_0 \]
because $\dot A$ is the average of the $\dot A^\alpha$ (this is a straightforward elementarity argument; for a formal proof, see Lemma~\ref{los-forcing}
in Subsection~\ref{boolultrapowers-po} below). Hence $\dot A$ witnesses that $\PP^\kappa / \D$ forces non-maximality of $\dot \A$.
\end{proof}

Let $\mu$ be an uncountable regular cardinal. Say that a sequence $\{ f_\beta : \beta < \mu \} \sub \omom$ is a {\em scale } if $f_\alpha \leq^* f_\beta$
for $\alpha < \beta < \mu$ and for all $g \in \omom$ there is $\alpha < \mu$ with $g \leq^* f_\alpha$. It is well-known and easy to see that
the existence of a $\mu$-scale is equivalent to $\bb = \dd = \mu$.

\begin{lem} \label{ultra-scales}
Let $\PP$ be ccc, $\mu \neq \kappa$ regular, and assume $\PP$ adjoins a scale $\{ \dot d_\beta : \beta < \mu \}$. Then $\PP^\kappa /\D$ forces that
$\{ \dot d_\beta : \beta < \mu \}$ is still a scale.  In particular, if $\PP^\kappa / \D \cong \PP \star \dot \QQ$, then $\PP$ forces that $\dot \QQ$ is an 
$\omom$-bounding forcing notion.
\end{lem}

\begin{proof}
Let $\dot y$ be a $\PP^\kappa / \D$-name for a real in $\omom$. By the preceding discussion, there are $f_{n,i} \in \PP^\kappa$ and $k_{n,i} \in \omega$,
$n,i \in \omega$, such that $\dot y$ is determined by $\{ [f_{n,i}] : n \in \omega \}$ and $\{ k_{n,i} : n \in \omega \}$, $i \in \omega$.
Letting $p_{n,i}^\alpha = f_{n,i} (\alpha)$, $\{ p_{n,i}^\alpha : n \in \omega \}$ and $\{ k_{n,i} : n \in \omega \}$, $i \in \omega$, determine $\PP$-names $\dot x^\alpha$ for reals,
and $\dot y = \la \dot x^\alpha : \alpha < \kappa \ra / \D$. 

By ccc-ness of $\PP$, for each $\alpha < \kappa$ there is $\beta_\alpha < \mu$ such that \[ \forces_\PP \dot x^\alpha \leq^* \dot d_{\beta_\alpha}. \]
Since $\mu \neq \kappa$ are both regular, we obtain $\beta < \mu$ such that 
$\{ \alpha < \kappa : \beta_\alpha \leq \beta \} \in \D$ (use the $\kappa$-completeness of $\D$ in case $\mu < \kappa$). Therefore
\[ \forces_{\PP^\kappa / \D}   \dot y  \leq^* \dot d_\beta  \]
because $\dot y$ is the average of the $\dot x^\alpha$. Hence $\{ \dot d_\beta : \beta < \mu \}$ remains a scale in the $\PP^\kappa / \D$-extension.
\end{proof}

\underline{\sf STRATEGY.} These two simple lemmas provide us with a scenario for proving the consistency of $\dd < \aa$. Let $\kappa < \mu < \lambda$ be regular
cardinals. Force $\bb = \dd = \mu$ with a ccc p.o. $\PP$. Then keep taking ultrapowers of $\PP$ for $\lambda$ many steps. By Lemma~\ref{ultra-scales}, 
$\bb = \dd =\mu$ should be preserved while, by the ZFC-inequality $\bb \leq \aa$, Lemma~\ref{ultra-mad} tells us that mad families of size less than $\lambda$
will be destroyed so that $\aa = \lambda$ in the final model. The problem, however, is what to do in limit steps of the procedure of iteratively taking
ultrapowers. To get a handle on this, we shall look at ultrapowers of iterations and iterations of ultrapowers in the next subsection.

%\bigskip

%%%

\subsection{Ultrapowers and iterations}
\label{ultrapowers-iterations}

We start with the discussion of two-step iterations. 

\begin{lem} \label{ultra-twostep}
Assume $\PP \embed \QQ$ are cBa's. Then $\PP^\kappa / \D \embed \QQ^\kappa / \D$ and, in fact, the projection mapping is given by 
$h^{\QQ^\kappa / \D}_{\PP^\kappa / \D} ( [f]) = \la h^\QQ_\PP (f(\alpha)) : \alpha < \kappa \ra / \D$ for $f \in \QQ^\kappa$.
Furthermore, projections in the diagram $\la \PP, \QQ, \PP^\kappa / \D , \QQ^\kappa / \D \ra$ are correct.
\end{lem}

\begin{proof}
Let $g \in \PP^\kappa$ and notice that
\[ \begin{array}{rcl} [g] \geq h^{\QQ^\kappa / \D}_{\PP^\kappa / \D} ( [f])  & \Loleriar & [g] \geq [f] \\
& \Loleriar & \{\alpha < \kappa : g(\alpha) \geq f(\alpha) \} \in \D \\
& \Loleriar & \{\alpha < \kappa : g(\alpha) \geq h^\QQ_\PP(f(\alpha)) \} \in \D \\
& \Loleriar & [g] \geq \la h^\QQ_\PP (f(\alpha)) : \alpha < \kappa \ra / \D \end{array} \]
This shows equality. Thus
\[ \begin{array}{rcl} [g] \bot h^{\QQ^\kappa / \D}_{\PP^\kappa / \D} ( [f])  & \Loleriar & [g] \bot \la h^\QQ_\PP (f(\alpha)) : \alpha < \kappa \ra / \D \\
& \Loleriar & \{\alpha < \kappa : g(\alpha) \bot h^\QQ_\PP(f(\alpha)) \} \in \D \\
& \Loleriar & \{\alpha < \kappa : g(\alpha) \bot  f(\alpha) \} \in \D \\
& \Loleriar & [g] \bot [f] \end{array} \]
and we obtain complete embeddability. To see correctness, let $q \in \QQ$, $p: = h^\QQ_\PP (q)$, $g(\alpha) = q$,
and $f (\alpha ) = p$ for all $\alpha$. Then
\[ h^{\QQ^\kappa / \D}_{\PP^\kappa / \D} (q) = h^{\QQ^\kappa / \D}_{\PP^\kappa / \D} ([g]) = \la h^\QQ_\PP (g(\alpha)) : \alpha < \kappa \ra / \D = 
\la f(\alpha) : \alpha < \kappa \ra / \D = [f] = p \]
as required.
\end{proof}

\begin{lem}  \label{ultra-Suslin}  
Let $\PP$ be a p.o. and let $\QQ$ be a Suslin ccc forcing notion. Then $(\PP \star \dot \QQ)^\kappa / \D \cong \PP^\kappa / \D \star \dot \QQ$.
\end{lem}

Note that the first $\dot \QQ$ is a $\PP$-name while the second is a $\PP^\kappa / \D$-name. Also, in elementary embedding notation,
the conclusion reads as $j_\D (\PP \star \dot \QQ) = j_\D (\PP) \star \dot \QQ$.

\begin{proof}
Taking $(p_\alpha ,  \dot x_\alpha) \in \PP \star \dot \QQ$,
$\alpha < \kappa$, we see that an arbitrary condition in $(\PP \star \dot \QQ)^\kappa / \D$ is of the form $\la (p_\alpha ,  \dot x_\alpha): \alpha < \kappa \ra / \D$.
Therefore, letting $f \in \PP^\kappa$ be defined by
$f(\alpha) = p_\alpha$ and setting $\dot y = \la \dot x_\alpha : \alpha < \kappa \ra / \D$, we see that we can identify $\la (p_\alpha , \dot x_\alpha): \alpha < \kappa \ra / \D$
with the condition $( [f] ,  \dot y) \in \PP^\kappa / \D \star \dot \DD$. This obviously defines an isomorphism between (dense subsets of) the two partial orders.
\end{proof}

We now move to transfinite iterations. Let $\mu$ be an ordinal. Say a sequence of cBas $\P = \la \PP_\gamma : \gamma \leq \mu \ra$ is an
{\em iteration} (see also~\cite{Br07}) if $\PP_\gamma \embed \PP_\delta$ for $\gamma < \delta$. Note that, for technical reasons which will become
obvious later on (see Lemma~\ref{ultra-iteration}), we do {\em not} require that $\PP_\delta$ is any kind of limit of $\PP_\gamma$, $\gamma < \delta$, for limit ordinals $\delta$.
For $\gamma < \delta$, let $h^\delta_\gamma : \PP_\delta \to \PP_\gamma$ be the projection. The {\em support} of $p \in \PP_\mu$ is defined by
\[ \supp (p) = \{ \delta : \mbox{ there is no } \gamma < \delta \mbox{ such that } h^\mu_\delta (p) = h^\mu_\gamma (p) \} \]
Note that $\delta + 1 \in \supp (p)$ iff $h^\mu_{\delta + 1} (p) < h^\mu_\delta (p)$. Similarly, for limit ordinals $\delta$,
$\delta \in \supp (p)$ iff $h^\mu_\delta (p) < h^\mu_\gamma (p)$ for all $\gamma < \delta$.

An iteration $\P$ has {\em finite supports} if there is a sequence $\la D_\gamma \sub \PP_\gamma : \gamma \leq \mu \ra$ of dense sets such that
\begin{itemize}
\item[(I)] $D_\gamma \sub D_\delta$ for $\gamma \leq \delta \leq \mu$,
\item[(II)] $h_\gamma^\mu (p) \in D_\gamma$ for $\gamma \leq \mu$ and $p \in D_\mu$, and
%\item[(III)] all $D_\gamma$ are closed under non-zero products (if $p,q \in D_\gamma$ and $p \cdot q \neq \zero$ then $p \cdot q \in D_\gamma$), and
\item[(III)] $\supp (p)$ is finite for all $p \in D_\mu$. 
\end{itemize}
While this is not the usual definition of a finite support iteration (fsi)
the following simple lemma implies that an iteration with finite supports is equivalent to an fsi (we leave the details of this to the reader).
Furthermore, our definition of the support is nonstandard as well: ``$\delta$ belongs to the support of $p$ in the traditional sense" is equivalent to $\delta + 1 \in \supp (p)$ in our sense,
and for limit ordinals $\delta$, $\delta \in \supp (p)$ for some $p$ means that  $\lim \dir_{\gamma < \delta} \PP_\gamma$ is a proper subforcing
of $\PP_\delta$.

\begin{lem}   \label{finite-supports}
Assume $\P = \la \PP_\gamma : \gamma \leq \mu \ra$ has finite supports. Let $\delta \leq \mu$ be a limit ordinal. Then $\lim \dir_{\gamma < \delta}
\PP_\gamma \embed \PP_\delta$ where, as usual, $\lim\dir$ denotes the direct limit of forcing notions.
\end{lem}

\begin{proof}
Let $p \in D_\delta$. If $\delta \notin \supp (p)$, then $p = h^\mu_\delta (p) = h^\mu_{\gamma_0} (p) $ for some $\gamma_0 < \delta$. Therefore
$p \in \PP_{\gamma_0} \sub \bigcup_{\gamma < \delta} \PP_\gamma = \lim \dir_{\gamma < \delta} \PP_\gamma$, and there is nothing to prove.

If $\delta \in \supp (p)$, then $p = h^\mu_\delta (p) < h^\mu_\gamma (p)$ for all $\gamma < \delta$. However, since supports are finite, there is
$\gamma_0 < \delta$ such that $h^\mu_{\gamma_0} (p) = h^\mu_\gamma (p)$ for all $\gamma$ with $\gamma_0 \leq \gamma < \delta$.
Set $p_0 = h^\mu_{\gamma_0} (p)$. So $p_0 \in \PP_{\gamma_0}  \sub \bigcup_{\gamma < \delta} \PP_\gamma = \lim \dir_{\gamma < \delta}
\PP_\gamma$. We claim $p_0$ is a reduction of $p$ to $ \lim \dir_{\gamma < \delta} \PP_\gamma$. Indeed, suppose $q \leq p_0$ belongs to
$ \lim \dir_{\gamma < \delta} \PP_\gamma$. Then there is $\gamma$ with $\gamma_0 \leq \gamma < \delta$ such that $q \in \PP_\gamma$. 
Since $h^\mu_\gamma (p) = p_0 \geq q$ in $\PP_\gamma$, $p$ and $q$ are compatible in $\PP_\delta$ (with common extension $p \cdot q$),
as required.
\end{proof}
 
We next investigate \underline{ultrapowers of iterations}.

\begin{lem}  \label{ultra-iteration}
Let $\P = \la \PP_\gamma : \gamma \leq \mu \ra$ be an iteration. Then $\P^\kappa / \D = \la (\PP_\gamma)^\kappa / \D : \gamma \leq \mu \ra$ also
is an iteration. Moreover, if $\P$ has finite supports then so does $\P^\kappa / \D$, as witnessed by $\la (D_\gamma)^\kappa / \D : \gamma \leq \mu \ra$.
\end{lem}

\begin{proof}
The first part is obvious by Lemma~\ref{ultra-twostep}. So assume that $\P$ has finite supports. Let $f \in (D_\mu)^\kappa$. Then $\supp ( f (\alpha))$
is finite for all $\alpha$, say $\supp (f (\alpha)) = \{ \gamma^\alpha_0 < \gamma^\alpha_1 < ... < \gamma^\alpha_{n_\alpha - 1} \}$
for all $\alpha$. By completeness of $\D$ there is $n \in \omega$ such that $\{ \alpha < \kappa : n_\alpha = n \}$ belongs to $\D$. 
For each $i < n$ there is $\gamma_i \leq \mu$ such that
\begin{itemize}
\item either $\{ \alpha < \kappa : \gamma_i^\alpha = \gamma_i \} \in \D$
\item or $\{ \alpha < \kappa : \gamma_i^\alpha < \gamma_i \} \in \D$ while $\{ \alpha < \kappa : \gamma_i^\alpha \leq \delta\} \notin \D$ for all $\delta < \gamma_i$.
\end{itemize}
In the latter case we necessarily have $\cf (\gamma_i) = \kappa$ by the $\kappa$-completeness of $\D$. 

We claim that $\supp ( [f]) = \{ \gamma_i : i < n \}$ and so is finite as required. (Therefore, $| \supp ([f]) | \leq n$, but equality does not necessarily hold because
the $\gamma_i$ are not necessarily distinct.) 

Indeed, 
\[ \begin{array}{rclr} \gamma \in \supp ( [f])  & \Loleriar &  \forall \delta < \gamma \; (    h^\mu_\gamma ([f]) < h^\mu_\delta ([f]) )  & \\
& \Loleriar &  \forall \delta < \gamma \; (  \{\alpha < \kappa : h^\mu_\gamma (f(\alpha)) < h^\mu_\delta( f(\alpha) )\} \in \D  )  &   \mbox{(by Lemma~\ref{ultra-twostep})} \\
& \Loleriar & \mbox{either } \{\alpha < \kappa : \gamma \in \supp (f (\alpha))  \} \in \D & \mbox{(first case)} \\
&& \mbox{or } \cf (\gamma) = \kappa \mbox{ and }  \forall \delta < \gamma \; (  \{\alpha < \kappa : (\delta,\gamma) \cap \supp (f(\alpha)) \neq \emptyset \} \in \D )  & \mbox{(second case)} \\
& \Loleriar & \gamma = \gamma_i  \mbox{ for some } i < n &  \end{array} \]
\end{proof}

We note that even if $\PP_\delta = \lim\dir_{\gamma < \delta} \PP_\gamma$ for all limits $\delta$, this is not necessarily the case
for $(\PP_\delta)^\kappa / \D$. Indeed, the proof above shows that for $\delta$ with $\cf (\delta) = \kappa$, $\lim\dir_{\gamma < \delta} (\PP_\gamma)^\kappa / \D$ 
will be a strict complete suborder of $(\PP_\delta)^\kappa / \D$. This is the reason for our notion of {\em iteration with finite supports}.

We want to iterate the procedure of taking ultrapowers -- and thus obtain \underline{iterations of ultrapowers}. 
To this end, we need to explain what to do in limit steps.
The following lemma (Lemma 7 of~\cite{Br07}) is a special case of the amalgamated limit from~\cite[Section 2]{Brta}.

\begin{lem}     \label{iu-limit}
Let $\mu$ and $\lambda$ be limit ordinals. Assume $\P^\alpha = \la \PP^\alpha_\gamma : \gamma\leq \mu \ra$, $\alpha < \lambda$, are iterations
such that $\PP^\alpha_\gamma \embed \PP^\beta_\gamma$ for $\alpha < \beta < \lambda$ and $\gamma \leq \mu$. Also assume $\la \PP^\lambda_\gamma :
\gamma < \mu \ra$ is an iteration such that $\PP^\alpha_\gamma \embed \PP^\lambda_\gamma$ for $\alpha < \lambda$ and $\gamma < \mu$. 
Furthermore, assume that for all $\alpha < \beta \leq \lambda$ and $\gamma < \delta \leq \mu$ with $(\beta,\delta) \neq (\lambda,\mu)$, the projections
in the diagram $\la \PP^\alpha_\gamma ,\PP^\beta_\gamma , \PP^\alpha_\delta , \PP^\beta_\delta \ra$ are correct.

Then there is a p.o. $\PP^\lambda_\mu$ such that $\P^\lambda = \la \PP^\lambda_\gamma : \gamma \leq \mu \ra$ is an iteration and $\PP^\alpha_\mu 
\embed \PP^\lambda_\mu$ for all $\alpha < \lambda$. Moreover, correctness is preserved in the sense that the projections in the
diagram $\la \PP^\alpha_\gamma ,\PP^\lambda_\gamma , \PP^\alpha_\mu , \PP^\lambda_\mu \ra$ are correct as well.
Furthermore, if for cofinally many $\alpha < \lambda$, $\PP^\alpha_\mu$ is the direct limit of $\la \PP^\alpha_\gamma : \gamma < \mu \ra$, then $\PP^\lambda_\mu$
is the direct limit of $\la \PP^\lambda_\gamma : \gamma < \mu \ra$, and, similarly, if for cofinally many $\gamma < \mu$, $\PP^\lambda_\gamma$ is the
direct limit of $\la \PP^\alpha_\gamma : \alpha < \lambda \ra$, then $\PP^\lambda_\mu$
is the direct limit of $\la \PP^\alpha_\mu : \alpha < \lambda \ra$.

Assume also all $\P^\alpha$, $\alpha < \lambda$, and $\la \PP^\lambda_\gamma : \gamma < \mu \ra$ have finite supports as witnessed
by $\la D^\alpha_\gamma : \gamma \leq \mu \ra$, $\alpha < \lambda$, and $\la D^\lambda_\gamma : \gamma < \mu \ra$ with
$D^\alpha_\gamma \sub D^\beta_\gamma$ for  $\alpha < \beta \leq \lambda$ and $\gamma \leq \mu$ and $(\beta,\gamma) \neq (\lambda,\mu)$,
and such that for all $\gamma < \delta < \mu$ and $\alpha < \lambda$, if $p \in D^\alpha_\delta$ and $q \in D^\lambda_\gamma$ are such that
$h^{\gamma,\lambda}_{\gamma,\alpha} (q) \leq h^{\delta,\alpha}_{\gamma,\alpha} (p)$, then $p \cdot q \in D^\lambda_\delta$.
Then there is $D^\lambda_\mu$ with $D^\alpha_\mu \sub D^\lambda_\mu$ for $\alpha < \lambda$ such that $\la D^\lambda_\gamma : \gamma \leq\mu\ra$ 
witnesses that $\P^\lambda$ has finite supports,
and such that for all $\gamma  < \mu$ and $\alpha < \lambda$, if $p \in D^\alpha_\mu$ and $q \in D^\lambda_\gamma$ are such that
$h^{\gamma,\lambda}_{\gamma,\alpha} (q) \leq h^{\mu,\alpha}_{\gamma,\alpha} (p)$, then $p \cdot q \in D^\lambda_\mu$.
\end{lem}

\begin{proof}  
$\PP^\lambda_\mu$ is the cBa generated by formal products of the form $p \cdot q$ where $p \in \PP^\alpha_\mu$ for some $\alpha < \lambda$,
$q \in \PP^\lambda_\gamma$ for some $\gamma < \mu$, and $h^{\gamma,\lambda}_{\gamma,\alpha} (q) = h^{\mu,\alpha}_{\gamma,\alpha} (p)$
(see Figure 1). This means that such formal products form a dense subset of $\PP^\lambda_\mu$. 
In this case we say that the pair $\alpha,\gamma$ witnesses $p \cdot q \in \PP^\lambda_\mu$. For formal products $p_0 \cdot q_0$ with 
witness $\alpha,\gamma$ and $p_1 \cdot q_1$ with witness $\beta,\delta$, where $\beta \geq \alpha$ and $\delta \geq \gamma$, we define
the partial order $\leq$ on $\PP^\lambda_\mu$ by $p_1 \cdot q_1 \leq p_0 \cdot q_0$ if $p_1 \leq p_0$ in $\PP^\beta_\mu$ and $q_1 \leq q_0$
in $\PP^\lambda_\delta$. 
\begin{gather*}
\begin{picture}(240,140)(0,0)
\put(20,20){\vector(1,0){220}}
\put(20,20){\vector(0,1){105}}
\put(205,20){\line(0,1){85}}
\multiput(120,20)(0,10){6}{\line(0,1){5}}
\multiput(120,90)(0,10){2}{\line(0,1){5}}
\put(20,105){\line(1,0){185}}
\multiput(20,80)(10,0){6}{\line(1,0){5}}
\multiput(160,80)(10,0){5}{\line(1,0){5}}
\put(205,10){\makebox(0,0){$\lambda$}}
\put(120,10){\makebox(0,0){$\alpha$}}
\put(120,115){\makebox(0,0){$p \in \PP^\alpha_\mu$}}
\put(10,105){\makebox(0,0){$\mu$}}
\put(10,80){\makebox(0,0){$\gamma$}}
\put(225,80){\makebox(0,0){$q\in \PP^\lambda_\gamma$}}
\put(220,115){\makebox(0,0){$p \cdot q \in \PP^\lambda_\mu$}}
\put(120,80){\makebox(0,0){$h^{\mu,\alpha}_{\gamma,\alpha} (p) = h^{\gamma,\lambda}_{\gamma,\alpha} (q)$}}
\end{picture} \\
\mbox{Figure 1}
\end{gather*}
We first verify that if a pair $\alpha, \gamma$ witnesses $p \cdot q \in \PP^\lambda_\mu$, then the pair $\beta, \delta$ is
also a witness for any $\beta \geq \alpha$ and $\delta \geq \gamma$. For indeed, by correctness, we have 
$h^{\gamma, \lambda}_{\gamma, \beta} (q) \leq h^{\gamma,\lambda}_{\gamma,\alpha} (q) = h^{\mu,\alpha}_{\gamma,\alpha} (p) = h^{\mu, \beta}_{\gamma, \beta} (p)$.
Thus, letting $p' = p \cdot h^{\gamma, \lambda}_{\gamma, \beta} (q)  \in \PP^\beta_\mu$, we see $h^{\mu,\beta}_{\gamma,\beta} (p') =
h^{\gamma, \lambda}_{\gamma, \beta} (q)$ so that the formal product $p' \cdot q$ belongs to $\PP^\lambda_\mu$ as witnessed by
the pair $\beta, \gamma$. Obviously, $p' \cdot q \leq p \cdot q$, and it is easy to see that any $p'' \cdot q'' \leq p \cdot q$ is compatible 
with $p' \cdot q$ so that $p' \cdot q$ and $p \cdot q$ in fact describe the same condition. By symmetry we can also replace $\gamma$ by $\delta$.

This fact provides us with an easier description of the ordering: given formal products $p_0 \cdot q_0$ and $p_1 \cdot q_1$ in $\PP^\lambda_\mu$, 
we may assume they have the same witness $\alpha, \gamma$, and we let $p_1 \cdot q_1 \leq p_0 \cdot q_0$ if $p_1 \leq p_0$ in $\PP^\alpha_\mu$ and
$q_1 \leq q_0$ in $\PP^\lambda_\gamma$.

Next we prove completeness of the embeddings. By symmetry it suffices to show that $\PP^\alpha_\mu \embed \PP^\lambda_\mu$ ($\PP^\lambda_\gamma
\embed \PP^\lambda_\mu$ is analogous). Fix $p \cdot q \in \PP^\lambda_\mu$ with witness $\beta, \gamma$. By the preceding discussion 
we may assume that $\beta \geq \alpha$. Let $p_0 = h^{\mu,\beta}_{\mu,\alpha} (p) \in \PP^\alpha_\mu$. We need to show that
any $p_1 \leq p_0$ in the partial order $\PP^\alpha_\mu$ is compatible with $p \cdot q$. Let $p' = p_1 \cdot p \in \PP^\beta_\mu$. Then
clearly $p' \leq p$ in $\PP^\beta_\mu$ and $h^{\mu,\beta}_{\gamma,\beta} (p') \leq h^{\mu,\beta}_{\gamma,\beta} (p) = h^{\gamma,\lambda}_{\gamma,\beta} (q)$
in $\PP^\beta_\gamma$. Thus, letting $q' = q \cdot h^{\mu,\beta}_{\gamma,\beta} (p') \in \PP^\lambda_\gamma$, we get that
$h^{\gamma,\lambda}_{\gamma,\beta} (q') = h^{\mu,\beta}_{\gamma,\beta} (p')$ and $p' \cdot q' \leq p \cdot q$ as required.

Preservation of correctness is obvious by the definition of (a dense subset of) $\PP^\lambda_\mu$ as the collection of formal products.
The statements about direct limits are straightforward for the same reason.

%The above arguments show in fact that $h^{\mu,\lambda}_{\mu,\beta} (p \cdot q) = p \cdot h^{\gamma,\lambda}_{\gamma,\beta} (q)$
%for any $\beta \geq \alpha$ where the pair $\alpha, \gamma$ witnesses $p \cdot q \in \PP^\lambda_\mu$.
%Similarly $h^{\mu,\lambda}_{\delta,\lambda} (p \cdot q) = h^{\mu,\alpha}_{\delta,\alpha} (p) \cdot q$ for any $\delta \geq \gamma$
%where $\alpha, \gamma$ witnesses $p \cdot q \in \PP^\lambda_\mu$. Thus, if $p \cdot q \in \PP^\lambda_\mu$, then
%$\supp (p \cdot q) = \supp (p) \cup \supp (q)$ so that $\P^\lambda$ has finite supports if $\la \PP^\lambda_\gamma : \gamma < \mu \ra$
%and all the $\P^\alpha$, $\alpha < \lambda$, do.
Finally assume all $\P^\alpha$, $\alpha < \lambda$, and $\la \PP^\lambda_\gamma : \gamma < \mu \ra$ have finite supports. 
Define $D^\lambda_\mu \sub \PP^\lambda_\mu$ by $p \cdot q \in D^\lambda_\mu$ if $p \in D^\alpha_\mu$, $q \in D^\lambda_\gamma$, and
$h^{\gamma,\lambda}_{\gamma,\alpha} (q) \leq h^{\mu,\alpha}_{\gamma,\alpha} (p)$ for some $\alpha < \lambda$ and $\gamma < \mu$.
Note that while this product is not of the form stipulated above, it is clearly equivalent to the condition $p' \cdot q$ where
$p' = p \cdot h^{\gamma,\lambda}_{\gamma,\alpha} (q)$ in the cBa $\PP^\lambda_\mu$.

We first need to see that this set is indeed dense. To this end, take $p \cdot q  \in \PP^\lambda_\mu$ with witness $\alpha, \gamma$.
Find $p' \leq p $ in $D^\alpha_\mu$ and then $q' \leq q \cdot h^{\mu,\alpha}_{\gamma,\alpha} (p')$ in $D^\lambda_\gamma$ and note that $p ' \cdot q' \in D^\lambda_\mu$ is the required
extension of $p \cdot q$. The inclusion relations (see (I) in the definition of ``finite supports") are obvious.

Next note that the argument for completeness of the embeddings presented above in fact shows that for $\delta < \mu$,
\[ h^{\mu,\lambda}_{\delta,\lambda} (p \cdot q) = \left\{ \begin{array}{ll} h^{\mu,\alpha}_{\delta,\alpha} (p) \cdot q & \mbox{ if } \delta \geq \gamma \\
   h^{\gamma,\lambda}_{\delta,\lambda} (q) & \mbox{ if } \delta \leq \gamma \end{array} \right. \]
Since $h^{\mu,\alpha}_{\delta,\alpha} (p) \in D^\alpha_\delta$ by (II) in the definition of ``finite supports", we see by assumption that
$h^{\mu,\lambda}_{\delta,\lambda} (p \cdot q) \in D^\lambda_\delta$ in either case, and (II) still holds. Furthermore, this implies
that $\supp (p \cdot q) \sub \supp (p) \cup \supp (p)$, so that supports are still finite (clause (III)). This completes the proof.
\end{proof}

We next need to argue that the limit construction of the previous lemma preserves the ccc. This is far from obvious.
In view of later applications (Theorem~\ref{a>d-meas} and Subsection~\ref{ultrapowers-further}) we do this in the special situation
when $\P^{\alpha +1}$ is the ultrapower of $\P^\alpha$ and the iterands in the $\P^\alpha$ are sufficiently simple.

We say a forcing notion $(\PP, \leq )$ is {\em Suslin $\sigma$-linked} if it is Suslin ccc and $\PP = \bigcup_n P_n$ where all
$P_n$ are linked (that is, any two elements of $P_n$ have a common extension) $\SIGMA^1_1$ sets. This implies 
``$P_n$ is linked" is a $\PI^1_1$ statement and therefore absolute. (Indeed, linkedness of $P_n$ is equivalent to
``$\forall x, y \; ( x \notin P_n \lor y \notin P_n \lor \neg ( x \bot y ))$".)  Also the statement ``$\PP =  \bigcup_n P_n$" is $\PI^1_2$ and absolute as well.

\begin{lem}  \label{iu-ccc}
Let $\mu$ and $\lambda$ be limit ordinals.
Let $\QQ = \bigcup_n Q_n$ be a Suslin $\sigma$-linked forcing notion. Assume $\P^\alpha = \la \PP^\alpha_\gamma : \gamma \leq \mu \ra$ are iterations, 
$\alpha \leq \lambda$, such that
\begin{itemize}
\item[{\rm{(a)}}] $\PP^\alpha_0 = \{ \zero, \one \}$ for all $\alpha$,
\item[{\rm{(b)}}] $\PP^0_\delta = \lim \dir_{\gamma < \delta} \PP^0_\gamma$ for limit $\delta$,
\item[{\rm{(c)}}] $\PP^\alpha_{\gamma + 1} = \PP^\alpha_\gamma \star \dot \QQ$ for all $\alpha$ and $\gamma$,
\item[{\rm{(d)}}] $\PP^{\alpha + 1}_\gamma = ( \PP^\alpha_\gamma)^\kappa / \D$ for all $\alpha$ and $\gamma$, and
\item[{\rm{(e)}}] $\PP^\beta_\gamma$ is built according to the proof of Lemma~\ref{iu-limit} for limit $\beta$ and $\gamma$.
\end{itemize}
Then
\begin{itemize}
\item[{\rm{(A)}}] if $\alpha \leq \beta$ and $\gamma \leq \delta$, then $\PP^\alpha_\gamma \embed \PP^\beta_\delta$,
\item[{\rm{(B)}}] projections in all diagrams $\la\PP^\alpha_\gamma , \PP^\alpha_\delta , \PP^\beta_\gamma, \PP^\beta_\delta \ra$,  $\alpha < \beta$ and $\gamma < \delta$, are correct,
\item[{\rm{(C)}}] all $\P^\alpha = \la \PP^\alpha_\gamma : \gamma \leq \mu \ra$ are iterations with finite supports as witnessed by $\la D^\alpha_\gamma : \gamma \leq \mu \ra$
   with $D^\alpha_\gamma \sub D^\beta_\gamma$ for  $\alpha < \beta$ and
   such that for all $\gamma < \delta$ and $\alpha < \beta$, if $\beta$ is a limit ordinal, and $p \in D^\alpha_\delta$ and $q \in D^\beta_\gamma$ are such that
   $h^{\gamma,\beta}_{\gamma,\alpha} (q) \leq h^{\delta,\alpha}_{\gamma,\alpha} (p)$, then $p \cdot q \in D^\beta_\delta$,
\item[{\rm{(D)}}] for limit $\delta$, $\PP^\alpha_\delta$ is the direct limit of the $\PP^\alpha_\gamma$, $\gamma < \delta$, unless $cf (\delta) = \kappa$,
\item[{\rm{(E)}}] for limit $\beta$, $\PP^\beta_\gamma$ is the direct limit of the $\PP^\alpha_\gamma$, $\alpha < \beta$, unless $cf(\beta) = \omega$, and
\item[{\rm{(F)}}] all $\PP^\alpha_\gamma$ satisfy property K (and thus are ccc).
\end{itemize}
\end{lem}

Let us note that (a) through (e) completely determine the iteration and that, by Lemma~\ref{ultra-Suslin}, there is no conflict between (c) and (d).
Also, (F) is the main point of this lemma, but we explicitly state (A) through (E) for later use (see the proof of Theorem~\ref{a>d-meas}).

\begin{proof}   
(A) Induction on $\beta \leq \lambda$. If $\beta $ is successor, this follows from Lemmas~\ref{ultra-embed} and~\ref{ultra-twostep}. If $\beta$ is limit, the successor case
is Lemma~\ref{Suslinccc-completeness} and the limit case, Lemma~\ref{iu-limit}.

(B) For successor $\beta$, this is Lemma~\ref{ultra-twostep}. If $\beta$ is limit, the successor case follows from the comment after the definition of correctness (in Subsection~\ref{prelim}) 
and the limit case, from Lemma~\ref{iu-limit}. 

(C) Induction on $\alpha \leq \lambda$. For $\alpha = 0$ this is well-known, but we recall the standard definition of the $D^0_\gamma$: $D^0_0 = \{ \one \}$,
$D^0_{\gamma + 1} = D^0_\gamma \star \dot \QQ$, and $D^0_\delta = \bigcup_{\gamma < \delta} D^0_\gamma$ for limit $\delta$. For successor $\alpha = \beta + 1$ see Lemma~\ref{ultra-iteration}:
$D^\alpha_\gamma = ( D^\beta_\gamma)^\kappa / \D$ for all $\gamma$. For limit $\alpha$, the $D^\alpha_\gamma$ are produced by recursion on $\gamma$
using Lemma~\ref{iu-limit}, which describes how to define $D^\alpha_\gamma$ for limit $\gamma$. If $\gamma = \delta + 1$ is successor, let $D^\alpha_\gamma =
D^\alpha_\delta \star \dot \QQ$ and check that the required properties are still satisfied.

(D) Assume $\delta$ is limit with $cf (\delta) \neq \kappa$ and make  induction on $\alpha \leq \lambda$. For $\alpha = 0$ this is (b), for successor $\alpha$, Lemma~\ref{ultra-iteration}, and
for limit $\alpha$, Lemma~\ref{iu-limit}.

(E) Assume $cf (\beta) > \omega$ and make induction on $\gamma \leq \mu$. For $\gamma = 0$ this is trivial by (a), and for limit $\gamma$ we use Lemma~\ref{iu-limit}.
So assume $\gamma = \delta + 1$ is successor. Then $\PP^\beta_\gamma = \PP^\beta_\delta \star \dot \QQ$ and $\PP^\beta_\delta$ is the direct limit of the
$\PP^\alpha_\delta$, $\alpha < \beta$, by induction hypothesis. Thus, if $\dot q$ is a $\PP^\beta_\delta$-name for an element of $\dot \QQ$, then by the
ccc (see (F)) and $cf (\beta) > \omega$, $\dot q$ is a $\PP^\alpha_\delta$-name for some $\alpha < \beta$. Hence if $(p, \dot q) \in \PP^\beta_\gamma$, then $(p, \dot q) \in \PP^\alpha_\gamma$
for some $\alpha < \beta$, as required.

(F) By recursion on $\alpha \leq \lambda$, we define $I^\alpha$, $E^\alpha_\gamma$, $\gamma \leq \mu$, and $s^\alpha$, such that
\begin{enumerate}
\item the $I^\alpha$ are linear orders, $I^0 = \mu$, and $I^\alpha \sub I^\beta$ for $\alpha \leq \beta$,
\item the $E^\alpha_\gamma $ are dense subsets of $D^\alpha_\gamma$ and $E^\alpha_\gamma \sub E^\beta_\delta$ for $\alpha \leq \beta$ and $\gamma \leq \delta$,
\item the $s^\alpha$ are functions with domain $E^\alpha_\mu$ and $s^\alpha \sub s^\beta$ for $\alpha \leq \beta$,
\item if $p \in E^\alpha_\gamma$, $\gamma \leq \mu$, then $s^\alpha (p) : I^\alpha \to \omega$ is a finite partial function with $\dom (s^\alpha (p) ) \cap \mu =
   \{ \delta : \delta + 1 \in \supp (p) \}$,
\item if $p \in E^\alpha_\gamma$, $\gamma \leq \mu$, $\delta < \gamma$ with $\delta + 1 \in \supp (p)$, then $h^{\gamma,\alpha}_{\delta,\alpha} (p)
   \forces p(\delta) \in \dot Q_{s^\alpha (p)(\delta) }$,
\item if $p,q \in E^\alpha_\gamma$, $\gamma \leq \mu$, $\delta < \gamma$, $s^\alpha (p)$ and $s^\alpha (q)$ agree on their common domain,
   and $r_0 \leq h^{\gamma,\alpha}_{\delta,\alpha} (p) ,  h^{\gamma,\alpha}_{\delta,\alpha} (q)$ in $\PP^\alpha_\delta$, 
   then there is $r \leq p,q$ in $\PP^\alpha_\gamma$ such that $h^{\gamma,\alpha}_{\delta,\alpha} (r) = r_0$.
\end{enumerate}
Note that for $\delta = 0$, the last clause means that $p$ and $q$ in $\PP^\alpha_\gamma$ are compatible if $s^\alpha (p)$ and $s^\alpha (q)$ agree
on their common domain. By the $\Delta$-system lemma, it is then immediate that $\PP^\alpha_\gamma$ has property K and the ccc follows.
So it suffices to carry out the recursion.

\underline{Basic step $\alpha = 0$}. Let $I^0: = \mu$. Define $E^0_\gamma$ and $s^0$ with the required properties by recursion on $\gamma \leq \mu$.

If $\gamma = 0$ there is nothing to do.

Assume $\gamma = \delta + 1$ is successor, and $E^0_\delta$ and $s^0 \re E^0_\delta$ have been defined. 
Let \[E^0_\gamma = E^0_\delta \cup \{ p = (p_0 , \dot x) \in \PP^0_\gamma : p_0 \in E^0_\delta \mbox{ and } p_0 \forces ``\dot x \in \dot Q_n" \mbox{ for some } n \}. \]
Clearly, $E^0_\gamma$ is dense in $D^0_\gamma$. Next, for such $p = (p_0 , \dot x) \in E^0_\gamma$, notice that $\delta + 1 \in \supp (p)$,
let $\dom (s^0 (p)) = \dom ( s^0(p_0)) \cup \{ \delta  \}$
and define $s^0 (p)$ such that $s^0 (p_0) \sub s^0 (p)$ and $s^0 (p) ( \delta ) = n$ where $p_0 \forces \dot x \in \dot Q_n$. Then (5) is
satisfied and we need to show (6).

If $p = (p_0 , \dot x) , q = (q_0 , \dot y) \in E^0_\gamma$, $s^0 (p)$ and $s^0 (q)$ agree on their common domain,
and $r_0 \leq p_0, q_0$, then $r_0 \forces \dot x , \dot y \in \dot Q_n$ where $n = s^0 (p) (\delta) = s^0 (q) (\delta)$ so that $r_0$ forces $ \dot x , \dot y$ to be compatible
and there is a name $\dot z$ such that $r_0 \forces \dot z \leq \dot x , \dot y $. Letting $r = ( r_0, \dot z)$, we see that $r$ is as required.

Finally assume $\gamma$ is a limit ordinal. Since $\PP^0_\gamma = \lim \dir_{\delta < \gamma} \PP^0_\delta$, 
$E^0_\gamma = \bigcup_{\delta < \gamma} E^0_\delta$ is dense in $D^0_\gamma$ (see (C)), and properties (3) through (6) hold vacuously.

\underline{Successor step $\alpha = \beta + 1$}. Then $\PP^\alpha_\gamma = ( \PP^\beta_\gamma)^\kappa / \D$ for all $\gamma \leq \mu$.
We let $I^\alpha = (I^\beta)^\kappa / \D$. Clearly, $I^\alpha$ is linearly ordered with $[v] \leq [w]$ if $\{ \xi < \kappa : v (\xi) \leq w (\xi) \} \in \D$
for $v,w : \kappa \to I^\beta$, and $I^\beta \sub I^\alpha$. Next let \[ E^\alpha_\gamma = \{ [f] \in \PP^\alpha_\gamma : f : \kappa \to \PP^\beta_\gamma \mbox{ 
and }\{ \xi < \kappa :  f(\xi) \in E^\beta_\gamma \} \in \D \}. \] Clearly $E^\alpha_\gamma$ is dense in $D^\alpha_\gamma$ and $E^\beta_\delta \sub E^\alpha_\gamma$ for $\delta
\leq \gamma$. For such $f$, and for all $\xi$ with $f (\xi) \in E^\beta_\gamma$, $\dom (s^\beta (f (\xi))) \sub I^\beta$ is finite, say
$ \dom (s^\beta  (f (\xi))) = \{ i^\xi_0 < ... < i^\xi_{m_\xi -1} \}$. By $\omega_1$-completeness of $\D$, there is $m$ such that
$\{ \xi < \kappa : m_\xi = m \}$ belongs to $\D$. Define $v_j : \kappa \to I^\beta$ by
\[  v_j (\xi) = \left\{ \begin{array}{ll} i^\xi_j & \mbox{ if this is defined} \\ 0 & \mbox{ otherwise} \end{array} \right. \]
for $j < m$. (Note that the first case occurs $\D$-almost everywhere.) Then $\{ [v_0] < ...<  [v_{m-1}] \} \sub I^\alpha$ and we let $\dom (s^\alpha ( [f]))
= \{ [v_j] : j < m \}$. Applying once more $\omega_1$-completeness of $\D$, we see that there are $n_j$, $j < m$, such that
$\{ \xi < \kappa : s^\beta (f (\xi)) (i^\xi_j) = n_j \} \in \D$. Thus we let $s^\alpha ( [f] ) ( [ v_j] ) = n_j$ for $j < m$. Clearly
$s^\beta \sub s^\alpha$.

To see (5), if $\delta < \gamma$ and $s^\alpha ( [f] ) (\delta) = n$, then $\{ \xi < \kappa : s^\beta (f (\xi )) (\delta) = n \} \in \D$. Also,
\[  h^{\gamma, \alpha}_{\delta +1, \alpha} ( [f]) =  \la h^{\gamma,\beta}_{\delta + 1 , \beta} (f (\xi )) : \xi < \kappa \ra / \D = 
\la ( h^{\gamma,\beta}_{\delta , \beta} (f (\xi )) , \dot x^\xi ) : \xi < \kappa \ra / \D = (  h^{\gamma, \alpha}_{\delta, \alpha} ( [f]) , \dot y) \in (\PP^\beta_\delta)^\kappa / \D \star \dot \QQ \]
where $\dot y = \la \dot x^\xi : \xi < \kappa \ra / \D$ (see Lemmata~\ref{ultra-twostep} and~\ref{ultra-Suslin}). 
By induction hypothesis (5) we know that $\{ \xi < \kappa : h^{\gamma,\beta}_{\delta , \beta} (f (\xi )) \forces \dot x^\xi \in \dot Q_n \} \in \D$.
Therefore $ h^{\gamma, \alpha}_{\delta, \alpha} ( [f]) \forces \dot y \in \dot Q_n$.

To prove (6), assume $\delta < \gamma$, $[f], [g] \in E^\alpha_\gamma$, $s^\alpha ([f])$ and $s^\alpha ([g])$ agree on their common domain, and
$[h_0] \leq h^{\gamma,\alpha}_{\delta,\alpha} ( [f]) , h^{\gamma,\alpha}_{\delta,\alpha} ( [g]) $ in $ \PP^\alpha_\delta$. Let
$\{ [v_0] < ... < [v_{n-1}] \}$ list the common domain of $s^\alpha ([f])$ and $s^\alpha ([g])$. Then
\[ \{ \xi < \kappa : V(\xi) = \{ v_0 (\xi) < ... < v_{n-1} (\xi) \} = \dom(  s^\beta (f (\xi)) ) \cap \dom ( s^\beta (g (\xi)) ) \mbox{ and }
s^\beta (f (\xi)) \re V (\xi) =  s^\beta (g (\xi))\re V (\xi) \} \]
belongs to $\D$. Also, the set $\{ \xi < \kappa : h_0 (\xi) \leq h^{\gamma,\beta}_{\delta,\beta} (f (\xi)), h^{\gamma,\beta}_{\delta,\beta} (g (\xi)) \}$ belongs to $\D$.
For $\xi$ which belong to both sets we find, by induction hypothesis (6), $h (\xi)\in \PP^\beta_\gamma$ with $h (\xi)  \leq f(\xi), g(\xi) $
and $h^{\gamma,\beta}_{\delta,\beta} (h (\xi)) = h_0 (\xi)$. So $[h] \leq [f] , [g]$ and $h^{\gamma,\alpha}_{\delta,\alpha} ( [h]) = [ h_0]$,
as required.

\underline{Limit step $\alpha$}. Let $I^\alpha = \bigcup_{\beta < \alpha} I^\beta$, equipped with the obvious ordering. As in the basic step, we
define $E^\alpha_\gamma$ and $s^\alpha$ by recursion on $\gamma \leq \mu$.

The cases $\gamma = 0$ and $\gamma = \delta + 1$ are identical to the basic step. The only difference is that, this time,
$E^\alpha_\gamma$ must contain all $E^\beta_\gamma$, and that $s^\alpha$ must extend all $s^\beta$, for $\beta < \alpha$.

So assume $\gamma$ is a limit ordinal, and $E^\alpha_\delta$ and $s^\alpha \re E^\alpha_\delta$ have been defined for $\delta < \gamma$.
Since supports are finite by (C), %(see also Lemmas~\ref{ultra-iteration} and~\ref{iu-limit}), 
we know that $\bar \PP^\alpha_\gamma \embed \PP^\alpha_\gamma$ where $\bar \PP^\alpha_\gamma : = 
\lim \dir_{\delta < \gamma} \PP^\alpha_\delta$, by Lemma~\ref{finite-supports}. 
By the proof of Lemma~\ref{iu-limit}, elements of $\PP^\alpha_\gamma$ are formal products $p \cdot \bar p$ with $p \in \PP^\beta_\gamma$
for some $\beta < \alpha$, $\bar p \in \bar \PP^\alpha_\gamma$, and $h^{\bar \gamma, \alpha}_{\bar \gamma ,\beta} (\bar p)
= h^{\gamma, \beta}_{\bar \gamma, \beta} (p)$ (where we use $\bar \gamma$ as an index for the direct
limit $\bar \PP^j_\gamma$ of the $\PP^j_\delta$, $\delta < \gamma$, which completely embeds into $\PP^j_\gamma$, where $j = \alpha, \beta$).
By strengthening $p$ and $\bar p$, if necessary, we may assume $p \in E^\beta_\gamma$.
By further strengthening $\bar p$, we may assume $\bar p \in \bar E^\alpha_\gamma : = \bigcup_{\delta < \gamma } E^\alpha_\gamma$.
In general, we will then only have $h^{\bar \gamma, \alpha}_{\bar \gamma ,\beta} (\bar p)
\leq  h^{\gamma, \beta}_{\bar \gamma, \beta} (p)$, but this does not concern us because the collection of formal products
satisfying this weaker condition is obviously forcing equivalent with the original $\PP^\alpha_\gamma$. Hence, if we let
$E^\alpha_\gamma$ consist of formal products $p \cdot \bar p$ with $p \in E^\beta_\gamma$ for some $\beta < \alpha$,
$\bar p \in \bar E^\alpha_\gamma$, and $h^{\bar \gamma, \alpha}_{\bar \gamma ,\beta} (\bar p)
\leq  h^{\gamma, \beta}_{\bar \gamma, \beta} (p)$, then $E^\alpha_\gamma$ is dense in $\PP^\alpha_\gamma$. Also $E^\alpha_\gamma \sub D^\alpha_\gamma$ by the proof 
of Lemma~\ref{iu-limit}.
Clearly $\bar E^\alpha_\gamma \sub E^\alpha_\gamma$ and $E^\beta_\gamma \sub E^\alpha_\gamma$ for $\beta < \alpha$.
For such $p \cdot \bar p \in E^\alpha_\gamma$, we define $s^\alpha ( p \cdot \bar p)$ by 
$\dom (s^\alpha ( p \cdot \bar p)) = \dom (s^\alpha (\bar p)) \cup \{ i \in \dom (s^\beta (p)) :i \geq \delta$ for all $\delta < \gamma\}$
and
\[ s^\alpha ( p \cdot \bar p) (i) = \left\{ \begin{array}{ll} s^\alpha (\bar p) (i) & \mbox{ for } i \in \dom (s^\alpha (\bar p)) \\
   s^\beta (p) (i) & \mbox{ for } i \in \dom (s^\beta ( p)) \mbox{ with } i > \delta \mbox{ for all } \delta < \gamma. \end{array} \right. \]
(5) is immediate by induction hypothesis (5) for $\bar p$.

To prove (6), let $\delta < \gamma$, $p \cdot \bar p, q \cdot \bar q \in E^\alpha_\gamma$, $s^\alpha ( p \cdot \bar p) $ and $s^\alpha (q \cdot \bar q)$
agree on their common domain, and $r_0 \leq h^{\gamma, \alpha}_{\delta,\alpha} (p \cdot \bar p) , h^{\gamma, \alpha}_{\delta,\alpha} (q\cdot \bar q)$ in $\PP^\alpha_\delta$.
It is immediate that $h^{\gamma,\alpha}_{\delta,\alpha} (p \cdot \bar p) = h^{\bar\gamma, \alpha}_{\delta,\alpha} (\bar p)$ and similarly for $q \cdot \bar q$.
Thus, by induction hypothesis (6), there is $\bar r \leq \bar p, \bar q$ in $\bar \PP^\alpha_\gamma$ such that $h^{\bar\gamma, \alpha}_{\delta,\alpha} (\bar r) = r_0$.
Let $\beta_p,\beta_q < \alpha$ be such that $p \in E^{\beta_p}_\gamma$ and $q \in E^{\beta_q}_\gamma$. Without loss of generality $\beta_p \leq \beta_q$,
and we let $\beta : = \beta_q$. By correctness $h^{\gamma,\beta}_{\bar\gamma,\beta} (p) = h^{\gamma,\beta_p}_{\bar\gamma,\beta_p} (p)$. So we know
that $h^{\bar\gamma, \alpha}_{\bar\gamma,\beta} (\bar r) \leq h^{\bar\gamma, \alpha}_{\bar\gamma,\beta} (\bar p) \leq 
h^{\bar\gamma, \alpha}_{\bar\gamma,\beta_p} (\bar p) \leq h^{\gamma, \beta_p}_{\bar\gamma,\beta_p} ( p) = h^{\gamma,\beta}_{\bar\gamma,\beta} (p)$
and $h^{\bar\gamma, \alpha}_{\bar\gamma,\beta} (\bar r) \leq h^{\bar\gamma, \alpha}_{\bar\gamma,\beta} (\bar q) \leq 
h^{\gamma,\beta}_{\bar\gamma,\beta} (q)$. Thus by induction hypothesis (6) there is $r \leq p,q$ in $\PP^\beta_\gamma$ such that
$h^{\gamma,\beta}_{\bar\gamma,\beta} (r) = h^{\bar\gamma, \alpha}_{\bar\gamma,\beta} (\bar r)$. 
Hence $h^{\gamma,\alpha}_{\delta,\alpha} (r \cdot \bar r) = h^{\bar\gamma, \alpha}_{\delta,\alpha} (\bar r) = r_0$, and $r \cdot \bar r \leq p \cdot \bar p,
q \cdot \bar q$ is as required.
\end{proof}
 
%\bigskip

%%%

\subsection{The consistency of $\dd < \aa$ from a measurable cardinal}
\label{a>d-measurable}

Recall that {\em Hechler forcing} $\DD$ consists of $(s,x) \in \omlom \times \omom$ with $s \sub x$ ordered by $(t, y) \leq (s,x)$ if
$t \supseteq s$ and $y \geq x$ everywhere. $\DD$ is a Suslin $\sigma$-centered forcing adding a dominating real. In particular it
is Suslin $\sigma$-linked.

We outline the proof of:

\begin{thm}[Shelah~\cite{Sh04}]  \label{a>d-meas}
Assume ZFC + ``there is a measurable cardinal" is consistent. Then so is $\dd < \aa$.
More explicitly, if GCH holds, $\kappa$ is measurable and $\lambda > \mu > \kappa$ are regular, then there is
a ccc forcing extension satisfying $ \bb = \dd = \mu$ and $\aa = \cc = \lambda$.
\end{thm}

\begin{proof}
As usual, let $\D$ be a $\kappa$-complete ultrafilter on $\kappa$.
By recursion on $\alpha \leq \lambda$ we construct iterations (with finite supports) $\P^\alpha = \la \PP^\alpha_\gamma : \gamma \leq \mu \ra$ such that
\begin{enumerate}
\item $\PP^\alpha_0 = \{ \zero, \one \}$ for all $\alpha$,
\item $\PP^0_\delta = \lim \dir_{\gamma < \delta} \PP^0_\gamma$ for limit $\delta$,
\item $\PP^\alpha_{\gamma + 1} = \PP^\alpha_\gamma \star \dot \DD$ for all $\alpha$ and $\gamma$,
\item $\PP^{\alpha +1}_\gamma = ( \PP^\alpha_\gamma)^\kappa / \D$ for all $\alpha$ and $\gamma$,
\item $\PP^\beta_\delta$ is built according to the proof of Lemma~\ref{iu-limit} for limit $\beta$ and $\delta$.
%\item if $\alpha \leq \beta$ and $\gamma \leq \delta$, then $\PP^\alpha_\gamma \embed \PP^\beta_\delta$,
%\item projections in all diagrams $\la\PP^\alpha_\gamma , \PP^\alpha_\delta , \PP^\beta_\gamma, \PP^\beta_\delta \ra$,  $\alpha < \beta$ and $\gamma < \delta$, are correct.
\end{enumerate}
%For $\beta = 0$ construct $\P^0$ according to the first two clauses, that is, $\P^0$ is a $\mu$-stage finite support iteration of Hechler forcing.
%If $\beta = \alpha + 1$, $\P^\beta$ is the ultrapower of $\P^\alpha$, see Lemma~\ref{ultra-iteration}. Thus (3) is satisfied, and (5) and (6) follow
%from Lemmas~\ref{ultra-embed} and~\ref{ultra-twostep}. Furthermore, (2) also holds for by definition, the induction hypothesis for (2), and Lemma~\ref{ultra-Suslin} we have
%\[ \PP^\beta_{\gamma + 1} = ( \PP^\alpha_{\gamma + 1 })^\kappa / \D = ( \PP^\alpha_{\gamma  } \star \dot \DD)^\kappa / \D = ( \PP^\alpha_{\gamma  })^\kappa / \D  \star \dot \DD
%= \PP^\beta_\gamma \star \dot \DD. \]
%Finally, for limit $\beta$, we define $\PP^\beta_\delta$ for successor $\delta = \gamma + 1$ by (2), and for limit $\delta$, by (4). 
%(5) and (6) follow from Lemma~\ref{Suslinccc-completeness} and the comment after the definition of correctness (in Subsection~\ref{prelim}) in the first
%case, and from Lemma~\ref{iu-limit}, in the second case. This completes the recursive construction.
This gives a complete description of the construction of the $\P^\alpha$. Also note that (1) through (5) correspond exactly to (a) through (e) in the
assumptions of Lemma~\ref{iu-ccc} for the special case $\QQ = \DD$. Therefore, all conclusions of the lemma hold and, in particular,
all $\PP^\alpha_\gamma$ are ccc and so is $\PP : = \PP^\lambda_\mu$. Furthermore, by (D) and (E) of this lemma, we see that
$\PP^\lambda_\mu$ is a direct limit of the $\PP^\lambda_\gamma$, $\gamma < \mu$, and of the $\PP^\alpha_\mu$, $\alpha < \lambda$.

We next note that $| \PP | = \lambda$. Indeed, by induction on $\alpha$, we have that $\max \{ |\alpha| , \mu \} \leq | \PP^\alpha_\mu | \leq \max \{ |\alpha| , \mu \}^\kappa$
for all $\alpha \leq \lambda$: $|\PP^0_\mu| = \mu = \mu^\kappa$ is obvious. For successor, the formula follows from $| \PP^{\alpha + 1}_\mu | \leq | \PP^\alpha_\mu|^\kappa$,
and, for limit, from $|\PP^\beta_\mu| = \sum_{\alpha < \beta} | \PP^\alpha_\mu |$. In particular, a standard argument with nice names gives $\cc \leq \lambda$.

The last iteration $\P^\lambda = \la \PP^\lambda_\gamma : \gamma \leq \mu \ra$ cofinally often adds a dominating real (more explicitly, $\PP^\lambda_{\gamma + 1}$ adds a Hechler
real dominating the reals of the $\PP^\lambda_\gamma$-extension), and since $\PP^\lambda_\mu$ is the direct limit of this iteration $\bb = \dd = \mu$ holds in the extension. 
(In fact, the sequence of $\mu$ Hechler reals is already added by $\P^0$, see also Lemma~\ref{ultra-scales}).

To see $\aa \geq \lambda$, we use Lemma~\ref{ultra-mad}: assume $\A$ is an a.d. family of size $\nu$ for some $\nu < \lambda$ in the generic
extension. If $\nu < \mu$ then $\A$ is not maximal because $\bb \leq \aa$ in ZFC. If $\nu \geq \mu > \kappa$, by the ccc, the regularity of $\lambda$, and the
fact that $\PP^\lambda_\mu$ is the direct limit of the $\PP^\alpha_\mu$, $\alpha < \lambda$,
there is an $\alpha < \lambda$ such that $\dot \A$ is a $\PP^\alpha_\mu$-name. By Lemma~\ref{ultra-mad}, we then see that $\PP^{\alpha + 1}_\mu =
( \PP^\alpha_\mu)^\kappa / \D$ forces that $\dot \A$ is not maximal. This completes the proof.
\end{proof}

%%%

\subsection{Further results}
\label{ultrapowers-further}

We collect a number of results obtained by related methods. 

In his original work~\cite{Sh04}, Shelah also obtains the consistency of $\uu < \aa$
by basically the same method.
Recall here that, given a nonprincipal ultrafilter $\U$ on $\omega$, $\B \sub \U$ is a {\em base} of $\U$ if for
every $A \in \U$ there is $B \in \B$ such that $B \sub A$. If we only require that $\B \sub \omoms$,
$\B$ is called a {\em $\pi$-base}. The {\em character} $\chi (\U)$ ({\em $\pi$-character} $\pi\chi (\U)$, respectively) of $\U$ is the least size of a 
base ($\pi$-base, resp.) of $\U$. We define
$\uu : = \min \{ \chi (\U) : \U \mbox{ is a nonprincipal ultrafilter on } \omega \}, \mbox{ the {\em ultrafilter number}}.$
%It is well-known that if $\rr$ is the reaping number, then 
%\[ \rr = \min \{ \pi\chi (\U) : \U \mbox{ is a nonprincipal ultrafilter on } \omega \} \]

\begin{thm}[Shelah~\cite{Sh04}]  \label{a>u-meas}
Assume ZFC + ``there is a measurable cardinal" is consistent. Then so is $\uu < \aa$.
More explicitly, if GCH holds, $\kappa$ is measurable and $\lambda > \mu > \kappa$ are regular, then there is
a ccc forcing extension satisfying $\uu = \bb = \dd = \mu$ and $\aa = \cc = \lambda$.
\end{thm}

This is proved by replacing Hechler forcing $\DD$ with {\em Laver forcing} $\LL_\U$ with an ultrafilter $\U$ in the iterated
forcing construction of the proof of Theorem~\ref{a>d-meas}. For details see also~\cite{Br07}. 

In two subsequent papers (\cite{Sh08}, \cite{Sh11}), Shelah used this construction for obtaining several
results about non-convexity of the {\em character spectrum} 
$ \Spec (\chi) = \{ \chi (\U) : \U \mbox{ is a nonprincipal ultrafilter on } \omega \}. $

\begin{thm}[Shelah~\cite{Sh08}]   \label{spec1}
Assume ZFC + ``there  is a measurable cardinal" is consistent. Then so is ``$\Spec (\chi)$ is not convex".
More explicitly, if GCH holds, $\kappa$ is measurable and $\lambda > \kappa > \mu$ are uncountable regular, then there is
a ccc forcing extension satisfying $\uu = \bb = \dd = \mu$,  $\cc = \lambda$, $\{ \mu , \lambda \} \sub \Spec (\chi)$ and $\kappa \notin \Spec(\chi)$.
\end{thm}

Using two measurables, this can be combined with Theorem~\ref{a>u-meas}~\cite[Theorem 1.1]{Sh08}. 
Furthermore, there is an extension to a variant of the $\pi$-character spectrum~\cite[Theorem 2.5]{Sh08}. Again, details can be
found as well in~\cite{Br07}. 

\begin{thm}[Shelah~\cite{Sh11}]   \label{spec2}
Given two disjoint sets $\Theta_1$ and $\Theta_2$ of regular cardinals such that $\theta^{< \theta} = \theta$ for all $\theta \in \Theta_1$ and
all members of $\Theta_2$ are measurable cardinals, there is a partial order forcing $\Theta_1 \sub \Spec (\chi)$ and $\Theta_2 \cap \Spec (\chi) = \emptyset$.
\end{thm}

This is a generalization of Theorem~\ref{spec1}, obtained by a product of a ccc forcing with iterated ultrapowers (to guarantee the
measurables in $\Theta_2$ are not characters) and a forcing not adding reals but adjoining ultrafilters of character in $\Theta_1$.

\begin{thm}[Shelah~\cite{Sh11}]   \label{spec3}
Assuming the consistency of infinitely many strongly compact cardinals, for any $A \sub \omega \setminus \{ 0 \}$ it is consistent that 
$\{ \aleph_n : n \in A \}$ is the set of characters below $\aleph_\omega$.
\end{thm}

This is proved by combining the previous theorem with a product of Levy collapses.
Some of these results answer questions originally addressed in~\cite{BS99}.

A family $\SSS \sub \omoms$ is a {\em splitting family} if for all $A \in \omoms$ there is $B \in \SSS$ with $| A\cap B| = |A \sem B| = \omega$
(we say {\em $B$ splits $A$}). The {\em splitting number} $\ss$ is the smallest cardinality of a splitting family. Iterating the
ultrapower construction as in Subsection~\ref{ultrapowers-iterations} in a matrix and also destroying splitting families
while preserving an unbounded family added by the first iteration, but taking direct limits in the limit step (this makes the
complex Lemmas~\ref{iu-limit} and~\ref{iu-ccc} unnecessary) one obtains:

\begin{thm}[Brendle and Fischer~\cite{BF11}]   \label{b<as}
Assume GCH holds, $\kappa$ is measurable and $\lambda > \mu > \kappa$ are regular cardinals. Then there is
a ccc forcing extension in which $\bb = \mu < \ss = \aa = \cc = \lambda$.
\end{thm}

This generalizes an older result of Shelah~\cite{Sh84} who proved that $\bb = \aleph_1 < \ss = \aa = \cc = \aleph_2$ is consistent
with a countable support iteration of proper forcing.

%%%%%%%%%%%%%%%%%%%%%%%

%

%

%

%

%

%

%

%%%%%%%%%%%%%%%

\section{Templates}
\label{templates}

In standard well-ordered iterations $\la \PP_\alpha , \dot \QQ_\alpha : \alpha < \mu \ra$ where $\mu$ is an ordinal,
the {\em initial segments} of the iteration $\PP_\alpha$ and the {\em iterands} $\dot \QQ_\alpha$ are handed down by the
same recursion, with the latter being $\PP_\alpha$-names, one defines $\PP_{\alpha + 1} = \PP_\alpha \star \dot \QQ_\alpha$,
and then one has to specify what to do in the limit step. However, while wellfoundedness is a necessity for recursively defining
the initial segments in many cases, this is not so for the iterands, and they may well be indexed by a non-wellfounded structure.
This is what {\em iterations along templates} do. The index set of their iterands is an arbitrary linear order $L$, and
their initial segments are produced by specifying a wellfounded subfamily $\I$ of $\P (L)$ and then recursively 
defining the segment $\PP \re A$ for $A \in \I$. Since we need to consider the same iterand $\QQ_x$ for $x \in L$ over various initial
extensions, definability of such iterands is crucial, and iteration along templates can be seen as a generalization of
finite support iteration of Suslin ccc forcing. See Subsection~\ref{template-method} for details (the exposition follows
largely~\cite{Br05}).

As an application, we provide a complete proof  of Shelah's proof of the consistency of $\dd < \aa$ in ZFC~\cite{Sh04} in Subsection~\ref{d<a-ZFC}
(Theorem~\ref{a>d-ZFC}), following roughly~\cite{Br02}. Subsection~\ref{templates-further} presents further results
obtained by the template method.

%bigskip

\subsection{The template method}
\label{template-method}

Let $\la L, \leq_L \ra$ be a linear order. For $x \in L$ let $L_x = \{ y \in L : y <_L x \}$ be the {\em initial segment} determined by $x$.

\begin{defi}
A {\em template} is a pair $(L, \bar \I = \{ \I_x : x \in L \} )$ such that $\la L, \leq_L \ra$ is a linear order, $\I_x \sub \P(L_x)$ for all
$x \in L$, and
\begin{enumerate}
\item $\I_x$ contains $\emptyset$, all singletons, and is closed under unions and intersections,
\item $\I_x \sub \I_y$ for $x <_L y$,
\item $\I : = \bigcup_{x \in L} \I_x \cup \{ L \}$ is wellfounded with respect to inclusion, as witnessed by the {\em depth function} $\Dp_\I : \I \to \On$.
\end{enumerate}
If $A \sub L$ and $x \in L$, we define $\I_x \re A : = \{ B \cap A : B \in \I_x \}$, the {\em trace} of $\I_x$ on $A$
and let $\bar \I \re A = \{ \I_x \re A : x \in A \}$.
\end{defi}

$L$ is meant to be an index set for the iterands of an iteration while $\I$ describes in a sense the support.
The wellfoundedness condition 3 is crucial because it allows for a recursive definition of the iteration.
By clause 1, a real $r_y$ added at stage $y$ will be generic over a real $r_x$ added at stage $x$ for $x <_L y$ while
by clause 2, if $r_x$ is generic over some initial stage of the iteration then so is $r_y$. Conditions 1 and 2 together
imply that $\I$ is also closed under intersections and unions.

Note that $(A, \bar \I \re A)$ is a template as well and, if $B \in \I$ and $C = B \cap A \in \I \re A$, then $\Dp_{\I\re A} (C) \leq \Dp_\I (B)$.
Thus, with every $A \sub L$ we can associate its {\em depth}, $\Dp (A) := \Dp_{\I\re A} (A)$. $\Dp$ has the following properties:

\begin{lem}   \label{depth}
\begin{enumerate}
\item If $B \in \I \re A$ then for any $C \sub B$, $C \in \I \re B$ iff $C \in \I \re A$. 
\item If $B \in \I \re A$, then $\Dp_{\I \re A} (B) = \Dp (B)$.
\item If $B \sub A$ then $\Dp (B) \leq \Dp (A)$ and if additionally $B \in \I \re A$ and $B \subsetneq A$ then $\Dp (B) < \Dp (A)$.
\item If $B = B_0 \cup \{ x \} \sub A$, $x \notin B_0$, and $B_0,B \in \I \re A$ then $\Dp (B) = \Dp (B_0) + 1$.
\end{enumerate}
\end{lem}

\begin{proof}
(1) If $C \in \I \re A$ then there is $D \in \I$ with $C = A \cap D$. Thus $C = B \cap D$ and $C \in \I \re B$ follows.

On the other hand if $C \in \I \re B$ then there is $D \in \I$ with $C = B \cap D$. Also there is $D' \in \I$ with $B = A \cap D'$.
Thus $C = A \cap (D' \cap D)$ and therefore $C \in \I \re A$.

(2) is immediate by (1), and (3) is obvious.

(4) $\Dp (B) \geq \Dp (B_0) + 1$ follows by (3), and $\Dp (B) \leq \Dp (B_0) + 1$ is shown by induction on $\Dp (B_0)$. Suppose
$\Dp (B_0) = \alpha$. Let $C \subsetneq B$ with $C \in \I \re B$ and $x \in C$. Put $C_0 = C \cap B_0 = C \sem \{ x \}$. 
Then $C_0 \in \I \re B_0$ and $C_0 \subsetneq B_0$, and therefore $\Dp (C_0) < \alpha$ by (3). Since also $C_0 , C \in \I \re A$, 
by induction hypothesis $\Dp (C)  \leq \Dp (C_0) + 1 \leq \alpha$. Thus $\Dp (B) \leq \alpha + 1$.
\end{proof}

Call a Suslin ccc forcing notion  {\em $\QQ$ correctness-preserving} if for every diagram $\la \PP_i \ra$ with correct projections,
the canonical projections in the diagram $\la \PP_i \star \dot \QQ \ra$ are correct as well. It is unclear whether there are not correctness-preserving
Suslin ccc forcing notions (see~\cite[Conjecture 1.5]{Br05} and~\cite[Conjecture 1]{Brta}).

\begin{defthm}   \label{template-main}
Assume $(L,\bar \I)$ is a template. Also assume $( \QQ_x : x \in L )$ is a family of correctness-preserving Suslin ccc forcing notions
whose definition lies in the ground model. By recursion-induction on $\Dp (A)$, $A \sub L$,
\begin{enumerate}
\item we define the partial order $\PP \re A$ (both the set of elements and the order relation),
\item we prove that $\PP \re D \sub \PP \re A$ and $\leq_{\PP \re D} = \leq_{\PP \re A} \cap ( \PP \re D)^2$ for $D \sub A$,
\item we prove that $\PP \re A$ is transitive, 
\item we describe how $\PP \re A$ is obtained from $\PP \re B$ where $B \subsetneq A$ with $B \in \I \re A$ (so $\Dp (B) < \Dp (A)$),
\item we prove that $\PP \re D \embed \PP \re A$ for $D \sub A$,
\item we prove that for $D \sub L$ with $\Dp (D) \leq \Dp (A)$, we have $\PP \re (A \cap D) = \PP \re A \cap \PP \re D$,
\item we prove correctness in the sense that if $A' , D \subsetneq A$,  and $D' = A' \cap D$, then projections in the diagram 
   $\la \PP \re D', \PP \re A', \PP \re D, \PP \re A\ra$ are correct. 
\end{enumerate}
$\PP \re L$ is called {\em iteration along a template}.
\end{defthm}

\begin{proof}[Definition and Proof]
(1) First consider the case $\Dp (A) = 0$. This is equivalent to $A = \emptyset$. Let $\PP \re \emptyset = \{ \emptyset \}$.

So assume $\Dp (A) > 0$. Let $\PP \re A$ consist of all finite partial functions $p$ with domain contained in $A$ and such that, putting
$x = \max(\dom (p))$, there is $B \in \I_x \re A$ (so $\Dp (B) < \Dp (A)$) such that $p \re (A \cap L_x) \in \PP \re B$ and $p(x)$ is a
$\PP \re B$-name for a condition in $\dot \QQ_x$ (where we construe $\dot \QQ_x$ as a $\PP \re B$-name as well). 

The ordering on $\PP \re A$ is given as follows. $q \leq_{\PP \re A} p$ if $\dom (q) \supseteq \dom (p)$ and, putting
$x = \max(\dom(q))$, there is $B \in \I_x \re A$ such that $q \re (A \cap L_x) \in \PP \re B$ and
\begin{itemize}
\item either $x \notin \dom (p)$, $p \in \PP \re B$, and $q \re (A \cap L_x) \leq_{\PP \re B} p$,
\item or $x \in \dom (p)$, $p \re (A \cap L_x) \in \PP \re B$, $q \re (A \cap L_x) \leq_{\PP \re B} p\re (A \cap L_x)$,
and $p(x)$ and $q(x)$ are $\PP \re B$-names for conditions in $\dot \QQ_x$ such that $q \re (A \cap L_x) \forces_{\PP \re B} q(x) \leq_{\dot \QQ_x} p(x)$.
\end{itemize}
Concerning the first alternative here, note that it is easy to see that given $p \in \PP \re A$ and $x >_L \max (\dom (p))$ with $x \in A$,
there is $B \in \I_x \re A$ such that $p \in \PP \re B$. (Indeed, let $y = \max (\dom (p))$. Then there is $B_0 \in \I_y \re A$ such that $p \re (A \cap L_y) \in
\PP\re {B_0}$ and $p(y)$ is a $\PP \re B_0$-name for a condition on $\dot \QQ_y$. Now $B = B_0 \cup \{ y \} \in \I_x \re A$ and $p \in \PP \re B$.)

(2) Let $D \sub A$ and $p \in \PP \re D$. Also let $x = \max (\dom (p))$. There is $E \in \I_x \re D$ such that $p \re (D \cap L_x) \in \PP \re E$ and $p(x)$
is a $\PP \re E$-name for a condition in $\dot \QQ_x$. Let $B \in \I_x \re A$ be such that $E = B \cap D$. Then $E \sub B$ and, by induction
hypothesis (2), $p \re (D \cap L_x) = p \re (A \cap L_x) \in \PP \re B$. By induction hypothesis (5), $p(x)$ is a $\PP \re B$-name as well.
Therefore $p \in \PP \re A$.

The inclusion for the order is proved similarly.

(3) We use completeness of the embeddings (induction hypothesis (5)) and closure of the template under unions.

Assume $r \leq_{\PP \re A} q \leq_{\PP \re A} p$. Let $y$ and $x$ be the maximal elements of $\dom (r)$ and $\dom (q)$, respectively. 
There are $B_y \in \I_y \re A$ and $B_x \in \I_x \re A$ witnessing the order relationship. In particular, 
$r \re (A \cap L_y) , q \re (A \cap L_y) \in \PP \re B_y$,  $q \re (A \cap L_x) , p \re (A \cap L_x) \in \PP \re B_x$, and 
$r \re (A \cap L_y) \leq_{\PP \re B_y} q \re (A \cap L_y)$, $q \re (A \cap L_x)\leq_{\PP \re B_x}p \re (A \cap L_x) $. Let
$B = B_y \cup B_x \in \I_y \re A$. We check that $B$ witnesses $r \leq_{\PP \re A} p$.

If $x < y$, we have $x \in B$, (2) gives us $q, p \re (A \cap L_x)  \in \PP \re B$ and therefore also $p \in \PP \re B$,
and $r \re (A \cap L_y) \leq_{\PP \re B} q \leq_{\PP \re B} p$. By induction hypothesis (3), $r \re (A \cap L_y) \leq_{\PP \re B} p$,
and $r \leq_{\PP \re A} p$ follows.

If $x = y$, $r \re (A \cap L_y) \leq_{\PP \re B} q \re (A \cap L_y) \leq_{\PP \re B} p \re (A \cap L_y)$, 
and thus by induction hypothesis (3), $r \re (A \cap L_y) \leq_{\PP \re B} p \re (A \cap L_y)$, so we are done if $x \notin \dom (p)$. 
Assume $x \in \dom (p)$. Then $r(y)$ and $q(y)$ are $\PP \re B_y$-names, and $q(x)$ and $p(x)$ are $\PP \re B_x$-names.
By induction hypothesis (5), they are all $\PP \re B$-names and $r \re (A \cap L_y) \forces_{\PP \re B} r(y) \leq q(y) \leq p(y)$, as required.

(4) We consider several cases.

\underline{Case 1.} There is $x =\max (A)$ such that $A_0 : = A \cap L_x = A \sem \{ x \} \in \I_x \re A$. Then $\PP \re A$
is easily seen to be the standard two-step iteration $\PP \re A_0 \star\dot \QQ_x$ where $\dot \QQ_x$ is a $\PP \re A_0$-name, for 
$p \in \PP \re A$ iff $p \re A_0 \in \PP \re A_0$ and $p(x)$ is a $\PP \re A_0$-name for a condition in $\dot \QQ_x$.
In particular, $\PP \re A_0 \embed \PP \re A$.

\underline{Case 2a.} There is $x = \max (A)$, but $A_0 = A \cap L_x \notin \I_x \re A$. Let $p \in \PP \re A$\footnote{Note that
in this case it is not necessarily true that $\PP \re A = \PP \re A_0 \star \dot \QQ_x$}. There is
$B_0 \in \I_x \re A$ such that $p \re A_0 \in \PP \re B_0$ and either $x \notin \dom (p)$ or $p(x)$ is a $\PP \re B_0$-name.
Note $B_0 \subsetneq A_0$. Let $B : = B_0 \cup \{ x\} \subsetneq A$. Clearly $\Dp (B) < \Dp (A)$. 
(This holds because $\Dp (B) = \Dp (B_0) +1$ and $\Dp (A) \geq \Dp (B_0) + \omega$ by parts 3 and 4 of Lemma~\ref{depth}.)
Also $B_0 = B \cap L_x \in \I_x \re B$ and thus $p \in \PP \re B$.
Since $\I_x \re A$ is closed under unions, the collection of $B \subsetneq A$ with $B \cap L_x \in \I_x \re A$ is directed. 
Also note that by induction hypothesis (5), if $B \sub B'$ are of this form, then $\PP \re B \embed \PP \re B'$.
Therefore $\PP \re A$ is the direct limit of the $\PP \re B$ with $B \subsetneq A$ and $B \cap L_x \in \I_x \re A$.

\underline{Case 2b.} $A$ has no maximum. Let $p \in \PP \re A$ be a condition. There is $x > \max (\dom (p))$ with
$x \in A$ and therefore, as remarked earlier, there is $B \in \I_x \re A$ such that $p \in \PP \re B$. Clearly $\Dp (B) < \Dp (A)$.
The collection of $B \in \I \re A$ with $B \in \I_x \re A$ for some $x \in A$ is directed and therefore, using again induction
hypothesis (5), we see that $\PP \re A$ is the direct limit of the $\PP \re B$ for such $B$.

(5) Let $D \sub A$. We split the proof into cases according to (4) for $A$.

\underline{Case 1.} Let $D_0 : = D \cap A_0 = D \cap L_x \in \I_x \re D$. Since $\Dp (A_0) < \Dp (A)$, we may use the induction hypothesis  (5)
and see $\PP \re D_0 \embed \PP \re A_0$. We know already $\PP \re A \cong \PP \re A_0 \star \dot \QQ_x$. If
$x \notin D$, then $D = D_0$, and $\PP \re D \embed \PP \re A_0 \embed \PP \re A$ follows. If $x \in D$, then
$\PP \re D \cong \PP \re D_0 \star \dot \QQ_x$ where $\dot \QQ_x$ is a $\PP \re D_0$-name. Since $\QQ_x$ is Suslin
ccc, $\PP \re D \embed \PP \re A$ follows from Lemma~\ref{Suslinccc-completeness}.

\underline{Case 2a.} Assume first $D_0 = D \cap L_x \in \I_x \re D$. So there is $B_0 \in \I_x \re A$ such that
$D_0 = D \cap B_0$. Put $B : = B_0 \cup \{ x \} \subsetneq A$. Then $D \sub B$ and $\PP \re D \embed \PP \re B \embed \PP \re A$
where the first $\embed$ is by induction hypothesis (5) (because $\Dp (B) < \Dp (A)$) and the second, by Case 2a of (4) above.

So assume $D_0 \notin \I_x \re D$. Suppose first that $x \in D$. By Case 2a of (4) applied to $D$ instead of $A$,
$\PP \re D$ is the direct limit of the $\PP \re E$ where $E \subsetneq D$ with $E \cap L_x \in \I_x \re D$. 
Each such $E$ is of the form $D \cap B$ where $B \subsetneq A$ and $B \cap L_x \in \I_x \re A$.
Conversely, any $D \cap B$ is such an $E$. Using the inductive hypothesis for correctness (7), we see that
projections in all diagrams of the form $\la  \PP \re (D \cap B),  \PP \re B, \PP (D \cap B'), \PP \re B' \ra$, where $B \sub B' \subsetneq A$
with $B \cap L_x, B' \cap L_x \in \I_x \re A $, are correct.
By Lemma~\ref{correctness-limdir}, this means, however, that the direct limit of the $\PP \re E$ completely embeds into the
direct limit of the $\PP \re B$, as required.

Suppose finally that $x \notin D$. Then $D = D_0$ and, since $D_0 \notin \I_x \re D$, we must be in Case 2 for $D$ and, 
depending on whether $D$ has a maximum or not, we are either in Case 2a or Case 2b. In the first case, by (4), if
$y = \max (D)$, then $\PP \re D$ is the direct limit of $\PP \re E$ where $E \subsetneq D$ with $E \cap L_y \in \I_y \re D$.
In the second case, again by (4), $\PP \re D$ is the direct limit of $\PP \re E$ where $E \subsetneq D$ with $E \in \I_y \re D$ for some $y \in D$.
In either case, such $E$ belongs to $\I_x \re D$ (though not all $E \in \I_x \re D$ are necessarily of this form). Since 
$\PP \re E \sub \PP \re D$ by (2) and the collection of $E \in \I_x \re D$ is directed, $\PP \re D$ must in fact be the direct
limit of $\PP \re E$ where $E \in \I_x \re D$. We note again that such $E$ agree with the sets of the form
$D \cap B$ where $B \subsetneq A$ and $B \cap L_x \in \I_x \re A$.
By correctness and the inductive hypothesis for (7), we can apply Lemma~\ref{correctness-limdir}
and see that $\PP \re D \embed \PP \re A$.

\underline{Case 2b.} If $D \in \I_x \re D$ for some $x \in A$, we are done because then $D \sub B$ for some $B \in \I_x \re A$, and 
$\PP \re D \embed \PP \re B \embed \PP \re A$ by induction hypothesis (5) and Case 2b of (4) above.

So assume $D \notin \I_x \re D$ for any $x \in A$. Again, we must be in Case 2 for $D$ and, as in the last paragraph of Case 2a in (5),
we see that $\PP \re D$ is the direct limit of $\PP \re E$ where $E \in \I_x \re D$ for some $x \in A$. Using again Lemma~\ref{correctness-limdir},
we conclude that $\PP \re D \embed \PP \re A$.

(6) $\PP \re (A \cap D) \sub \PP \re A \cap \PP \re D$ is immediate by part (2). So assume $p \in  \PP \re A \cap \PP \re D$. Let $x = \max (\dom (p))$.
There are $B \in \I_x \re A$ and $E \in \I_x \re D$ such that $p \re L_x \in \PP \re B \cap \PP \re E$ and $p(x)$ is both a $\PP \re B$-name and
a $\PP \re E$-name, and thus a $\PP \re B \cap \PP \re E$-name. (To see the latter simply note that since $p(x)$ is a name for a real,
being a $\PP \re B$-name means that all Boolean values $\Bool p(x)(i) = j \Boor$ belong to $\PP \re B$, and similarly for $\PP \re E$.
Hence the Boolean values must belong to $\PP \re B \cap \PP \re E$.) Since $\Dp (B) < \Dp (A)$ and $\Dp (E) < \Dp (D) \leq \Dp (A)$, we may apply
the induction hypothesis (6) and get $\PP \re (B \cap E) = \PP \re B \cap \PP \re E$. Note that $B \cap E \in \I_x \re (A \cap D)$.
Therefore $p \in \PP \re (A \cap D)$, as required.

(7) Again we split into cases according to (4) for $A$.

\underline{Case 1.} $x = \max (A)$, $A_0 = A \cap L_x \in \I_x \re A$, and $\PP \re A = \PP \re A_0 \star \dot \QQ_x$.
If $x \notin A '$, we get $A' \sub A_0$ and thus  $\PP \re A' \embed \PP \re A_0$, and correctness follows from induction hypothesis (7). Similarly if $x \notin D$. So
we may assume $x \in A' \cap D = D'$. Then let $A'_0 = A' \cap L_x, D_0 = D \cap L_x, D_0' = D' \cap L_x$, apply the induction hypothesis (7) 
to the diagram $\la \PP \re D_0', \PP\re A_0', \PP \re D_0, \PP \re A_0 \ra$ and use that $\QQ_x$ is correctness-preserving.
(This is the only place where this assumption is needed.) 

\underline{Case 2a.} $x = \max (A)$, $A_0 = A \cap L_x \notin \I_x \re A$, and $\PP \re A$ is the direct limit of the $\PP \re B$ where 
$B \subsetneq A$ and $B \cap L_x \in \I_x \re A$. Let $p \in \PP \re A' \sub \PP \re A$. We need to show that the
projections agree, that is, that $h_{A'D'} (p) = h_{AD} (p)$. 

First assume that $D_0 = D \cap L_x \in \I_x \re D$. Then, by the discussion in (5) (Case 2a), $D \sub B$ for a $B$ as above and,
enlarging $B$ if necessary, we may assume $p \in \PP \re B$. Let $B' = A' \cap B$. By (6) we know that $p \in \PP \re B'$. By $\Dp (B) < \Dp (A)$ and
induction hypothesis (7), $h_{A'D'} (p) = h_{B'D'} (p) = h_{BD} (p) = h_{AD} (p)$, as required. If $A_0 ' = A' \cap L_x \in \I_x \re A'$, then,
by the symmetry of correctness, the same argument works.

So assume $D_0 \notin \I_x \re D$ and $A_0 ' \notin \I_x \re A'$. By the discussion in (5) (Case 2a), we know that $\PP \re D$ 
($\PP \re A'$, respectively) is the direct limit of $\PP \re (D \cap B)$ ($\PP \re (A' \cap B)$, resp.)
where $B$ is as above. Again fix such $B$ such that $p \in \PP \re B$. (6) gives us $p \in \PP \re (A' \cap B)$. Using
$ \Dp (B) < \Dp (A)$ and the induction hypothesis (7), we see that $ h_{A'\cap B , D' \cap B} (p) = h_{B, D \cap B} (p)$. 
Again by induction hypothesis (7), we have that $h_{B_0, D \cap B_0} (p) = h_{B,D \cap B} (p)$ for any $B \sub B_0 \subsetneq A$ with
$B_0 \cap L_x \in \I_x \re A$. Since $\PP \re A$ and $\PP \re D$ are the direct limits of such $\PP \re B_0$ and $\PP \re (D \cap B_0)$,
respectively, $h_{AD} (p) = h_{B,D \cap B} (p)$ follows (see Lemma~\ref{correctness-limdir}). In case $D_0 ' = D' \cap L_x \in \I_x \re D'$,
argue as in the previous paragraph to see that we may assume $D' \sub B$ and thus obtain $h_{A'D'} (p) = h_{A' \cap B, D'} (p)$, while if
$D_0 ' \notin \I_x \re D'$, $h_{A'D'} (p) = h_{A' \cap B, D' \cap B} (p)$ follows as in the previous sentence. 
In either case, $h_{A'D'} (p) = h_{AD} (p)$, and we are done.

\underline{Case 2b.} Depending on whether $D \in \I_x \re D$ or $A' \in \I_x \re A'$ for some $x \in A$, we repeat the previous argument, referring to Case 2b of (5).
\end{proof}

While the definition of the iteration along a template looks complicated, clause (4) should be seen as saying that
such iterations are recursively built up using the two simple operations of two-step iteration and direct limit -- as are
standard finite support iterations (fsi). Note in this context that an fsi is the special case where $L = \mu$ is an ordinal and
$\I_\alpha = \{ \beta \cup F : \beta \leq \alpha $ and $F \in [\alpha]^{<\omega} \}$ for $\alpha < \mu$.

\begin{lem}   \label{template-linked}
Assume $(L, \bar \I)$ is a template, and the $\QQ_x$, $x \in L$, are correctness-preserving Suslin $\sigma$-linked partial orders coded in the ground model,
$\QQ_x = \bigcup_n Q_{x,n}$ with each $Q_{x,n}$ being linked. Then, for any $A \sub L$, $\PP \re A$ is a ccc p.o.
\end{lem}

\begin{proof}
We argue in three steps.

\underline{Step 1.} By induction on $\Dp (A)$, we show that given $p \in \PP \re A$, there is $q \leq_{\PP \re A} p$ such that for all $x \in \dom (q)$
there are $B \in \I_x \re A$ and $n = n_{q,x}$ such that $q \re (A \cap L_x) \in \PP \re B$ and $q \re (A \cap L_x) \forces_{\PP \re B} q(x) \in \dot Q_{x,n}$.

Indeed, let $p \in \PP \re A$. Also let $x = \max (\dom (p))$. There is $B \in \I_x \re A$ such that $p \re (A \cap L_x) \in \PP \re B$
and $p(x)$ is a $\PP \re B$-name for a condition in $\dot \QQ_x$. Thus we may find $r \in \PP \re B$ and $n \in \omega$ with $r \leq_{\PP \re B} p \re (A \cap L_x)$ and
such that $r \forces_{\PP \re B} p(x) \in \dot Q_{x,n}$. Since $\Dp (B) < \Dp (A)$, there is $q_0 \in \PP \re B$ with 
$q_0 \leq_{\PP \re B} r$ satisfying the induction hypothesis. Let $q \in \PP \re A$ be such that $\dom (q) = \dom (q_0) \cup \{ x \}$,
$q \re (A \cap L_x) = q_0$ and $q(x) = p(x)$. Then $q$ is as required.

\underline{Step 2.} Assume $p,q \in \PP \re A$ are as in Step 1, that is, the $n_{p,x}$ and $n_{q,x}$ exist for all $x \in \dom (p)$ and $x \in \dom (q)$,
respectively. Also suppose that $n_{p,x} = n_{q,x}$ for all $x \in \dom (p) \cup \dom (q)$. Then $p$ and $q$ are compatible.

This is proved by building a common extension by recursion on $\dom (p) \cup \dom (q)$. For $x = \min ( \dom (p) \cup \dom (q))$,
$r_x \in \PP \re (A \cap L_x)$ is the trivial condition. Assume $r_x \in \PP \re (A \cap L_x)$ has been produced for some $x \in \dom (p) \cup
\dom (q)$. Let $y$ be the successor of $x$ in $\dom (p) \cup \dom (q)$ or let $y = \infty$ if $x = \max (\dom (p) \cup \dom (q))$.
In the latter case also let $L_\infty = L$. If $x \in \dom (p) \sem \dom (q)$, let $r_y \in \PP \re (A \cap L_y)$ be such that
$\dom (r_y) = \dom (r_x) \cup \{ x \}$, $r_y \re (A \cap L_x) = r_x$, and $r_y (x) = p(x)$. If $x \in \dom (q) \sem \dom (p)$, define
$r_y$ analogously. If $x \in \dom (p) \cap \dom (q)$, find $r_y \re (A \cap L_x) \leq r_x$ and $r_y (x)$ such that
$r_y \re (A \cap L_x) \forces_{\PP \re (A \cap L_x)} r_y (x) \leq p(x), q(x)$. This is possible because $n_{p,x} = n_{q,x}$.
Letting $r = r_\infty$, we see that $r$ is a common extension of $p$ and $q$.

\underline{Step 3.} ccc-ness now follows by a straightforward $\Delta$-system argument.
\end{proof}

\begin{lem}   \label{template-countable}
Let $(L, \bar \I)$ be a template. Also assume the $\QQ_x$ are as in the previous lemma. Let $A \sub L$.
\begin{enumerate}
\item If $p \in \PP \re A$, then there is a countable $C \sub A$ such that $p \in \PP \re C$.
\item If $\dot f$ is a $\PP \re A$-name for a real, then there is a countable $C \sub A$ such that $\dot f$ is a $\PP \re C$-name.
\end{enumerate}
\end{lem}

\begin{proof}
This is proved by a simultaneous induction on $\Dp (A)$.

(1) Assume $p \in \PP \re A$. Let $x = \max (\dom (p))$. There is $B \in \I_x \re A$ such that $p \re (A \cap L_x) \in \PP \re B$ and
$p(x)$ is a $\PP \re B$-name. By induction hypothesis (1), there is a countable $C_0 \sub B$ such that $p \re (A \cap L_x) \in \PP \re C_0$.
By induction hypothesis (2), since $p(x)$ is a name for a real, there is a countable $C_1 \sub B$ such that $p(x)$ is a
$\PP \re C_1$-name. Let $C = C_0 \cup C_1 \cup \{ x \}$. Then $C$ is countable and $p \in \PP \re C$.

(2) Assume $\dot f$ is a $\PP \re A$-name. For $i\in\omega$, let $\{ p_{n,i} : n \in \omega \}$ be a maximal antichain of conditions
deciding $\dot f (i)$. This uses the ccc-ness proved in the previous lemma. By part (1), there are countable $C_{n,i} \sub A$ such that $p_{n,i} \in \PP \re C_{n,i}$. 
Let $C = \bigcup_{n,i} C_{n,i}$. Then $\dot f$ is a $\PP \re C$-name.
\end{proof}

\begin{cor}
Let $(L, \bar \I)$ be a template. Also assume the $\QQ_x$ are as in Lemma~\ref{template-linked}. 
Then $\PP \re L$ is the direct limit of the $\PP \re A$ where $A \sub L$ is countable.
\end{cor}

\begin{proof}
By the previous lemma, $\PP \re L = \bigcup \{ \PP \re A : A \sub L$ is countable$\}$. Since the collection of countable subsets of
$L$ is directed, $\PP \re L = \lim \dir \{ \PP \re A : A \sub L$ is countable$\}$ follows.
\end{proof}

An easy consequence of this is for example that the limit of an fsi of Suslin ccc partial orders
can be represented as the direct limit of its countable fragments.\footnote{When iterating along a wellorder, ccc-ness is
preserved, so the $\sigma$-linkedness of Lemma~\ref{template-linked} is not needed.} More explicitly,
if $\la \PP_\alpha , \dot \QQ_\alpha : \alpha < \mu \ra$ is such an iteration, then $\PP_\mu = \lim \dir \{ \PP_A : A \sub \mu$ is countable$\}$
where $\PP_A$ is obtained by only iterating the $\dot \QQ_\alpha$ with $\alpha \in A$.

%%%

\subsection{The consistency of $\dd < \aa$ in ZFC}
\label{d<a-ZFC}

For showing the consistency of $\dd < \aa$ in ZFC, the ultrapower argument is replaced by an isomorphism-of-names argument.
Recall the following folklore result (see e.g.~\cite[Proposition 3.1]{Br02} for a proof).

\begin{prop}  \label{Cohen-mad}
Assume CH, and let $\lambda = \lambda^\omega$ be a cardinal. In the forcing extension obtained by adding $\lambda$
Cohen reals, every mad family has either size $\aleph_1$ or size $\cc=\lambda$.
\end{prop}

\underline{\sf STRATEGY.} The point of the proof of Proposition~\ref{Cohen-mad} is that using CH and a $\Delta$-system argument, if $\aleph_2 \leq \kappa < \lambda$,
and $\{ \dot A_\alpha : \alpha < \kappa \}$ is a name for an a.d. family, then one can produce 
another name $\dot A_\kappa$ isomorphic to $\aleph_2$ many $\dot A_\alpha$ and such that
$\dot A_\kappa$ is a.d. from all $\dot A_\alpha$. (More explicitly, let $K \in [\kappa]^{\aleph_2}$ be such that the names $\dot A_\alpha$, $\alpha \in K$,
are isomorphic and their supports $B_\alpha$ form a $\Delta$-system with root $R$, then we choose the isomorphic $\dot A_\kappa$ such that
its support $B_\kappa$ contains $R$ and $B_\kappa \sem R$ is disjoint from $\bigcup_{\alpha < \kappa} B_\alpha$. It is then easy to see that $\dot A_\kappa$
is forced to be a.d. from all $\dot A_\alpha$.)

By global homogeneity of Cohen forcing, the isomorphism 
producing $\dot A_\kappa$ comes from an automorphism of the whole forcing, but this is more than what is needed.
To obtain the consistency of $\dd < \aa$, it suffices to build a partial order forcing $\bb = \dd = \aleph_2$
and having sufficient local homogeneity to allow for the construction of $\dot A_\kappa$. This is exactly what the
template method achieves. Mad families of size $\aleph_1$ are ruled out in this scenario by $\bb = \aleph_2$.

\begin{lem}
Hechler forcing is a correctness-preserving Suslin ccc forcing notion: assume projections in $\la \PP_i \ra$ are correct. Then
so are projections in $\la \PP_i \star \dot \DD \ra$.
\end{lem}

\begin{proof}
By Lemma~\ref{correctness-quotient}, we may assume without loss of generality that $\PP_{0 \land 1} = \{ \zero, \one \}$.
Let $q_1 = (p_1 , ( s_1, \dot f_1 )) \in \PP_1 \star \dot \DD$ and fix $(s,f) \leq h_{1, 0 \land 1} ( q_1)$.
We may suppose that $s \supseteq s_1$. Then, given any $t \supseteq s$ such that $t$ dominates $f$
on its domain, there is $p_1' \leq p_1$ forcing $t \geq \dot f_1$ on its domain. (This is so because
$(t, g) \leq (s,f)$ with $g \re [|t| , \infty) = f \re [|t| , \infty)$ is compatible with $q_1$.) 

Now assume $q_0 = (p_0, (s_0, \dot f_0)) \in \PP_0 \star \dot \DD$ extends $(s,f)$. Thus $s_0 \supseteq s$
and $s_0$ dominates $f$ on its domain. By the previous paragraph find $p_1 ' \leq p_1$ forcing $s_0 \geq \dot f_1$
on its domain. By correctness, $p_0$ and $p_1 '$ are compatible in $\PP_{0 \lor 1}$ and, clearly,
the common extension $p_0 \cdot p_1 '$ forces that $(s_0, \dot f_0)$ and $(s_1, \dot f_1)$ are compatible.
Thus $(s,f) \leq h_{0\lor 1 , 0} (q_1)$ in $\PP_0 \star \dot \DD$.
\end{proof}

We now introduce the template for the proof of the main theorem (Theorem~\ref{a>d-ZFC}).

Let $\mu$ and $\lambda$ be cardinals. As usual, $\lambda^*$ denotes (a disjoint copy of) $\lambda$ with the reverse
ordering. Elements of $\lambda$ will be called {\em positive}, and members of $\lambda^*$ are {\em negative}.
Choose a partition $\lambda^* = \bigcup_{\alpha < \omega_1} S^\alpha$ such that each $S^\alpha$ is coinitial in $\lambda^*$. 
Define $L = L (\mu,\lambda)$ as follows. Elements of $L$ are non-empty finite sequences $x$ such that
$x(0) \in \mu$ and $x(n) \in \lambda^* \cup \lambda$ for $n > 0$. The order is naturally given by $x < y$ if
\begin{itemize}
\item either $x \subsetneq y$ and $y( |x|) \in \lambda$,
\item or $y \subsetneq x$ and $x ( |y|) \in \lambda^*$,
\item or $x(0) < y(0)$,
\item or, letting $n:= \min \{ m : x(m) \neq y(m) \} > 0$, $x(n) < y(n)$ in the natural ordering of $\lambda^* \cup \lambda$.
\end{itemize}
It is immediate that this is indeed a linear ordering. We identify sequences of length one with their range so that
$\mu \sub L$ is cofinal. 

Write $\gamma^*$ for the element of $\lambda^*$ corresponding to $\gamma \in \lambda$. Also define $\abs (z) \in \lambda$
for any $z \in \lambda \cup \lambda^*$ by $\abs (z)  = z$ for $z \in \lambda$ and $z = \abs(z)^*$
for $z \in \lambda^*$.

Say $x \in L$ is {\em relevant } if $|x| \geq 3$ is odd, $x(n)$ is negative for odd $n$ and positive for even $n$, 
$x ( |x| -1) < \omega_1$, and whenever $0<n < m$ are even such that $x(n), x(m) < \omega_1$, then there are $\beta < \alpha$
such that $x (n-1) \in S^\alpha$ and $x (m-1) \in S^\beta$. For relevant $x$, set $J_x = [ x \re ( |x| - 1), x)$, the interval of nodes
between $x \re ( |x| -1)$ and $x$ in the order of $L$. Notice that if $x < y$ are relevant, then either $J_x \cap J_y = \emptyset$
or $J_x \subsetneq J_y$ (in which case we also have $|y| \leq |x|$, $x \re ( |y| -1 ) = y \re ( |y| -1)$ and $x (|y| -1) \leq
y ( |y| -1)$).

For $x \in L$, let $\I_x$ consist of finite unions of 
\begin{itemize}
\item sets of the form $L_\alpha$, where $\alpha \leq x$ and $\alpha \in \mu$, 
\item sets of the form $J_y$, where $y \leq x$ is relevant, and
\item $\emptyset$ and singletons.
\end{itemize}

\begin{lem}
$(L, \bar\I = \{ \I _x : x \in L \} )$ is a template.
\end{lem}

\begin{proof}
By definition $\I_x$ contains singletons, is closed under unions, and $\I_x \sub \I_y$ for $x \leq y$. 
Closure under intersections follows easily from the discussion immediately preceding the definition of $\I_x$.
Hence it suffices to show that $\I := \bigcup_{x\in L} \I_x \cup \{ L \}$ is wellfounded.

Assume $A_n$, $n \in \omega$, is a decreasing chain from $\I$. Let $\alpha_n$ be such that $L_{\alpha_n}$ occurs
in $A_n$ as a component. The $\alpha_n$ must be decreasing and therefore eventually constant. This means it suffices to
consider the $J_x$ components of the $A_n$ and we may as well assume without loss of generality that $A_0 = J_{x^{\la\ra}}$,
and that there is a finitely branching tree $T \sub \omlom$ such that $A_n = \bigcup_{\sigma \in T \cap \omega^n}
J_{x^\sigma} \cup F_n$ where the $F_n \sub L$ are finite, such that $\sigma \sub \tau$
implies $J_{x^\tau} \sub J_{x^\sigma}$, and such that the $J_{x^\sigma}$, $\sigma \in T \cap \omega^n$,
are pairwise disjoint. Now note that if $f \in [T]$ is a branch, then the sequence $\{ x^{f \re n} : n \in \omega \}$ must eventually
stabilize. Indeed, if $| x^{f \re n} | \to \infty$, then $\{ \alpha : x^{f\re n} ( | x^{f \re n} | -2 ) \in S^\alpha$ for some $n \}$ would constitute a 
decreasing sequence of ordinals, by the definition of ``relevant", a contradiction. Therefore $| x^{f \re n} |$ is eventually constant.
But then the decreasing sequence $x^{f \re n} ( | x^{f \re n} | - 1)$ must be eventually constant as well, and so must be $x^{f \re n}$.
Since $T$ is a finitely branching tree this means that the total number of $x^\sigma$ is finite which in turn implies that the
sequence of the $A_n$ eventually stabilizes.
\end{proof}

Note that, ordered by inclusion, $L$ is a tree of countable height. Countable subtrees $A ,B \sub L$ are called {\em isomorphic}
if there is a bijection $\varphi = \varphi_{A,B} : A \to B$ such that for all $x, y \in A$ and all $n \in \omega$,
\begin{itemize}
\item $| \varphi (x) | = |x|$,
\item $\varphi (x) \re n = \varphi (x \re n)$,
\item $x < y$ iff $\varphi (x) < \varphi (y)$,
\item $x(n)$ is positive iff $\varphi (x) (n)$ is positive, 
\item $\QQ_x = \QQ_{\varphi (x)}$, and
\item $\varphi$ maps $\I \re A$ to $\I \re B$.
\end{itemize}
Since the trace of $\I$ on any countable set is countable, there are at most $\cc$ many isomorphism types of trees. Note that,
in view of the last two clauses, if $A$ and $B$ are isomorphic, then so are $\PP \re A$ and $\PP \re B$, for the
partial order only depends on the structure of the template and on the iterands. If only the first four clauses hold,
we call the trees {\em weakly isomorphic}.

\begin{thm}[Shelah~\cite{Sh04}]  \label{a>d-ZFC}
Assume CH. Let $\lambda > \mu > \aleph_1$ be regular cardinals with $\lambda^\omega = \lambda$. Then there is
a ccc forcing extension satisfying $ \bb = \dd = \mu$ and $\aa = \cc = \lambda$.
\end{thm}

\begin{proof}
Take the template $(L , \bar \I)$ introduced above. Let $\PP = \PP \re L$ be the iteration of Hechler forcing
along this template, that is, $\QQ_x = \DD$ for all $x \in L$ in Definition and Theorem~\ref{template-main}.
Using the description of $\PP \re A$ as a two-step iteration or direct limit in (4) of the latter, it is easy to prove
by induction on $\Dp (A)$, for $A \sub L$ with $|A| \geq 2$, that $\PP \re A$ has size $|A|^\omega$ and that there are $|A|^\omega$ many
$\PP\re A$-names for reals. Thus $|\PP| = \lambda^\omega = \lambda$ and $\PP$ forces $\cc \leq \lambda$.

Also, letting $\dot d_\alpha$, $\alpha < \mu$, be the $\PP$-name of the Hechler generic added at stage $\alpha$,
we see that the $\dot d_\alpha$ form a scale of length $\mu$. Indeed, if $\alpha < \beta$, then since
$\dot d_\beta$ is generic over $\PP \re L_\beta$ and $\alpha \in L_\beta$, $\dot d_\beta$
dominates $\dot d_\alpha$. Furthermore, if $\dot x$ is an arbitrary name for a real, by Lemma~\ref{template-countable},
there is a countable $A \sub L$ such that $\dot x $ is a $\PP \re A$-name. Choosing $\alpha$ such that
$A \sub L_\alpha$ and recalling $L_\alpha \in \I$, we see that $\dot d_\alpha$ dominates $\dot x$.
Thus $\bb = \dd = \mu$ follows.

We are left with showing $\aa \geq \lambda$. Since $\bb \leq \aa$ in ZFC, we already know $\aa \geq \mu$.
Thus let $\dot \A$ be a name for an almost disjoint family of size $< \lambda$ and $\geq \mu$, say
$\dot \A = \{ \dot A^\alpha : \alpha < \kappa \}$ where $\kappa \geq \omega_2 \cdot 2$ (the latter choice is for later
pruning arguments). By Lemma~\ref{template-countable}, there are countable $B^\alpha \sub L$ such that
the $\dot A^\alpha$ are $\PP \re B^\alpha$-names. More explicitly, letting $\{ p_{n,i}^\alpha : n \in \omega \}$, $i\in\omega$,
be maximal antichains and $\{ k^\alpha_{n,i} \in \{ 0,1 \} : i , n \in \omega \}$ be such that
$p_{n,i}^\alpha  \forces i \in \dot A^\alpha$ iff $k^\alpha_{n,i} = 1$ and $p_{n,i}^\alpha  \forces i \notin \dot A^\alpha$ iff $k^\alpha_{n,i} = 0$,
we have $\{ p_{n,i}^\alpha : i , n \in \omega \} \sub \PP \re B^\alpha$. We may also assume all $B^\alpha$'s are trees.
Letting $B: = \bigcup_{\alpha < \kappa} B^\alpha$ we see that $|B| < \lambda$. By CH and the $\Delta$-system lemma
we may also assume that $\{ B^\alpha : \alpha < \omega_2 \}$ forms a $\Delta$-system with root $R$ and that
\begin{itemize}
\item $\varphi_{\alpha,\beta} : B^\alpha \to B^\beta$ is an isomorphism of trees (as defined above) fixing $R$ pointwise,
\item the induced isomorphism $ \psi_{\alpha,\beta} : \PP \re B^\alpha \to \PP \re B^\beta$ maps $p^\alpha_{n,i}$ to $p^\beta_{n,i}$,
\item there are numbers $k_{n,i}$ such that $k^\alpha_{n,i} = k_{n,i}$ for all $\alpha < \omega_2$,
\item there is some $\theta_0 < \omega_1$ such that whenever $\alpha < \omega_2$, $x \in B^\alpha$, $j$ odd, and $x(j) \in \lambda^*$, then
   $x(j) \in S^\theta$ for some $\theta < \theta_0$.
\end{itemize}
Note that we then have $\varphi_{\beta, \alpha} = \varphi^{-1}_{\alpha,\beta}$ and $\varphi_{\alpha, \gamma} = \varphi_{\beta,\gamma} \circ \varphi_{\alpha,\beta}$, and 
similarly for the $\psi_{\alpha,\beta}$.
Further notice that the second and third clauses immediately imply that $\psi_{\alpha,\beta}$ also maps the name $\dot A^\alpha$ to $\dot A^\beta$.

For $\alpha < \omega_2$, write $B^\alpha = \{ x^\alpha_s : s \in T \}$ where $T \sub (\omega_1^* \cup \omega_1)^{< \omega}$ is the canonical tree 
weakly isomorphic to any $B^\alpha$. This means in particular that $|s | = | x^\alpha_s |$, that $s(n)$ is positive iff $x^\alpha_s (n)$
is positive, and that $\varphi_{\alpha, \beta} (x^\alpha_s) = x^\beta_s$. Let $S \sub T$ be the subtree corresponding to the root $R$,
that is, $s \in S$ iff $x^\alpha_s \in R$ for any $\alpha < \omega_2$. So, for $\alpha \neq \beta$, $x^\alpha_s = x^\beta_s$ iff $s \in S$. List
the immediate successors of $S$ in $T$ as $\{ t_n : n \geq 1 \}$, i.e., $\{ t_n : n \geq 1 \} = \{ t \in T \sem S : t \re ( |t| - 1) \in S \}$.
For $\alpha < \beta < \omega_2$ define
\[ F ( \{ \alpha, \beta \} ) = \left\{ \begin{array}{ll} n & \mbox{ if  $\abs( x^\alpha_{t_n} ( |t_n|-1) ) > \abs (x^\beta_{t_n} ( |t_n|-1))$ } \\
%   & \mbox{ or $t_n ( |t_n| - 1) \in \omega_1^*$ and $x^\alpha_{t_n} ( |t_n|-1) < x^\beta_{t_n} ( |t_n|-1) $ } \\
   & \mbox{ if such $n$ exists and is minimal with this property } \\
   0 & \mbox{ otherwise } \end{array} \right. \]
Note that, by wellfoundedness of the ordinals, for every $n \geq 1$, any subset of $\omega_2$ homogeneous in color $n$ must be finite.
Hence, by the Erd\H os-Rado Theorem, we obtain a subset of size $\omega_2$ homogeneous in color $0$ and may as well assume
that $\omega_2$ itself is $0$-homogeneous. Using further pruning arguments, we may additionally suppose that 
if $s \in S$ and pairs $(\zeta , \xi) \in (\omega_1^*)^2 \cup \omega_1^2$ with $\abs (\zeta) < \abs (\xi)$ and $s\ha \zeta , s \ha \xi \in T \sem S$ (so $s \ha \zeta = t_n, s \ha \xi  = t_m$,  
for some $n \neq m \geq 1$), then for all $\alpha < \beta < \omega_1$,
%\begin{itemize}
$\abs (x^\alpha_{s \ha \zeta} ( |s|) ) < \abs ( x^\beta_{s \ha \zeta} ( |s|) ) $,  all $\abs ( x^\alpha_{s \ha \zeta}
   (|s|))$ are larger than $\omega_1$ $(\star)$, and
   \begin{itemize}
   \item either $\abs ( x^\beta_{s \ha \zeta} (|s|) ) < \abs ( x^\alpha_{s \ha \xi} (|s|) )$ (this is the case when $\sup_{\alpha < \omega_1} \abs ( x^\alpha_{s \ha \zeta} (|s|) ) <
      \sup_{\alpha < \omega_1} \abs ( x^\alpha_{s \ha \xi} (|s|))$),
   \item or $\abs ( x^\alpha_{s \ha \xi} (|s|)) < \abs ( x^\beta_{s \ha \zeta} (|s|) ) $ (this is the case when $\sup_{\alpha < \omega_1} \abs ( x^\alpha_{s \ha \zeta}  (|s|) ) =
      \sup_{\alpha < \omega_1} \abs ( x^\alpha_{s \ha \xi} (|s|) )$).
   \end{itemize}
%\item if $\zeta$ is negative, then $x^\alpha_{s \ha \zeta} ( |s|) > x^\beta_{s \ha \zeta} ( |s|) $,  and  if $\zeta > \xi$ then
%   \begin{itemize}
%   \item either $x^\beta_{s \ha \zeta} (|s|) > x^\alpha_{s \ha \xi} (|s|)$ (this is the case when $\inf_{\alpha < \omega_1} x^\alpha_{s \ha \zeta} (|s|) >
%      \inf_{\alpha < \omega_1} x^\alpha_{s \ha \xi} (|s|)$),
%   \item or $x^\alpha_{s \ha \xi} (|s|) > x^\beta_{s \ha \zeta} (|s|)$ (this is the case when $\inf_{\alpha < \omega_1} x^\alpha_{s \ha \zeta}  (|s|)=
%     \inf_{\alpha < \omega_1} x^\alpha_{s \ha \xi} (|s|)$).
%   \end{itemize}
%\end{itemize}
Define $x^\kappa_s \in L$ by recursion on the length of $s \in T$, as follows. If $s \in S$, then let $x^\kappa_s = x^\alpha_s$ for any $\alpha < \omega_1$
(in particular, $|x^\kappa_s| = |x^\alpha_s|= |s|$). If $s \in S$ and $s \ha \zeta \notin S$, we will have $| x^\kappa_{s \ha \zeta} | = | s \ha \zeta| + 2$.
First let $x^\kappa_{s \ha \zeta} ( |s|) $ be the limit of the $x^\alpha_{s \ha \zeta} ( |s|)$ (so it is either the sup or the inf, depending on whether $\zeta$ is
positive or negative). Next find $\gamma < \lambda$ with $\gamma > \omega_1$ and $\gamma^* \in S^{\theta_0}$, such that for all $s$ and $\zeta$,
%\begin{itemize}
and all $y \in B$ with $y \re ( |s| + 1) = x^\kappa_{s \ha \zeta} \re
   ( |s| + 1)$, we have $\abs ( y (|s| + 1) ) < \gamma$. 
%\item if $x^\kappa_{s \ha \zeta} ( |s| ) = \inf_{\alpha < \omega_1} x^\alpha_{s \ha \zeta} (|s|)$, then for all $y \in B$ with $y \re ( |s| + 1) = x^\kappa_{s \ha \zeta} \re
%   ( |s| + 1)$, we have $y (|s| + 1) < \gamma$.
%\end{itemize}
It is clear that such a $\gamma$ exists because $\lambda > |B|$ is regular. If $\zeta$ (and $x^\kappa_{s \ha \zeta} ( |s| )$ ) is positive, let $x^\kappa_{s \ha \zeta} ( |s| + 1) = \gamma^*$, and if
$\zeta$ is negative, $x^\kappa_{s \ha \zeta} ( |s| + 1) = \gamma$ $(\star\star)$. To complete the definition of $x^\kappa_{s \ha \zeta}$ define
\[ x^\kappa_{s \ha \zeta} ( |s| + 2) = \left\{ \begin{array}{ll} x^0_{s \ha \zeta} (|s|) & \mbox{ if } |s| > 0 \\
   \xi + 2n + 1 &  \mbox{ if } |s| = 0 \mbox{ and } \zeta = \xi + n \mbox{ with $\xi$ limit or $\xi = 0$} \end{array} \right.   \]
Finally, for the remaining $t \in T$, stipulate again that $| x^\kappa_t| = | t | + 2$, find $s \subsetneq t$ with $s \in S$ maximal, put $x^\kappa_t \re ( |s| + 3)
= x^\kappa_{s \ha t(|s|) }$ and $x^\kappa_t (j+2) = x^0_t (j)$ for $j > |s|$.

Let $B^\kappa = \{ x^\kappa_s : s\in T \}$. Notice that $B^\kappa$, though very tree-like, is not a tree like the $B^\alpha$'s.
For $\alpha < \omega_1$ define $\varphi_{\alpha, \kappa} : B^\alpha \to B^\kappa$ by $\varphi_{\alpha,\kappa} (x^\alpha_s) = x^\kappa_s$ for $s \in T$.
We proceed to show that  $\PP \re B^\alpha$ and $\PP \re B^\kappa$ are isomorphic by the  map $\psi_{\alpha,\kappa}$ induced by $\varphi_{\alpha,\kappa}$,
which almost maps $\I \re B^\alpha$ to $\I \re B^\kappa$ (in the sense explained below). 
It suffices to consider the case $\alpha = 0$ as $\PP \re B^\alpha \cong \PP \re B^0$. Clearly, $\varphi = \varphi_{0,\kappa} $ is order-preserving, and it
is sufficient to figure out the effect of $\varphi$ and its inverse on the $L_\beta$ and the $J_x$.

First fix $\beta$ and consider $L_\beta$. Note that there is $\beta_0 \leq \beta$ such that $\varphi  (L_{\beta_0} \cap B^0) = L_{\beta_0} \cap B^\kappa =
L_\beta \cap B^\kappa$. For any $s \in T$ with $x_s^0 \in L_\beta$ yet $x_s^\kappa \notin L_\beta$, we must have
$x_s^\kappa (0) > \beta \geq x_s^0 (0) \geq \beta_0$ and $x^\kappa_s (0) = \sup_{\alpha < \omega_1} x^\alpha_s (0)$. In particular,
for all such $s$, $x^\kappa_s (0)$ must have the same value, say $\gamma_0$. Also
$x^\kappa_s (1) = \gamma^*$ and $x^\kappa_s (2) = \xi + 2n + 1  < \omega_1$ where $s (0) = \xi +n$ with $\xi$ limit. 
If, for some $s \in T$, $x_s^0 (0) = \beta$, let $\eta = \xi + 2n + 1$ where $s(0) =\xi + n$ with $\xi$ limit. If there is no such $s$
and $\xi + n = \sup \{ s(0) + 1: x^0_s (0) < \beta \}$, $\xi$ limit, let $\eta = \xi + 2n $. 
Then we see that $L_\beta \cap B^0$ is mapped to
$(L_\beta \cup J_x) \cap B^\kappa$ via $\varphi$, where $|x| = 3$, $x(0) = \gamma_0$, $x(1) = \gamma^*$, and $x(2) = 
\eta $ (note that this $x$ is indeed relevant).  

Next assume $x$ is relevant and consider $J_x$. Assume that  $J_x \cap B^0 \neq \emptyset$. Then there must be $s \in T$ such that
$|s| = |x|-1$ and $x^0_s = x \re ( |x| -1 )$. In case $s \in S$, we have $x^\kappa_s = x^0_s$ and $J_x \cap B^0$ is mapped to
$J_x \cap B^\kappa$ via $\varphi$ because, by $(\star)$, we must have $y \in R$ for any $y \in B^0$ with $|y| = |x|$,
$y \re ( |x| - 1) = x^0_s$ and $y (|x| -1) \leq x (|x| - 1) < \omega_1$. In case $s \in T \sem S$, let $j_0 < |s|$ be maximal with $s \re j_0 \in S$.
Define $y$ by $|y| = |x| + 2$, $y \re ( |y| - 1) = x^\kappa_s$ and $y (|y| - 1) = x ( |x| -1)$ and note that $J_x \cap B^0$ gets mapped to
$J_y \cap B^\kappa$ via $\varphi$ provided we can show that $y$ is relevant. In case $j_0 > 0$, this follows because whenever
$x^0_s (j) > \omega_1$ where $j \geq j_0$ is even then also $x^\kappa_s (j + 2) = x^0_s (j) > \omega_1$, and, if $j_0$ is even,
we additionally have $x_s^\kappa (j_0) = \sup_{\alpha < \omega_1} x^\alpha_s (j_0) > \omega_1$ while, if $j_0$ is odd,
we additionally have $x^\kappa_s (j_0 + 1) = \gamma > \omega_1$. In case $j_0 = 0$ this is true because 
$x^\kappa_s (1) \in S^{\theta_0}$ and $\theta_0$ is larger than all the $\theta$ for which $x^\kappa_s (j) \in S^\theta$ where $j > 1$ is odd.
%We leave it to the reader to verify that similar arguments show that if $J_x \cap B^\kappa \neq \emptyset$ then there is a
%relevant $y$ such that $\varphi ( J_y \cap B^0) = J_x \cap B^\kappa$. 

On the other hand, assume $J_x \cap B_\kappa \neq \emptyset$. Letting again $s \in T$ with $|s| = |x| -1$ and $x^\kappa_s = x \re ( |x| - 1)$,
we conclude as in the previous paragraph in case $s \in S$. Let $s \in T \sem S$ and $j_0 < |s|$ maximal with $s \re j_0 \in S$.
If additionally $| x| \geq 5$, let $y$ be such that $|y| = |x| - 2$, $y \re ( |y| - 1) = x^0_s$, and $y ( |y| - 1) = x ( |x| - 1)$, and check 
that $\varphi$ maps $J_y \cap B^0$ to $J_x \cap B^\kappa$ as in the previous paragraph. In case $|x| = 3$, we must be in the situation, explained above,
that $\varphi$ maps $L_\beta \cap B^0$ to $(L_\beta \cup J_x) \cap B^\kappa$ for some $\beta$. This is the only case where the templates $\I \re B^0$
and $\I \re B^\kappa$ are not identified via $\varphi$,
for $\varphi^{-1} (J_x \cap B^\kappa)$ need not belong to $\I \re B^0$. However, note that only big sets in the template matter for the definition
of the iteration, and since $\varphi^{-1} (J_x \cap B^\kappa) \sub \varphi^{-1} ( (L_\beta \cup J_x) \cap B^\kappa ) = L_\beta \cap B^0$, it is easy to see
that we can conclude that $\PP \re B^0 \cong \PP \re B^\kappa$, as witnessed by $\psi_{0,\kappa}$ (for a more formal argument see~\cite[Lemma 1.7]{Br02}).

As mentioned already this means that $\psi_{\alpha,\kappa}$ is an isomorphism of $\PP \re B^\alpha$ and $\PP \re B^\kappa$,
and we can define $\dot A^\kappa$ as the $\psi_{\alpha,\kappa}$-image of $\dot A^\alpha$ (where $\alpha < \omega_1$ is arbitrary). 
More explicitly, $p^\kappa_{n,i} = \psi_{\alpha,\kappa} ( p^\alpha_{n,i})$, and $p_{n,i}^\kappa  \forces i \in \dot A^\kappa$ iff $k_{n,i} = 1$ 
and $p_{n,i}^\kappa  \forces i \notin \dot A^\kappa$ iff $k_{n,i} = 0$. 

By $(\star\star)$, it is then also clear that if $\beta < \kappa$ is arbitrary, we can find $\alpha < \omega_1$ such that 
$B^\alpha \cup B^\beta$ and $B^\kappa \cup B^\beta$ are order isomorphic via the mapping $\varphi '$ fixing nodes of $B^\beta$ and
sending the $x^\alpha_s$ to the corresponding $x^\kappa_s$ via $\varphi_{\alpha,\kappa}$. The point is that if $\alpha < \omega_1$
is large enough, then for any $s \in T \sem S$, there are no elements of $B^\beta$ between $x^\alpha_s$ and $x^\kappa_s$
(so that, in fact, this is true for all but countably many $\alpha$). Also $\varphi '$ almost maps
$\I \re  ( B^\alpha \cup B^\beta) $ to $\I \re ( B^\kappa \cup B^\beta )$ in the sense explained above so that
$\psi ': \PP \re ( B^\alpha \cup B^\beta) \to \PP \re ( B^\kappa \cup B^\beta )$ is an isomorphism fixing the name $\dot A^\beta$
and mapping the name $\dot A^\alpha$ to the name $\dot A^\kappa$. Since $ \PP \re ( B^\alpha \cup B^\beta)$ forces that
$\dot A^\alpha \cap \dot A^\beta$ is finite, $\PP \re ( B^\kappa \cup B^\beta )$ forces that $\dot A^\kappa \cap \dot A^\beta$ is finite.
As this is true for any $\beta$, $\dot \A$ is not maximal, and the proof is complete.
\end{proof}

Note that the template framework is in a sense more general than the framework of Subsection~\ref{ultrapowers-iterations},
for the proof of Theorem~\ref{a>d-meas} in Subsection~\ref{a>d-measurable} can also be done using templates.
See~\cite[Section 2]{Br02} for a discussion of this.

%%%

\subsection{Further results}
\label{templates-further}

By modifying the template of Theorem~\ref{a>d-ZFC}, Shelah also proved $\aa$ may be a singular cardinal of uncountable cofinality.

\begin{thm}[Shelah~\cite{Sh04}]  \label{a-singular-ZFC}
Assume GCH. Let $\mu > \aleph_1$ be regular and $\lambda> \mu$ singular of uncountable cofinality. Then there is
a ccc forcing extension satisfying $ \bb = \dd = \mu$ and $\aa = \cc = \lambda$.
\end{thm}

Embedding Hechler's forcing for adding a mad family of size $\aleph_\omega$ into the template framework,
the author obtained a model in which $\aa = \aleph_\omega$ has countable cofinality.

\begin{thm}[Brendle~\cite{Br03}] \label{a-cc-ZFC}
Assume CH and let $\lambda$ be a singular cardinal of countable cofinality. Then there is a ccc forcing extension
satisfying $\aa = \lambda$. In particular, $\aa = \aleph_\omega$ is consistent.
\end{thm}

A subgroup $G$ of $\Sym (\omega)$ is called {\em cofinitary} if any non-identity member of $G$ fixes only
finitely many numbers. The cardinal invariant $\aa_g$, the minimal size of a maximal cofinitary group, is a relative of $\aa$,
and similar results about it can be proved with the same techniques. For example:

\begin{thm}[Fischer and T\"ornquist~\cite{FT15}]
Assume CH and let $\lambda$ be a singular cardinal of countable cofinality. Then there is a ccc forcing extension
satisfying $\aa_g = \lambda$. In particular, $\aa_g = \aleph_\omega$ is consistent.
\end{thm}

More recently, the template technique has been used to obtain several consistency results about the {\em
${1\over2}$-independence number} $\ii_{1\over 2}$ in~\cite{BHKLSta} (see there for a definition), e.g.:

\begin{thm}[Brendle, Halbeisen, Klausner, Lischka, and Shelah~\cite{BHKLSta}]
Assume CH and let $\lambda$ be a singular cardinal of countable cofinality. Then there is a ccc forcing extension
satisfying $\ii_{1\over 2} = \lambda$. In particular, $\ii_{1 \over 2} = \aleph_\omega$ is consistent.
\end{thm}

The interest of these three results stems from the fact that for most cardinal invariants of the continuum,
it is known that they must have uncountable cofinality. Indeed, of the cardinals in the list of~\cite[p. 480]{Bl10},
all cardinals in Cicho\'n's diagram (see Subsection~\ref{COB-EUB}) except for $\cov (\N)$ (see~\cite[Sections 2.1 and 5.1]{BJ95} for the proofs),
as well as $\ee$, $\gg$, $\hh$, $\mm$, $\pp = \tt$, $\ss$, and $\uu$ have uncountable cofinality, with some of them being regular (this
is either straightforward or proved in~\cite{Bl10}; for $\ee$ this follows from Kada's results~\cite{Ka98}). 
Apart from the above and some more variations of $\aa$,
the only cardinal which is known to consistently have countable cofinality is $\cov (\N)$~\cite{Sh00}.
It is open whether the independence number $\ii$ or the reaping number $\rr$   (see~\cite{Bl10} for definitions) can have countable cofinality.

Replacing Hechler forcing in the template framework by other Suslin ccc forcings, one obtains a number of related
consistency results about the order relationship of cardinal invariants of the continuum. See~\cite[Section 4]{Br02}
for details. For example:

\begin{thm}[Brendle~\cite{Br02}] 
Assume CH. Let $\lambda > \mu > \aleph_1$ be regular cardinals with $\lambda^\omega = \lambda$. Then there is
a ccc forcing extension satisfying $ \add (\N) = \cof (\N) = \mu$ and $\aa = \cc = \lambda$.
\end{thm}

In all the template models discussed so far $\ss = \aleph_1$ (this is so because iterations of Suslin ccc forcing notions
keep $\ss$ small, see~\cite[Theorem 3.6.21]{BJ95}), 
and the question arose as to whether one could also increase $\ss$ in the template framework. Incorporating the
ultrapower construction from Section~\ref{ultrapowers}, Mej\'ia introduced iterations of non-definable ccc partial orders along templates and proved:

\begin{thm}[Mej\'ia~\cite{Me15}]
Assume GCH and let $\theta < \kappa < \mu < \lambda$ be uncountable regular cardinals with $\kappa$ measurable.
Then there is a ccc p.o. forcing $\ss = \theta$, $\bb = \dd = \mu$, and $\aa = \cc = \lambda$.
\end{thm}

The large cardinal assumption in fact can be removed:

\begin{thm}[Fischer and Mej\'ia~\cite{FM17}]
Assume GCH and let $\theta < \mu < \lambda$ be uncountable regular cardinals. Then there is a ccc p.o. forcing
$\ss = \theta$, $\bb = \dd = \mu$, and $\aa = \cc = \lambda$.
\end{thm}

As remarked in Subsection~\ref{ultrapowers-further} (see Theorem~\ref{a>u-meas}), Shelah also used the technique
of iterating ultrapowers of ccc forcing notions to obtain the consistency of $\aleph_1 < \uu < \aa$, assuming the consistency of a measurable
cardinal. It is not known whether this can be done in ZFC alone. The problem is that while non-definable ccc forcings of the type
$\LL_\U$ can be incorporated into the matrix-like framework discussed in Subsection~\ref{ultrapowers-iterations}, it is not clear
how to do this with the more complex template framework in Subsection~\ref{template-method}. However,
using a countable support iteration of proper forcing, Guzm\'an and Kalajdzievski~\cite{GKta} recently proved the
consistency of $\uu = \aleph_1 < \aa = \cc = \aleph_2$. On the other hand, whether $\dd = \aleph_1 < \aa $ is consistent is a
famous old open problem of Roitman's from the seventies.

%%%%%%%%%%%%%%%%%%%%%%%

%

%

%

%

%

%

%%%%%%%%%%%%%%%

\section{Boolean ultrapowers}
\label{boolultrapowers}

Assume $\PP$ is a ccc partial order, $\kappa$ is a strongly compact cardinal, and $\BB$ is a $\kappa^+$-cc and $< \kappa$-distributive cBa.
Given a $\kappa$-complete ultrafilter $\D$ on $\BB$ we may form the Boolean ultrapower $\Ult_\D (\PP, \BB)$ (see Subsection~\ref{boolultrapowers-po}). This is again a ccc partial order,
$\PP$ completely embeds into $\Ult_\D (\PP,\BB)$, and much of the basic theory is very similar to the ultrapowers of Section~\ref{ultrapowers}.
In particular Boolean ultrapowers of iterations are again iterations. Since there is considerable
freedom in selecting both the cBa $\BB$ and the ultrafilter $\D$, this method turns out to be more powerful and there is a lot of control as to
what can be achieved by just taking the Boolean ultrapower once. Accordingly, all results obtained with this method (Theorems~\ref{CM-compact}
through~\ref{RS3}) are obtained by finitely many Boolean ultrapowers, and sophisticated limit constructions as in Subsection~\ref{ultrapowers-iterations} become
unnecessary.

In Subsection~\ref{compactCichon} we present a proof of the result, due to Goldstern, Kellner, and Shelah~\cite{GKS19}, 
saying that it is consistent that all cardinal invariants in Cicho\'n's diagram 
simultaneously assume distinct values, assuming the consistency of four strongly compact cardinals (Theorem~\ref{CM-compact}). 
For this, the combinatorial properties $\COB$ and $\EUB$ and their behavior under Boolean ultrapowers is central (see Subsection~\ref{COB-EUB}).
Further results using Boolean ultrapowers can be found in Subsection~\ref{boolultrapowers-further}.

%\bigskip

\subsection{Boolean ultrapowers of partial orders}
\label{boolultrapowers-po}

Assume $\kappa$ is a strongly compact cardinal. 
Let $\BB$ be a $\kappa^+$-cc and $< \kappa$-distributive cBa. 
Then every $\kappa$-complete filter on $\BB$ can be extended to a $\kappa$-complete ultrafilter~\cite{KT64}. 
Let $\D$ be a $\kappa$-complete ultrafilter on $\BB$. For a p.o. $\PP$ define
\[ \F = \F (\PP , \BB) = \{ f : \dom (f) \mbox{ is a maximal antichain in } \BB , \ran (f) \sub \PP \} \]
For $f,g \in \F$, the {\em Boolean value} of $f = g$ is defined by
\[ \Bool f = g \Boor = \bigvee \{ a \land b : a \in \dom (f) , b \in \dom (g) , f(a) = g(b) \} \] 
Similarly we define Boolean values of other statements, e.g. $\Bool f \leq g \Boor$ etc. For $f \in \F$,
\[ [f] = f / \D = \{ g \in \F : \Bool f = g \Boor \in \D \} \]
is the {\em equivalence class} of $f$ modulo $\D$. The {\em Boolean ultrapower} $\Ult_\D (\PP,\BB)$ consists of
all such equivalence classes. It is partially ordered by
\[ [f] \leq [g] \; \mbox{ iff } \; \Bool f \leq g \Boor \in \D \]
As in the discussion of ultrapowers in Section~\ref{ultrapowers}, we identify $p \in \PP$ with the class $[f]$ of the constant function
$f (\one ) = p$ and think of $\PP$ as a subset of $\Ult_\D (\PP , \BB)$. In fact, since $\BB = \P (\kappa)$ is $\kappa^+$-cc and (trivially) $< \kappa$-distributive,
the ultrapower $\PP^\kappa / \D = \Ult_\D ( \PP, \P (\kappa))$ is a special case, and the following lemmata are generalizations of the corresponding
results in Section~\ref{ultrapowers}.

\begin{lem}   \label{Boolultra-embed}
If $\PP$ is $\kappa$-cc then $\PP \embed \Ult_\D (\PP,\BB)$.
\end{lem}

\begin{proof}
Like the proof of Lemma~\ref{ultra-embed}.
\end{proof}

Note that if $A \cap \D \neq \emptyset$ for all maximal antichains $A \sub \BB$ (in particular, if $\BB$ is $\kappa$-cc), then $\PP = \Ult_\D (\PP, \BB)$
so that the conclusion of the lemma trivially holds. On the other hand, if there is a maximal antichain $A \sub \BB$ with $A \cap \D = \emptyset$, then
the converse of Lemma~\ref{Boolultra-embed} holds (see the comment after the proof of Lemma~\ref{ultra-embed}).

\begin{lem}   \label{Boolultra-chain}
If $\PP$ is $\nu$-cc for some $\nu < \kappa$ then so is $\Ult_\D (\PP,\BB)$.
\end{lem}

\begin{proof}
This is like the proof of Lemma~\ref{ultra-chain}, but we provide the argument for the sake of completeness.
Let $f_\gamma \in \F (\PP,\BB)$, $\gamma < \nu$. By $< \kappa$-distributivity of $\BB$, the maximal
antichains $\dom (f_\gamma)$ have a common refinement $A$ and we may assume $\dom (f_\gamma) = A$
for all $\gamma < \nu$. By the $\nu$-cc of $\PP$, for all $a \in A$ there are $\gamma < \delta < \nu$ such that
$f_\gamma (a)$ and $f_\delta (a)$ are compatible. By $\kappa$-completeness of $\D$, there are $\gamma < \delta$
such that $\bigvee \{ a \in A : f_\gamma (a)$ and $f_\delta (a)$ are compatible$\}$ belongs to $\D$. Thus $[f_\gamma]$
and $[f_\delta] $ are compatible, as required.
\end{proof}

Again, if there is a maximal antichain $A \sub \BB$ with $A \cap \D = \emptyset$ and $\Ult_\D (\PP,\BB)$ is $\kappa$-cc then
$\PP$ is $\nu$-cc for some $\nu < \kappa$ (see the comment after the proof of Lemma~\ref{ultra-chain}).

For the remainder of this section, assume $\PP$ is ccc. Then so is $\Ult_\D (\PP,\BB)$ and $\PP$ completely embeds into $\Ult_\D (\PP,\BB)$.
As in Section~\ref{ultrapowers}, we obtain a natural description of $\Ult_\D (\PP,\BB)$-names for reals in terms of $\PP$-names for reals.
Let $A$ be a maximal antichain in $\BB$, and let $\{ p_n^a : n \in \omega \}$, $a \in A$, be $|A|$ many maximal antichains in $\PP$.
Defining $f_n : A \to \PP$ by $f_n (a) = p^a_n$ for $a \in A$  we obtain a maximal antichain $\{ [f_n] : n \in \omega \}$ in
$\Ult_\D (\PP,\BB)$. Furthermore, by distributivity of $\BB$, all maximal antichains of $\Ult_\D (\PP , \BB)$ are of this form.
Next assume we have $|A|$ many $\PP$-names $\dot x^a$ for reals in $\omom$, $a \in A$, given by
maximal antichains $\{ p^a_{n,i} : n \in \omega \}$ and numbers $\{ k_{n,i}^a : n \in \omega \}$, $ i\in \omega$ and $a \in A$,
such that
\[ p^a_{n,i} \forces_\PP \dot x^a (i) = k^a_{n,i}. \]
Then, letting $f_{n,i} (a) = p_{n,i}^a$ and defining $k_{n,i}$ to be the unique $\ell$ such that $\bigvee \{ a \in A : k^a_{n,i} = \ell \} \in \D$,
we obtain an $\Ult_\D (\PP,\BB)$-name $\dot y$ for a real given by
\[ [f_{n,i} ] \forces_{\Ult_\D (\PP,\BB) } \dot y (i) = k_{n,i} . \]
This is the {\em average} or {\em mean} of the $\dot x^a$, and we will usually write $\dot y = \la \dot x^a : a \in A \ra / \D$.
Using again the distributivity of $\BB$ we see that every $\Ult_\D (\PP,\BB)$-name for a real is of this form. The following {\L}o\'s-style
fact about the forcing relation is straightforward, but we include a proof for the sake of completeness.

\begin{lem} \label{los-forcing}
Suppose $B \sub \omom$ is a Borel set coded in the ground model.  Let $\dot y = \la \dot x^a : a \in A \ra / \D$ be an $\Ult_\D (\PP,\BB)$-name for a real.
Also assume $A'$ is a maximal antichain of $\BB$ refining $A$ and $f = \la p^{a'} : a ' \in A' \ra / \D \in \Ult_\D (\PP,\BB)$. Then
\[
f \forces_{\Ult_\D (\PP,\BB)}  \dot y \in B \; \Loleriar \; \bigvee \{ a' \in A ' : p^{a'} \forces_\PP \dot x^{a'} \in B \} \in \D
\]
where $\dot x^{a '} = \dot x^a$ for $a' \leq a$.
\end{lem}

\begin{proof}
We make induction on the complexity of $B$. First assume $B = [ s ]$ is a basic clopen set, $s \in \omlom$. 
Clearly $\Bool \dot y (i) = s (i) \Boor = \la  \Bool \dot x^{a'} (i) = s(i) \Boor : a' \in A' \ra / \D$ for each $i < |s|$. Thus we see:
\[ \begin{array}{rcl} f \forces \dot y \in [s] & \Loleriar & f \leq \bigwedge \{ \Bool \dot y (i) = s (i) \Boor : i < | s | \} \\
& \Loleriar & \bigvee \left\{ a' \in A' : p^{a'} \leq \bigwedge \{ \Bool \dot x^{a'} (i) = s(i) \Boor : i < | s| \} \right\} \in \D \\
& \Loleriar & \bigvee \left\{ a' \in A' : p^{a'} \forces \dot x^{a'} \in [s] \right\} \in \D \end{array} \]
We next deal with the complement:
\[ \begin{array}{rcl} f \forces \dot y \notin B & \Loleriar & \forall g = \la q^{a''} : a'' \in A'' \ra / \D \leq f \; \left[  \neg \; ( g \forces \dot y \in B ) \right] \\
& \Loleriar & \forall g = \la q^{a''} : a'' \in A'' \ra / \D \leq f \; \left[  \bigvee \left\{ a'' \in A'' : q^{a''} \forces \dot x^{a''} \in B \right\} \notin \D  \right] \\
& \Loleriar & \forall A'' \mbox{ refining } A' \; \left[ \bigvee \left\{ a'' \in A'' : \exists q^{a''} \leq p^{a''} \; ( q^{a''} \forces \dot x^{a''} \in B ) \right\} \notin \D  \right]\\
& \Loleriar & \forall A'' \mbox{ refining } A' \; \left[  \bigvee \left\{ a'' \in A'' : \forall q^{a''} \leq p^{a''} \; \neg ( q^{a''} \forces \dot x^{a''} \in B ) \right\} \in \D \right] \\
& \Loleriar & \bigvee \left\{ a' \in A' :  p^{a'} \forces \dot x^{a'} \notin B  \right\} \in \D  \end{array} \]
where the second equivalence is by induction hypothesis, and we let $p^{a''} = p^{a'}$ for $a'' \leq a'$. Finally, assume 
$B = \bigcap_n B_n$. Then
\[ \begin{array}{rcl} f \forces \dot y \in B & \Loleriar & \forall n \; ( f \forces \dot y \in B_n ) \\
& \Loleriar & \forall n \; \left( \bigvee \left\{ a' \in A ' : p^{a'} \forces_\PP \dot x^{a'} \in B_n \right\} \in \D  \right) \\
& \Loleriar &  \bigvee \left\{ a' \in A ' : \forall n \; ( p^{a'} \forces_\PP \dot x^{a'} \in B_n ) \right\} \in \D \\
& \Loleriar & \bigvee \left\{ a' \in A ' :  p^{a'} \forces_\PP \dot x^{a'} \in B \right\} \in \D  \end{array} \]
where the second equivalence holds again by induction hypothesis.
\end{proof}

\begin{lem}
Assume $\PP \embed \QQ$. Then $\Ult_\D (\PP,\BB) \embed \Ult_\D (\QQ,\BB)$.
Furthermore, projections in the diagram $\la \PP ,\QQ, \allowbreak \Ult_\D (\PP,\BB), \Ult_\D (\QQ,\BB) \ra$ are correct.
\end{lem}

\begin{proof}
Like the proof of Lemma~\ref{ultra-twostep}.
\end{proof}

The following result is not needed, but we include it to show that much of the theory can be developed like for ultrapowers (Section~\ref{ultrapowers}).

\begin{lem}
Let $\PP$ be a p.o. and let $\QQ$ be a Suslin ccc forcing notion. Then $\Ult_\D (\PP \star \dot \QQ,\BB) \cong \Ult_\D (\PP,\BB) \star \dot \QQ$.%\marginpar{needed?}
\end{lem}

\begin{proof}
Like the proof of Lemma~\ref{ultra-Suslin}.
\end{proof}

\begin{lem}
Let $\P = \la \PP_\gamma : \gamma \leq \mu \ra$ be an iteration. Then $\Ult_\D ( \P , \BB) = \la \Ult_\D (\PP_\gamma ,\BB) : \gamma \leq \mu \ra$ also
is an iteration. Moreover, if $\P$ has finite supports then so does $\Ult_\D (\P ,\BB)$.
\end{lem}

Here, ``iteration with finite supports" refers to the notion before Lemma~\ref{finite-supports} and not to the traditional one. 

\begin{proof}
Like the proof of Lemma~\ref{ultra-iteration}.
\end{proof}

In the main result of this section (Theorem~\ref{CM-compact}), we will apply the Boolean ultrapower operation (finitely often)
to an iteration. By the previous lemma, the result is again an iteration, though this is not really relevant for us.

%%%

\subsection{The properties $\COB$ and $\EUB$}
\label{COB-EUB}

We introduce and present the basic properties of two principles, $\COB$ and $\EUB$, which are important for preservation of
cardinal invariants. They will be used again in Section~\ref{submodels}.

Suppose we have a binary Borel relation $R$ on the Baire space $\omom$ (or the Cantor space $\twoom$) such that
\begin{itemize}
\item for all $x \in \omom$ there is $y \in \omom$ with $x R y$,
\item for all $y \in \omom$ there is $x \in \omom$ with $\neg ( x R y)$.
\end{itemize}
If $x R y$ holds, we say that {\em $y$ $R$-dominates $x$}, and if $\neg (x R y)$, {\em $x$ is $R$-unbounded over $y$}. 
We associate two cardinals with this relation $R$, the {\em unbounding number}
\[ \bb (R) : = \min \{ | F | : F \sub \omom \mbox{ is not $R$-dominated by a single } y \in \omom \} \]
and the {\em dominating number}
\[ \dd (R) : =  \min \{ | F | : F \sub \omom \mbox{ and all } x  \in \omom \mbox{ are $R$-dominated by a member of } F  \} \]
A typical example is when $R = \leq^*$, the {\em eventual domination ordering}: say $x \leq_n y$ if for all $k \geq n$, $x(k) \leq y(k)$ holds.
$\leq^* = \bigcup_n \leq_n$, and $\bb (\leq^*) = \bb$ ($\dd (\leq^*) = \dd$, respectively) is the usual unbounding (dominating, resp.) number.
We shall see more examples shortly.

Given such a relation $R$, a ccc partial order $\PP$, and cardinals $\lambda \leq \nu$ with $\lambda$ regular, we say {\em $\PP$ forces a
$<\lambda$-directed $R$-cone of bounds of size $\nu$}, $\COF (R, \PP, \lambda, \nu)$ in symbols, if there are a
$< \lambda$-directed partial order $\la S, \leq \ra$ of size $\nu$ and $\PP$-names $(\dot z_s : s \in S)$ for reals such that
for every $\PP$-name $\dot y$ for a real there is $s \in S$ such that for all $t \geq s$,
\[ \forces_\PP \dot y R \dot z_t \]
The connection between $\COF$ and the values of the cardinals $\bb (R)$ and $\dd (R)$ in the forcing extension is given by:

\begin{lem}   \label{COBcards}
$\COF (R,\PP, \lambda, \nu)$ implies that $\bb (R) \geq \lambda$ and $\dd (R) \leq \nu$ in the $\PP$-generic extension.
\end{lem}

\begin{proof}
Clearly the $(\dot z_s : s \in S)$ form a witness for $\dd (R)$ in the extension. On the other hand, by $< \lambda$-directedness,
any $\dot F \sub \omom$ of size $< \lambda$ will be bounded.
\end{proof}

We next discuss the relationship between $\COF$ for a partial order and its ultrapower.

\begin{lem}   \label{COBultra}
Assume $\lambda \leq \nu$ are regular and $\COF (R,\PP, \lambda, \nu)$.
\begin{enumerate}
\item If $\kappa < \lambda$ or $\nu < \kappa$, then $\COF (R,\Ult_\D(\PP,\BB), \lambda, \nu)$.
\item If $\lambda < \kappa$ and $\kappa \leq \nu$, then $\COF (R,\Ult_\D(\PP,\BB), \lambda, \max(\nu,\mu)^\kappa)$ where $\mu = | \BB |$.
\end{enumerate}
\end{lem}

\begin{proof}
(1) This is similar to Lemma~\ref{ultra-scales}, but we sketch the argument.
Let $(\dot z_s : s \in S)$ be $<\lambda$-directed and cofinal of size $\nu$ in the $\PP$-generic extension. 
We show this property is preserved in the $\Ult_\D (\PP, \BB)$-generic extension. To see this, let $\dot y = \la \dot x^a : a \in A \ra  / \D$
be a $\Ult_\D(\PP,\BB)$-name for a real in $\omom$. For each $a \in A$, find $s(a) \in S$ such that for all $t \geq s(a)$,
\[ \forces_\PP \dot x^a R \dot z_t \]
If $\kappa < \lambda$, directedness of $S$ gives us $s \in S$ bigger than all $s(a)$. If $\nu < \kappa$, from the completeness of $\D$
we obtain $s \in S$ such that $\bigvee \{ a \in A  :  s(a) = s \} \in \D$. In either case, for all $t \geq s$,
\[ \forces_{\Ult_\D(\PP,\BB)} \dot y R \dot z_t \]
follows by Lemma~\ref{los-forcing}.

(2) Again let $(\dot z_s : s \in S)$ be $<\lambda$-directed and cofinal of size $\nu$ in the $\PP$-generic extension. 
Letting $U = \Ult_\D (S , \BB)$ we easily see that $U$ is $< \lambda$-directed of size $\leq  \max(\nu,\mu)^\kappa$.
For $[u] \in U$, let $\dot y_{[u]} = \la \dot z_{u(a)} : a \in A \ra / \D$ where $A \sub \BB$  is a maximal antichain and $u : A \to S$.
We claim that $(\dot y_{[u]} : [u] \in U)$ is $R$-dominating in the $\Ult_\D(\PP,\BB)$-generic extension.
Taking $\dot x = \la \dot x^a: a \in A \ra / \D$ arbitrarily, there is $u : A \to S$ such that for all $a\in A$ and all $t \geq u(a)$,
\[ \forces_\PP \dot x^a R \dot z_t \]
In particular, if $[v] \in U$ with $[v] \geq [u]$ (i.e. $\Bool v \geq u \Boor \in \D$), we see that
$\bigvee \{ a \in A : \forces_\PP \dot  x^a R \dot z_{v(a)} \} \in \D$ and therefore, by Lemma~\ref{los-forcing},
\[ \forces_{\Ult_\D(\PP,\BB)} \dot x R \dot y_{[v]} \]
as required.
\end{proof}

Given a Borel relation $R$, a ccc partial order $\PP$, and a limit ordinal $\nu$, we say {\em $\PP$ forces an eventually $R$-unbounded sequence
of length $\nu$}, $\UBD (R, \PP, \nu)$ in symbols, if there are $\PP$-names $( \dot x_\alpha : \alpha < \nu )$ for reals such that for all
$\PP$-names $\dot y$ for reals there is $\alpha < \nu$ such that for all $\beta \geq \alpha$, 
\[ \forces_\PP \neg (\dot x_\beta R \dot y) \]
(i.e. $\dot y$ does not $R$-dominate $\dot x_\beta$). 

Note that for every Borel relation $R$ on the Baire space, we have the {\em dual relation}, $R^\perp$, given by 
\[ x  R^\perp y \;\;\; \Loleriar \;\;\; \neg ( y R x ) \]
It is well-known and easy to see that $\bb (R^\perp) = \dd (R)$ and $\dd (R^\perp) = \bb (R)$.  
Using duality we see that $\EUB$ is a special case of $\COB$.

\begin{lem}   \label{COB-EUB-equi}
$\COB (R^\perp, \PP, \nu, \nu)$ and $\EUB  (R,\PP,\nu)$ are equivalent.
\end{lem}

\begin{proof}
Indeed, if $S$ is $< \nu$-directed of size $\nu$, then $S$ has a cofinal subset isomorphic to $\nu$, and we may as well assume $S = \nu$.
Now, $\COB  (R^\perp, \PP, \nu, \nu)$ means that there are $\PP$-names $(\dot x_\alpha : \alpha < \nu)$ for reals  such that for every $\PP$-name $\dot y$ for a real 
there is $\alpha < \nu$ such that for all $\beta \geq \alpha$,
\[ \forces_\PP \neg ( \dot x_\beta R \dot y) \]
which is exactly $\UBD (R,\PP, \nu)$.
\end{proof}

Using the earlier results about $\COB$, we infer:

\begin{cor}   \label{EUBprops}
\begin{enumerate}
\item Assume $\nu$ is regular and $\UBD (R, \PP, \nu)$. Then $\PP$ forces that $\bb (R) \leq \nu$ and $\dd (R) \geq \nu$.
\item Assume $\nu$ is regular and $\UBD (R, \PP, \nu)$. If $\nu \neq \kappa$, then $\UBD (R, \Ult_\D(\PP,\BB), \nu)$.
\end{enumerate}
\end{cor}

We end this subsection with a couple of relations which we shall use in the next subsection as well as in Section~\ref{submodels}.
\begin{itemize}
\item Say a function $\varphi : \omega \to \omloms$ is a {\em slalom} if $|\varphi (n)| = n$ for all $n$. The slaloms can be identified with the Baire space.
   For a slalom $\varphi$ and $x \in \omom$, let \[ x R_1 \varphi \; \mbox{ if for all but finitely many } n, \; x(n) \in \varphi (n) \]
   It is well-known~\cite[Theorem 2.3.9]{BJ95} that $\bb (R_1) = \Add N$ and $\dd (R_1) = \Cof N$.
\item Let $x \in \omom$. There is a canonical way to associate a null $G_\delta$ set $N_x$ with $x$. 
   More explicitly, let $(U_i^n : i \in \omega )$ list all clopen subsets of $\twoom$ of measure $\leq 2^{-n}$, and
   put $N_x = \bigcap_m \bigcup_{n \geq m} U^n_{x(n)}$. %\marginpar{details!}
   For $x,y \in \omom$, let \[ x R_2 y \; \mbox{ if } \; y \notin N_x \]
   Then clearly $\bb (R_2) = \Cov N$ and $\dd (R_2) = \Non N$.
\item For $x, y \in \omom$, let \[ x R_3 y\;  \mbox{ if }\;  x \leq^* y \]
   Then clearly $\bb (R_3) = \bb$ and $\dd (R_3) = \dd$.
\item For $x, y \in \omom$, let \[ x R_4 y \; \mbox{ if } \; x \neq^* y \; \mbox{ if for all but finitely many } n, \; x(n) \neq y (n) \]
   It is well-known~\cite[Theorems 2.4.1 and 2.4.7]{BJ95} that $\bb (R_4) = \Non M$ and $\dd (R_4) = \Cov M$.
\end{itemize}
These cardinals can be  displayed in {\em Cicho\'n's diagram} (see~\cite[Chapter 2]{BJ95} or~\cite[Section 5]{Bl10} for details), where cardinals grow as one moves
up or right (see Figure 2).
\begin{figure}[ht]
\begin{center}
\setlength{\unitlength}{0.2000mm}
\begin{picture}(750.0000,180.0000)(10,15)
\thinlines
\put(655,180){\line(1,0){30}}
\put(475,180){\line(1,0){45}}
\put(195,180){\line(1,0){45}}
\put(375,180){\line(1,0){30}}
\put(315,100){\line(1,0){100}}
\put(515,20){\line(1,0){30}}
\put(335,20){\line(1,0){45}}
\put(235,20){\line(1,0){30}}
\put(55,20){\line(1,0){45}}
\put(160,30){\line(0,1){140}}
\put(580,30){\line(0,1){140}}
\put(440,30){\line(0,1){60}}
\put(440,110){\line(0,1){60}}
\put(300,110){\line(0,1){60}}
\put(300,30){\line(0,1){60}}
\put(550,170){\makebox(80,20){$\dd (R_1) = \cof(\mathcal{N})$}}
\put(550,10){\makebox(60,20){$\non(\mathcal{N})$}}
\put(410,170){\makebox(60,20){$\cof(\mathcal{M})$}}
\put(410,90){\makebox(60,20){$\mathfrak{d}$}}
\put(410,10){\makebox(80,20){$\dd (R_4)= \cov(\mathcal{M})$}}
\put(270,10){\makebox(60,20){$\add(\mathcal{M})$}}
\put(270,90){\makebox(60,20){$\mathfrak{b}$}}
\put(270,170){\makebox(80,20){$\bb (R_4) = \non(\mathcal{M})$}}
\put(130,170){\makebox(60,20){$\cov(\mathcal{N})$}}
\put(130,10){\makebox(80,20){$\bb (R_1) = \add(\mathcal{N})$}}
\put(20,10){\makebox(40,20){$\aleph_1$}}
\put(680,170){\makebox(40,20){$\mathfrak{c}$}}
\end{picture} \\[5mm]
Figure 2: Cicho\'n's diagram
\end{center}
\end{figure}

The following result forms the basis for the main results of both Subsections~\ref{compactCichon} and~\ref{submodelCichon}.

\begin{thm}[Goldstern, Mej\'ia, and Shelah~\cite{GMS16}, see also~\cite{GKS19}] \label{preforcing}
Assume GCH, and let $\lambda_1 \leq \lambda_2 \leq \lambda_3 \leq \lambda_4 \leq \lambda_5$ be uncountable regular cardinals. 
There is a ccc p.o. $\PP^{pre}$, the {\em preparatory forcing}, such that for $1 \leq i \leq 4$,
\begin{itemize}
\item $\EUB (R_i, \PP^{pre}, \nu)$ for every $\nu$ with $\lambda_i \leq \nu \leq \lambda_5$,
\item $\COB (R_i, \PP^{pre}, \lambda_i, \lambda_5)$
\end{itemize}
In particular, $\PP^{pre}$ forces 
\[ \add (\N) = \lambda_1 \leq \cov (\N) = \lambda_2 \leq \add (\M) = \bb = \lambda_3 \leq \non (\M) = \lambda_4 \leq  \cov (\M) = \cc = \lambda_5 \]
\end{thm}

\begin{proof}[Proof Sketch]
We first note that the values for the cardinal invariants follow from Lemma~\ref{COBcards} and part 1 of Corollary~\ref{EUBprops}.

We sketch the proof for the particular case $\lambda_4 = \lambda_5$ and then make some comments on the general case.\footnote{This 
special case has been known at least since the 90's and can be proved by an elaboration of the methods of~\cite{Br91}.} Make
a finite support iteration $(\PP_\alpha , \dot \QQ_\alpha : \alpha < \lambda_5)$ of length $\lambda_5$ of ccc partial orders, going through 
\begin{enumerate}
\item eventually different reals forcing $\EE$ cofinally often, as well as through
\item all subforcings of localization forcing of size $< \lambda_1$, 
\item all subforcings of random forcing of size $< \lambda_2$, and
\item all subforcings of Hechler forcing of size $< \lambda_3$
\end{enumerate}
using a book-keeping argument. Then $\PP^{pre} = \PP_{\lambda_5}$. For definitions of the particular forcing notions and their properties, 
see~\cite[Chapter 3 and 7.4.B]{BJ95}.
Standard arguments show that the family of partial generics added by item $i+1$ guarantees $\COB (R_i , \PP^{pre} , \lambda_i , \lambda_5)$
for $1 \leq i \leq 3$. For example, the partial Hechler generics form a witness for $\COB (R_3, \PP^{pre}, \lambda_3, \lambda_5)$. 
Similarly, the eventually different reals witness $\COB (R_4, \PP^{pre}, \lambda_4, \lambda_5)$ (since $\lambda_4 = \lambda_5$).

Fix regular uncountable $\nu \leq \lambda_5$. Let $(\dot c_\alpha : \alpha < \nu)$ be the sequence of Cohen reals added in the limit stages of the initial segment
$\PP_\nu$ of the iteration. They clearly witness $\EUB (R_i , \PP_\nu, \nu)$ for $1 \leq i \leq 3$. If $\lambda_i \leq \nu$, then standard preservation
arguments (basically going back to~\cite{Br91}) show that $\EUB (R_i , \PP_\alpha, \nu)$ holds for all $\alpha \geq \nu$, and $\EUB (R_i , \PP^{pre}, \nu)$ follows. If $i = 4$
this is trivial (by $\lambda_4 = \lambda_5$). Moreover, preservation for limit ordinals $\alpha$ is a standard argument. For successor ordinals $\alpha$,
if $i = 1$, use the fact that all $\dot \QQ_\alpha$ either are of size $< \lambda_1$ or carry a 
finitely additive measure (see~\cite{Ka89} for why $\sigma$-centered forcings and subforcings of random forcing carry such a measure),
and that this preserves $\EUB (R_1, \PP_\alpha, \nu)$. If $i = 2$, use that all $\dot \QQ_\alpha$ are either of size $< \lambda_2$ or
$\sigma$-centered and thus preserve $\EUB (R_2, \PP_\alpha, \nu)$. If $i = 3$, use that all $\dot \QQ_\alpha$ either are of size $< \lambda_3$ or are
$\EE$, which preserves $\EUB (R_3, \PP_\alpha, \nu)$ by a compactness argument (see~\cite{Mi81}).  
This completes the argument in the special case.

In case $\lambda_4 < \lambda_5$, one would like to go through all subforcings of eventually different reals forcing of size $< \lambda_4$ instead. It is not clear,
however, why this should preserve $\EUB (R_3, \PP_\alpha, \nu)$ for $\alpha \geq \nu$. For this reason the subforcings of $\EE$ have to be very
carefully chosen in a sophisticated argument, see~\cite{GMS16} and~\cite{GKS19} for details.
\end{proof}

By this theorem, all cardinal invariants on the left-hand side of Cicho\'n's diagram, except for $\add (\M) = \bb$, can be separated simultaneously. 
The latter equality always holds in finite support iterations because they force $\bb \leq \non (\M) \leq \cov (\M)$, and $\add (\M) = \min \{ \bb , \cov (\M) \}$
is a theorem of ZFC. It is harder to separate
the dual cardinals on the right-hand side. We shall present two methods for doing this, in Subsections~\ref{compactCichon} and~\ref{submodelCichon},
using the fact that we already achieved a strong form of separation, namely $\COB$, on the left-hand side.

%%%

\subsection{Compact cardinals and Cicho\'n's maximum}
\label{compactCichon}

Assume $\mu > \kappa$ is a regular cardinal (where $\kappa$ is strongly compact as before).
Let $\BB$ be the completion of $\Fn (\mu, \kappa , < \kappa)$, forcing with partial functions from $\mu$ to $\kappa$ of size $< \kappa$.
Note that $\BB$ is $\kappa^+$-cc (because $2^{< \kappa} = \kappa$) and $< \kappa$-distributive.
Let $A \sub \Fn (\mu, \kappa , < \kappa)$ be a maximal antichain and  $w : A \to \kappa$. Let $\supp (A) = \bigcup \{ \dom (a) : a \in A \}$,
the {\em support} of $A$. Clearly $| \supp (A) | \leq \kappa$. If the maximal antichain $A'$ refines $A$ we canonically
extend $w$ to $A'$ by letting $w(a') = w(a)$ where $a$ is the unique element of $A$ above $a'$, for $a' \in A'$.
If $w : A \to \kappa$ and $w': A' \to \kappa$ are two such functions we get the {\em Boolean value}
\[ \Bool w < w' \Boor = \bigvee  \{ a'' \in A'' : w (a'') < w' (a'') \} \]
where $A''$ is a common refinement of $A$ and $A'$. 
For $\delta < \mu$ let $A_\delta$ be the maximal antichain of singleton partial functions $\{ \la \delta , \xi \ra \}$, $\xi < \kappa$,
and define $v_\delta : A_\delta \to \kappa$ by $v_\delta ( \{ \la \delta , \xi \ra \} ) = \xi$. 

\begin{lem}   \label{compact-ext}
The Boolean values $\Bool v_\delta > w \Boor$ with $\delta > \sup (\supp (\dom (w)))$  form a $\kappa$-complete filter on $\BB$
and therefore can be extended to a $\kappa$-complete ultrafilter $\D$.
\end{lem}

\begin{proof}
Let $\nu < \kappa$ and $(w_\xi, \delta_\xi)$, $\xi < \nu$, be pairs such that $\delta_\xi >  \sup (\supp (\dom (w_\xi)))$, and let $A_\xi = \dom (w_\xi)$.
We need to show that $\bigwedge_{\xi < \nu}  \Bool v_{\delta_\xi} > w_\xi \Boor \neq \zero$. Enumerate $\{ \delta_\xi : \xi < \nu \}$ in increasing order
and without repetitions as $\{ \delta^\zeta : \zeta < \gamma \}$ for some $\gamma \leq \nu$. Let $C^\zeta = \{ \xi : \delta_\xi = \delta^\zeta \}$.
Construct a decreasing chain $\{ q^\zeta : \zeta \leq \gamma \}$ of conditions in $\Fn (\mu, \kappa , < \kappa)$ as follows.
$q^0$ is the trivial condition and for limit $\zeta$, $q^\zeta$ is the union of the $q^\eta, \eta < \zeta$. Assume $q^\zeta$ has been constructed
such that $\dom (q^\zeta) \sub \delta^\zeta$ and let $q^{\zeta + 1}$ be an extension such that $\dom (q^{\zeta + 1} ) \sub \delta^\zeta + 1$,
$q^{\zeta + 1} \re \delta^\zeta$ extends an element $a_\xi \in A_\xi$ for each $\xi \in C^\zeta$ and $q^{\zeta +1} (\delta^\zeta) =
\sup \{ w_\xi (a_\xi) : a_\xi \in C^\zeta \} + 1$. Clearly $q^\gamma$ is an extension of $\bigwedge_{\xi < \nu}  \Bool v_{\delta_\xi} > w_\xi \Boor$.

By strong compactness of $\kappa$, we can now extend this filter base to a $\kappa$-complete ultrafilter on $\BB$.
\end{proof}

We assume from now on that $\D$ is constructed as in this lemma.

\begin{lem} \label{EUBultra}
Assume $\UBD (R, \PP, \kappa)$.  Then $\UBD (R, \Ult_\D(\PP,\BB), \mu)$.
\end{lem}

\begin{proof}
This is like the proof of Lemma~\ref{COBultra}, but we additionally need to use the special property of the ultrafilter $\D$ given by Lemma~\ref{compact-ext}.
Let $( \dot x_\alpha : \alpha < \kappa )$ be the eventually unbounded sequence  forced by $\PP$.
For $\delta < \mu$  let $\dot y_\delta = \la \dot x_{v_\delta (a)} : a \in A_\delta \ra / \D$. We claim that $\Ult_\D (\PP, \BB)$ forces
$(\dot y_\delta : \delta < \mu)$ is an eventually unbounded sequence. Indeed, let $\dot z = \la \dot z^a : a \in A \ra / \D$ be 
a $\Ult_\D (\PP, \BB)$-name for a real in $\omom$ and let $\delta > \sup (\supp (A))$.
There is a function $w : A \to \kappa$ such that for all $a \in A$ and all $\beta \geq w(a)$,
\[ \forces_\PP \neg \dot x_\beta R \dot z^a \]
By Lemma~\ref{compact-ext}, we know that $\Bool v_\delta > w \Boor \in \D$. A fortiori 
$\bigvee \{ a \in A' : \forces_\PP \neg \dot x_{v_\delta (a)} R \dot z^a \} \in \D$ where $A'$ is a common refinement of $A_\delta$ and $A$
and therefore, by Lemma~\ref{los-forcing},
\[ \forces_{\Ult_\D (\PP, \BB)} \neg \dot y_\delta R \dot z \]
as required.
\end{proof}

$\UBD$ and $\COF$ are tools to compute the cardinal invariants $\bb (R)$ and $\dd(R)$ in $\Ult_\D(\PP,\BB)$-generic extensions.

\begin{cor} \label{cor-strat} % \marginpar{keep?}
Let $\lambda \leq \nu$ be regular and assume $\UBD (R, \PP, \lambda)$ and $\COF (R,\PP, \lambda, \nu)$ hold.
\begin{enumerate}
\item If $\lambda < \kappa \leq \nu \leq \mu = \mu^\kappa$ and, additionally, $ \UBD (R, \PP, \kappa)$ holds, 
then $\UBD (R, \Ult_\D(\PP,\BB), \lambda)$, $ \UBD (R, \allowbreak \Ult_\D(\PP,\BB), \allowbreak \mu)$, and $\COF (R,\Ult_\D(\PP,\BB), \lambda, \mu)$,
and therefore $\bb (R) = \lambda$ and $\dd(R) = \mu$ in the $\Ult_\D(\PP,\BB)$-generic extension.
\item If $\kappa < \lambda$, and, additionally, $ \UBD (R, \PP, \nu)$ holds, 
then these properties are preserved by $\Ult_\D(\PP,\BB)$, and $\bb (R) = \lambda$ and $\dd(R) = \nu$ in the $\Ult_\D(\PP,\BB)$-generic extension.
\end{enumerate}
\end{cor}

\begin{proof}
This is immediate by Lemmas~\ref{COBcards}, \ref{COBultra}, \ref{EUBultra}, and Corollary~\ref{EUBprops}.
\end{proof}

\underline{\sf STRATEGY.} Assume $\PP$ forces $\bb (R) < \dd (R)$ as witnessed by $\EUB$ and $\COB$. 
For the case the strongly  compact cardinal $\kappa$ of the ground model lies between these two values,
the first part of this corollary describes a method for further increasing $\dd (R)$ while keeping the value of $\bb (R)$ by taking the
Boolean ultrapower of $\PP$. For the case $\kappa$ is below the smaller cardinal, the Boolean ultrapower will preserve both cardinals by the second part of the corollary.
This provides us with a scenario for obtaining a model in which all cardinals in Cicho\'n's diagram are distinct. 
Namely, force the left-hand cardinals to be distinct, with strongly compact cardinals in between them, and then keep stretching
the right-hand cardinals while keeping the ones on the left, by repeatedly taking Boolean ultrapowers.

\begin{thm}[Goldstern, Kellner, and Shelah~\cite{GKS19}]   \label{CM-compact}
Assume the existence of four strongly compact cardinals is consistent. Then so is the statement that all cardinals in Cicho\'n's diagram
are distinct. More explicitly, assume GCH and let $\aleph_1 < \kappa_9 < \lambda_1 < \kappa_8 < \lambda_2 < \kappa_7 < \lambda_3 < \kappa_6
< \lambda_4 < \lambda_5 < \lambda_6 < \lambda_7 < \lambda_8 < \lambda_9$ be regular cardinals such that the $\kappa_i$ are strongly compact.
Then there is a ccc p.o. $\PP$ forcing
\[ \begin{array}{c} \add (\N) = \lambda_1 < \cov (\N) = \lambda_2 < \add (\M) = \bb = \lambda_3 < \non (\M) = \lambda_4 <  \hskip 2truecm \\  \\
\hskip 2truecm < \cov (\M) = \lambda_5 < \dd = \cof (\M) = \lambda_6 < \non (\N) = \lambda_7 < \cof (\N) = \lambda_8 < \cc = \lambda_9 \end{array} \]
\end{thm} 

\begin{proof}
Assume $\BB_j$ is the completion of $\Fn (\lambda_j, \kappa_j , < \kappa_j)$
for $6 \leq j \leq 9$, and let $\D_j$ be the $\kappa_j$-complete ultrafilter on $\BB_j$ obtained from Lemma~\ref{compact-ext}.

Let $\PP^5 : = \PP^{pre}$ be the ccc partial order from Theorem~\ref{preforcing}. 
Next let $\PP^j = \Ult_{\D_j} (\PP^{j-1}, \BB_j)$ for $6 \leq j \leq 9$. We claim that $\PP := \PP^9$ is as required by the theorem. Since the proof is
the same for the four relations $R_i$, $1 \leq i \leq 4$, we do it for $i = 3$, that is, we show that $\PP^9$ forces $\bb  = \bb (R_3) = \lambda_3$
and $\dd = \dd (R_3) = \lambda_6$. First, using part 2 of Lemma~\ref{COBultra}, part 2 of Corollary~\ref{EUBprops}  and Lemma~\ref{EUBultra}, we obtain
\begin{itemize}
\item $\COB (R_3, \PP^6, \lambda_3, \lambda_6)$
\item $\EUB (R_3, \PP^6 , \lambda_3)$ and $\EUB (R_3, \PP^6, \lambda_6)$
\end{itemize}
Then, using part 1 of Lemma~\ref{COBultra} and part 2 of Corollary~\ref{EUBprops}, we get
\begin{itemize}
\item $\COB (R_3, \PP^j, \lambda_3, \lambda_6)$
\item $\EUB (R_3, \PP^j , \lambda_3)$ and $\EUB (R_3, \PP^j, \lambda_6)$
\end{itemize}
for $7 \leq j \leq 9$. By Lemma~\ref{COBcards} and part 1 of Corollary~\ref{EUBprops}, $\bb (R_3) = \lambda_3$ and $\dd (R_3) = \lambda_6$ follow.
\end{proof}

With more work resulting in a somewhat different preparatory forcing, the large cardinal assumption can be reduced to three strongly compact cardinals instead of four, 
see~\cite{BCM21}. For a simple proof of a weaker version of Theorem~\ref{CM-compact},
based on the preparatory forcing with $\lambda_4 = \lambda_5$ whose proof is sketched above (Theorem~\ref{preforcing}),
using three strongly compact cardinals, and showing the consistency of 
\[  \aleph_1 < \add (\N)   < \cov (\N) <  \bb  < \dd  < \non (\N) < \cof (\N)  < \cc   \]
we refer the reader to~\cite{KTT18}.

%%%

\subsection{Further results}
\label{boolultrapowers-further}

The method of the previous subsection can be used to obtain some other results where many cardinal invariants simultaneously assume
distinct values. 

\begin{thm}[Kellner, Shelah, and T\u anasie~\cite{KST19}]   \label{CM-compact-a}
Assume GCH and let $\aleph_1 < \kappa_9 < \lambda_1 < \kappa_8 < \lambda_2 < \kappa_7 < \lambda_3 < \kappa_6
< \lambda_4 < \lambda_5 < \lambda_6 < \lambda_7 < \lambda_8 < \lambda_9$ be regular cardinals such that the $\kappa_i$ are strongly compact.
Then there is a ccc p.o. forcing
\[ \begin{array}{c} \add (\N) = \lambda_1 <  \add (\M) = \bb = \lambda_2 < \cov (\N) = \lambda_3 < \non (\M) = \lambda_4 <  \hskip 2truecm \\  \\
\hskip 2truecm < \cov (\M) = \lambda_5 <  \non (\N) = \lambda_6 < \dd = \cof (\M) = \lambda_7 < \cof (\N) = \lambda_8 < \cc = \lambda_9 \end{array} \]
\end{thm} 

The difference between this result and Theorem~\ref{CM-compact} is the order relationship of $\add (\M) = \bb$ and $\cov (\N)$ and, dually,
of $\cof (\M) = \dd$ and $\non (\N)$. The proof is based on a different, more involved, preparatory forcing, see also~\cite{KST19}. 

Mixing the technique of~\cite{BCM21} with small partial orders forcing specific values to $\mm$ (the smallest cardinal for which Martin's Axiom fails) and to
the pseudointersection, distributivity, and groupwise density numbers, $\pp$, $\hh$, and $\gg$ (see~\cite[Section 6]{Bl10} for definitions), in the iteration leading
to the preparatory forcing, and then taking ultrapowers, one obtains:

\begin{thm}[Goldstern, Kellner, Mej\'ia, and Shelah~\cite{GKMSta3}]   \label{CM-compact-more}
Assume GCH and let $\aleph_1 \leq \lambda_1 \leq \lambda_2 \leq \lambda_3 \leq \kappa_9 < \lambda_4 < \kappa_8 < \lambda_5 < \kappa_7 < \lambda_6 \leq
\lambda_7 \leq \lambda_8 \leq \lambda_9 \leq \lambda_{10} \leq \lambda_{11} \leq \lambda_{12}$ be cardinals such that the $\kappa_i$ are strongly compact,
the $\lambda_i$ are regular for $i \neq 9, 12$, and $\cf(\lambda_9) \geq \lambda_6$ and $\cf(\lambda_{12}) \geq \lambda_3$.
Then there is a p.o. preserving cofinalities and forcing 
\[ \begin{array}{c} \aleph_1 \leq \mm = \lambda_1 \leq \pp = \lambda_2 \leq \hh = \gg = \lambda_3 <\add (\N) = \lambda_4 <  \cov (\N) = \lambda_5 
< \add (\M) = \bb = \lambda_6 \leq  \non (\M) = \lambda_7 \leq  \hskip 2truecm \\  \\
\hskip 2truecm \leq \cov (\M) = \lambda_8 \leq   \dd = \cof (\M) = \lambda_9 \leq \non (\N) = \lambda_{10} \leq \cof (\N) = \lambda_{11} \leq \cc = \lambda_{12} \end{array} \]
\end{thm}

Composing this with collapses one gets for example:

\begin{thm}[Goldstern, Kellner, Mej\'ia, and Shelah~\cite{GKMSta3}]
Assume GCH and there are three strongly compact cardinals. Then there is a p.o. forcing
\[ \begin{array}{c} \aleph_1 < \mm = \aleph_2 < \pp = \aleph_3 < \add (\N) = \aleph_4 <  \cov (\N) = \aleph_5 
< \add (\M) = \bb = \aleph_6 <  \non (\M) = \aleph_7 < \hskip 2truecm \\  \\
\hskip 2truecm < \cov (\M) = \aleph_8 <   \dd = \cof (\M) = \aleph_9 < \non (\N) = \aleph_{10} < \cof (\N) = \aleph_{11} < \cc = \aleph_{12} \end{array} \]
\end{thm}

In still unpublished work~\cite{RSta}, Raghavan and Shelah have obtained a number of consistency results about higher cardinal invariants,
that is, cardinal invariants describing the higher Baire space $\lambda^\lambda$ for regular uncountable $\lambda$ like $\bb (\lambda)$ or $\dd (\lambda)$,
using the Boolean ultrapower technique.

\begin{thm}[Raghavan and Shelah~\cite{RSta}]   \label{RS1}
For any regular $\lambda > \omega$, $\dd (\lambda) < \aa (\lambda)$ is consistent relative to a supercompact cardinal. More specifically, 
suppose that $\aleph_0 < \lambda = \lambda^{< \lambda} < \kappa$ and that $\kappa$ is supercompact. Then there is a forcing extension in
which $\kappa < \bb (\lambda) = \dd (\lambda) < \aa (\lambda)$.
\end{thm}

For a proof sketch, as in the previous subsection, let $\BB$ be the completion of $\Fn (\mu, \kappa, < \kappa)$ for some regular $\mu > \kappa^+$.
Using the fact that $\kappa$ is supercompact, one builds a $\kappa$-complete optimal ultrafilter $\D$ on $\BB$. ``Optimality" is a technical notion
related to, but not the same as, the one in~\cite[Definition 5.8]{MS16}. Let $\PP$ be the
$\kappa^+$-stage iteration of $\lambda$-Hechler forcing with supports of size less than $\lambda$. This is the canonical p.o.
for forcing $\bb (\lambda) = \dd (\lambda) = \kappa^+$. Next, let $\QQ = \Ult_\D (\PP,\BB)$. Forcing with $\QQ$ preserves
$\bb (\lambda) = \dd (\lambda) = \kappa^+$ (this is basically the same argument as the proof of Lemma~\ref{ultra-scales},
see also part 2 of Corollary~\ref{cor-strat}) and makes
$\aa (\lambda) = \mu$ (this uses the combinatorial properties of $\D$ and is the core of the argument).

Further results of theirs include:

\begin{thm}[Raghavan and Shelah~\cite{RSta}]   \label{RS2}
Suppose that $\aleph_0 < \lambda = \lambda^{< \lambda} < \kappa$ and that $\kappa$ is supercompact. Then there is a forcing extension in
which $\kappa < \bb (\lambda) < \dd (\lambda) < \aa (\lambda)$.
\end{thm}

\begin{thm}[Raghavan and Shelah~\cite{RSta}]   \label{RS3}
Suppose that $\lambda  < \kappa$, that $\kappa$ is supercompact, and that $\lambda$ is Laver indestructible supercompact. Then there is a forcing extension in
which $\lambda$ is still supercompact and $\kappa < \uu (\lambda) < \aa (\lambda)$.
\end{thm}

%%%%%%%%%%%%%%%%%%%%%%%

%

%

%

%

%

%

%%%%%%%%%%%%%%%

\section{Submodels}
\label{submodels}

Assume $\PP$ is a ccc partial order, $\kappa$ is a regular uncountable cardinal, and $N $ is a $< \kappa$-closed elementary substructure of
$H (\chi)$ containing $\kappa$ and $\PP$, where $\chi$ is a large enough regular cardinal.
Then the restriction of $\PP$ to $N$, $\PP \cap N$, is again a ccc partial order, and $\PP \cap N$ completely embeds into $\PP$.
In fact, if we choose $N$ sufficiently carefully (typically $N$ is the union of a chain of submodels), then
$\PP \cap N$ reflects the combinatorial properties of $\PP$, albeit with possibly different cardinals as witnesses.
We will introduce this method and discuss its effect on $\COB$ and $\EUB$ in Subsection~\ref{submodel-method}.

In Subsection~\ref{submodelCichon} we present a proof of the result, due to Goldstern, Kellner, Mej\'ia, and Shelah~\cite{GKMSta1}, 
saying that on the basis of ZFC  it is consistent that all cardinal invariants in Cicho\'n's diagram 
simultaneously assume distinct values (Theorem~\ref{CM-ZFC}). 
Further results using the submodel method are presented in Subsection~\ref{submodel-further}.

\subsection{The submodel method}
\label{submodel-method}

Assume (for the whole subsection) $\kappa$ is a regular uncountable cardinal, 
$\PP$ is a $\kappa$-cc partial order, and $N \preceq H (\chi)$ is $< \kappa$-closed with $\kappa, \PP \in N$,
where $\chi$ is a large enough regular cardinal. Then: 

\begin{lem}  \label{submodelbasic}
\begin{enumerate}
\item For every antichain $A \sub \PP$, $A \in N$ if and only if $A \sub N$. 
\item $\PP \cap N$ is $\kappa$-cc.
\item $\PP \cap N \embed \PP$.
\end{enumerate}
\end{lem}

\begin{proof}
(1) First let $A \in N$. Since $|A| < \kappa \leq |N|$ by the assumptions on $\PP$ and $N$, $A \sub N$ follows. If, on the other hand, $A \sub N$, then
by the $\kappa$-cc, $|A| < \kappa$ and by $< \kappa$-closure $A \in N$.

(2) Let $A \sub \PP \cap N$ be an antichain. By elementarity we see $A$ is an antichain of $\PP$, and $|A| < \kappa$ follows.

(3) Let $A \sub \PP \cap N$ be a maximal antichain. Again, $A$ is an antichain of $\PP$, and by (1) $A \in N$. Clearly, $N$ thinks
that $A$ is a maximal antichain of $\PP$, and therefore $A$ is maximal in $\PP$ by elementarity.
\end{proof}

An easy consequence of this lemma is that $G \sub \PP \cap N$ is a $\PP \cap N$-generic filter over $V$ iff it is a $\PP$-generic filter over $N$. Furthermore,
by (3), any $\PP \cap N$-generic filter $G$ over $V$ (or, equivalently, $H (\chi)$) can be extended to a $\PP$-generic filter $G^+$ over $V$ ($H (\chi)$, respectively).
We see immediately that $G^+ \cap N = G$ and $N [G] = N [G^+] \preceq H (\chi)^{V [G^+]}$. Using this we can establish a correspondence between
$\PP \cap N$-names $\dot x \in V$ for reals and $\PP$-names $\dot y \in N$ for reals such that $\dot x [G] = \dot y [G^+]$ and for all $p \in \PP \cap N$
and sufficiently absolute (e.g. Borel) $\varphi$, \[ p \forces_\PP \varphi (\dot y) \mbox{ iff } p \forces_{\PP \cap N} \varphi (\dot x) \]
In particular, $N[G^+] \cap \omom = V[G] \cap \omom$.

To see this, recall that a $\PP \cap N$-name $\dot x$ for a real in $\omom$ is given by maximal antichains $\{ p_{n,i} : n \in \lambda_i \} \sub \PP \cap N$ and numbers
$\{ k_{n,i} : n \in \lambda_i \}$, $i \in \omega$ and $\lambda_i < \kappa$, such that
\[ p_{n,i} \forces_{\PP \cap N} \dot x (i) = k_{n,i} \]
Since $N$ is $< \kappa$-closed, $\{ \{ p_{n,i} : n \in \lambda_i \} , \{ k_{n,i} : n \in \lambda_i \}: i\in \omega \} \in N$, and by $\PP \cap N \embed \PP$, $\dot x$ can be construed as a
$\PP$-name in $N$. On the other hand, if $\dot x \in N$ is  a $\PP$-name, that is, $\{ \{ p_{n,i} : n \in \lambda_i \} , \{ k_{n,i} : n \in \lambda_i \}: i\in \omega \} \in N$,
then by part 1 of Lemma~\ref{submodelbasic}, $\{ p_{n,i} : n \in \lambda_i \}  \sub N$ for all $i$ and $\dot x$ is a $\PP \cap N$-name.

We now investigate how $\COB$ and $\EUB$ for $\PP \cap N$ relate to $\COB$ and $\EUB$ for $\PP$.\footnote{$\COB$ and $\EUB$ were originally
defined for ccc forcing in Subsection~\ref{COB-EUB} but this does not really matter.} For a partial order $\la S, \leq \ra$,
let $\comp (S)$, the {\em completeness of $S$}, be the least $\lambda$ such that $S$ is not $\lambda$-directed.

\begin{lem}   \label{submodelCOB1}
\begin{enumerate}
\item Assume $\COB (R, \PP, \lambda, \nu)$ as witnessed by the partial order $\la S, \leq \ra \in N$. 
   Then $\COB (R, \PP \cap N, \lambda ', \nu ')$ whenever $\lambda ' \leq \comp (S \cap N)$ and $\nu ' \geq \cof (S \cap N)$.
\item Assume $\EUB (R, \PP, \nu)$ with $\nu \in N$. Then $\EUB (R, \PP \cap N, \cof (\nu \cap N) )$.
\end{enumerate}
\end{lem}

\begin{proof}
(1) Assume $( \dot z_s : s \in S )$ witnesses $\COB (R, \PP, \lambda, \nu )$. It suffices to show that $( \dot z_s : s \in S \cap N)$ witnesses $\COB (R, \PP \cap N, \lambda ', |S\cap N|)$.
For then we can replace $S \cap N$ by any $< \lambda '$-directed superset of a cofinal subset.
Let $\dot y$ be a $\PP \cap N$-name for a real. By the previous discussion, we know that $\dot y \in N$ may be construed as a $\PP$-name for a real.
Hence, by $\COB (R, \PP, \lambda, \nu )$ and elementarity, there is $s \in S \cap N$ such that, in $N$, for all $t \geq s$,
\[ \forces_\PP \dot y R \dot z_t \]
Again by the previous discussion, this means that for all $t \geq s$ in $N$ we have
\[ \forces_{\PP \cap N}  \dot y R \dot z_t \]
as required.

(2) This follows from (1) and Lemma~\ref{COB-EUB-equi}.
\end{proof}

\begin{lem}    \label{submodelCOB2}
Assume $\kappa \leq \theta \leq \mu$, $\theta$ regular. Let $\lambda \leq \nu$ with $\lambda$ regular.
Next let $(N_i : i < \theta )$ be an increasing sequence of $< \theta$-closed elementary submodels of $H (\chi)$ with $|N_i| = \mu$,
$N_i \in N_{i+1}$, and $\mu \cup \{ \mu, \lambda , \nu, \PP \} \in N_0$. Put $N = \bigcup_i N_i$. 
\begin{enumerate}
\item Assume $\EUB (R, \PP, \nu)$. Then:
\begin{enumerate}
\item If $\nu \leq \mu$ then $\EUB (R, \PP \cap N, \nu)$.
\item If $\nu > \mu$ then $\EUB (R, \PP \cap N, \theta)$.
\end{enumerate}
\item Assume $S \in N_0$ witnesses $\COB (R, \PP, \lambda, \nu)$. Then $\comp (S \cap N) \geq \min \{ \theta, \lambda\}$ and $\cof (S \cap N) \leq \min \{ \mu, \nu \}$
and therefore $\COB (R, \PP \cap N, \min \{ \theta, \lambda \} , \min \{ \mu, \nu \} )$. Moreover:
\begin{enumerate}
\item If $\nu \leq \mu$ then $\COB (R, \PP \cap N,\lambda, \nu)$.
\item If $\mu < \lambda$ then $\comp (S \cap N) = \cof (S \cap N) = \theta$ and thus $\COB (R, \PP \cap N,\theta, \theta)$.
\end{enumerate}
\end{enumerate}
\end{lem}

\begin{proof}
(1) Notice that if $\nu > \mu$ is regular, then $\cof (\nu \cap N) = \theta$. Hence this follows from part 2 of Lemma~\ref{submodelCOB1}.

(2) $\comp (S \cap N) \geq \min \{ \theta, \lambda \}$ holds because if $A \sub S \cap N$ with $|A| < \min \{ \theta, \lambda \}$,
then $A \in N$ by $< \theta$-closure of $N$, and $N$ thinks that $A$ has an upper bound by elementarity.
$\cof (S \cap N) \leq |S \cap N| = \min \{ \mu, \nu \}$ is obvious. $\COB (R, \PP \cap N, \min \{ \theta, \lambda \} , \min \{ \mu, \nu \} )$ then follows
from part 1 of Lemma~\ref{submodelCOB1}.

(a) If $\nu \leq \mu$, then $S \sub N_0 \sub N$, so that $S = S \cap N$ and $\COB (R, \PP \cap N, \lambda, \nu)$ is immediate.

(b) Assume $\mu < \lambda$. We know already  $\comp (S \cap N) \geq \min \{ \theta, \lambda \} = \theta$. By elementarity and $N_i \in N_{i+1}$, for every
$i < \theta$ there is $s \in S \cap N_{i+1}$ such that $s \geq t$ for all $t \in S \cap N_i$. Hence $\cof (S \cap N) \leq \theta$ follows,
and we must actually have $\comp (S \cap N) = \cof (S \cap N) = \theta$.
\end{proof}

\begin{lem}    \label{submodelCOB3}
Assume additionally to the assumptions of part 2 of the previous lemma that $\mu < \lambda$, that $\theta ' \geq \theta$ and $\mu ' > \mu$
are regular cardinals in $N_0$, and that $(M_j : j < \theta ') \in N_0$ is a family of $< \mu'$-closed elementary submodels of $H(\chi)$ with $|M_j | = \mu'$.
Put $M = \bigcup_j M_j$ and assume $M \preceq H(\chi)$ is $< \theta$-closed. Then $\COB (R , \PP \cap M \cap N, \theta, \theta ')$.
\end{lem}

\begin{proof} 
$\comp (S \cap M \cap N) \geq \theta$ is straightforward from $< \theta$-closure of $M \cap N$. So it suffices to show $\cof (S \cap M \cap N) \leq \theta '$.

Clearly $\comp (S \cap M_j) \geq \min \{ \mu ', \lambda \}$ and $\cof (S \cap M_j) \leq \min \{ \mu ' , \nu \}$, and therefore  $\COB (R, \PP \cap M_j , 
\min \{ \mu' , \lambda \}, \allowbreak \min \{ \mu ' , \nu \} )$. 
By part 2 (b) of the previous lemma we see that $\comp (S \cap M_j \cap N) = \cof (S \cap M_j \cap N) = \theta$. Let $T_j \sub S \cap M_j \cap N$
be cofinal in $S \cap M_j \cap N$ of size $\theta$ and let $T = \bigcup_{j < \theta '} T_j$. Then $|T| \leq \theta '$, and $T$ is cofinal in $S \cap M \cap N$.
Thus $\cof (S \cap M \cap N) \leq \theta '$, and $\COB (R , \PP \cap M \cap N, \theta, \theta ')$ follows.
\end{proof}

\underline{\sf STRATEGY.} Assume $\PP$ forces $\bb (R) < \dd (R) = \nu$ as witnessed by $\EUB$ and $\COB$.
Let $\theta < \theta ' < \mu < \mu' : = \bb (R)$ be arbitrary regular cardinals. Building first $(N_i ': i < \theta ')$ of size $\mu '$
according to Lemma~\ref{submodelCOB2}, and then analogously $(N_j : j < \theta )$ of size $\mu$ such that $(N_i ': i < \theta ') \in N_0$, 
we see by the two previous lemmata that $\PP \cap N' \cap N$ forces $\bb (R) = \theta$ and $\dd (R) = \theta '$. 
If we then further intersect $\PP$ with $N ''$ such that  $(N_i ': i < \theta ') , (N_i: i < \theta ) \in N_0 ''$ with $\mu '' > \theta '$,
we will not change these values anymore. This provides us with a scenario for obtaining a model for Cicho\'n's maximum: force the left-hand
cardinals to be distinct, of large enough value, and then ``collapse" dual pairs of cardinals to a priori given values, by repeatedly restricting
to appropriate elementary submodels.

%%%

\subsection{Cicho\'n's maximum in ZFC}
\label{submodelCichon}

We are ready to present a ZFC-proof of the consistency result of Theorem~\ref{CM-compact}.

\begin{thm}[Goldstern, Kellner, Mej\'ia, and Shelah~\cite{GKMSta1}]    \label{CM-ZFC}
Assume GCH and $(\lambda_i^* : 1 \leq i \leq 9)$ is a $\leq$-increasing sequence of uncountable cardinals with $\lambda_i^*$ regular
for $i \leq 8$ and $\lambda_9^*$ of uncountable cofinality. Then there is a ccc partial order $\PP$ forcing that
\[ \begin{array}{c} \aleph_1 \leq \add (\N) = \lambda_1^* \leq  \cov (\N) = \lambda_2^* \leq \add (\M) = \bb = \lambda_3^* \leq  \non (\M) = \lambda_4^* \leq  \hskip 2truecm \\  \\
\hskip 2truecm \leq \cov (\M) = \lambda_5^* \leq   \dd = \cof (\M) = \lambda_6^* \leq \non (\N) = \lambda_7^* \leq \cof (\N) = \lambda_8^* \leq \cc = \lambda_9^* \end{array} \]
\end{thm}

\begin{proof}
For notational convenience we will rename
\[ \lambda_1^* \leq  \lambda_2^* \leq \lambda_3^* \leq\lambda_4^* \leq\lambda_5^* \leq \lambda_6^* \leq \lambda_7^* \leq \lambda_8^* \leq \lambda_9^* \]
as
\[ \theta_7 \leq \theta_5 \leq \theta_3 \leq \theta_1 \leq \theta_0 \leq \theta_2 \leq \theta_4 \leq \theta_6 \leq \theta_8  \]
and choose regular cardinals 
\[  \mu_7 < \mu_6 < \mu_5 < \mu_4 < \mu_3 < \mu_2 < \mu_1 < \mu_0 < \lambda_5 \]
strictly larger than $\theta_8$. (Thus all cardinals with the possible exception of $\theta_8$ are regular and $\cf (\theta_8) \geq \aleph_1$.)
Next let $\lambda_i := \mu_{8 - 2i}$ for $1 \leq i \leq 4$. Let $ \PP^{pre}$ be the ccc partial order from Theorem~\ref{preforcing}
for $\lambda_i$, $1 \leq i \leq 5$.

We will construct a complete subforcing $\PP = \PP^{pre} \cap N^*$ of $\PP^{pre}$ which forces $(\bb (R_i), \dd (R_i)) = (\theta_{9-2i}, \theta_{8-2i})$
for $1 \leq i \leq 4$ and $\cc = \theta_8$.

Fix $N_{n,\alpha}$ for $0 \leq n \leq 7$ and $\alpha < \theta_n$, $N_n := \bigcup_{\alpha < \theta_n} N_{n,\alpha}$, as well as $N_8$ such that
\begin{itemize}
\item all $N_{n,\alpha}$ as well as $N_8$ are elementary submodels of $H (\chi)$ containing the sequence of cardinals, $\PP^{pre}$, as well as
   the witnesses $S_i$ of $\COB (R_i, \PP^{pre}, \lambda_i, \lambda_5)$ (from Theorem~\ref{preforcing}) for $1 \leq i \leq 4$
\item $N_{n,\alpha}$ contains $(N_{m,\beta} : m < n, \beta < \theta_m )$ and $(N_{n,\beta} : \beta < \alpha)$, and $N_8$ contains $(N_{m,\beta} : m \leq 7, \beta < \theta_m )$
\item the $N_{n,\alpha}$ are $< \mu_n$-closed of size $\mu_n$, and $N_8$ is $< \aleph_1$-closed of size $\theta_8$
\end{itemize}
Let $N^* := \bigcap_{0 \leq m \leq 8} N_m$. For $0 \leq n \leq 8$, let $\PP_n := \PP^{pre} \cap \bigcap_{0 \leq m \leq n} N_m$, and let $\PP : = \PP_8 = \PP^{pre} \cap N^*$.
Note that $\bigcap_{0 \leq m \leq n} N_m$ is again an elementary submodel of $H (\chi)$, and therefore each $\PP_n$ is a complete subforcing of $\PP^{pre}$.

We first show that for all $i$ with $1 \leq i \leq 4$, we have $\EUB (R_i , \PP, \theta_{9-2i})$ and $\EUB (R_i , \PP, \theta_{8-2i})$. Since the proof is the same for
all $i$, we do it for $i = 3$. By Theorem~\ref{preforcing} we have 
\begin{itemize}
\item $\EUB (R_3, \PP^{pre}, \mu_2)$ and $\EUB (R_3, \PP^{pre}, \mu_1)$.
\end{itemize}
By part 1 of Lemma~\ref{submodelCOB2}, we successively get 
\begin{itemize}
\item $\EUB (R_3, \PP_i, \mu_2)$ and $\EUB (R_3, \PP_i, \mu_1)$ for $i = 0,1$,
\item $\EUB (R_3, \PP_2, \mu_2)$ and $\EUB (R_3, \PP_2, \theta_2)$,
\item $\EUB (R_3, \PP_j, \theta_3)$ and $\EUB (R_3, \PP_j, \theta_2)$ for $3 \leq j \leq 8$,
%\item $\EUB (R_3, \PP_8, \theta_3)$ and $\EUB (R_3, \PP_8, \theta_2)$.
\end{itemize}
as required.

Next, we prove that for all $i$ with all $1 \leq i \leq 4$, we have $\COB (R_i, \PP, \theta_{9-2i}, \theta_{8-2i})$. Again we consider only the case $i=3$. 
By Theorem~\ref{preforcing} we have $\COB (R_3, \PP^{pre}, \mu_2, \lambda_5)$.
We first show $\COB (R_3, \PP_3, \theta_3, \theta_2)$: for $s \in \theta_0 \times \theta_1 \times \theta_2$ let $M_s := N_{0, s(0)} \cap
N_{1, s(1)} \cap N_{2,s(2)}$. Then $M_s$ is a $< \mu_2$-closed elementary submodel of $H(\chi)$ of size $\mu_2$. 
Also $( M_s : s \in \theta_0 \times \theta_1 \times \theta_2) \in N_{3,0}$, $M: = \bigcup_s M_s = N_0 \cap N_1 \cap N_2$, and
$|\theta_0 \times \theta_1 \times \theta_2 | = \theta_2$. Note that $M$ is $< \theta_1$-closed and thus $< \theta_3$-closed.
Therefore we may apply Lemma~\ref{submodelCOB3} with $\theta = \theta_3$, 
$\theta ' = \theta_2$, $\mu = \mu_3$, $\mu ' = \lambda = \mu_2$, $\nu = \lambda_5$, $N_\alpha = N_{3,\alpha}$ for $\alpha < \theta_3$, and $N = N_3$
to obtain $\COB (R_3, \PP_3, \theta_3, \theta_2)$ (note here that $\PP_3 = \PP^{pre} \cap M \cap N$). An easy application of  part 2 (a) of Lemma~\ref{submodelCOB2} yields that
$\COB (R_3, \PP, \theta_3, \theta_2)$ still holds.

By Lemma~\ref{COBcards} and part 1 of Corollary~\ref{EUBprops}, $\PP$ forces $\bb (R_i)= \theta_{9-2i}$ and $  \dd (R_i) = \theta_{8-2i}$.

Finally note that $| \PP | = | N^* | = | N_8 | = \theta_8 = \theta_8^{\aleph_0}$ and therefore by standard arguments $\PP$ forces $\cc \leq \theta_8$.
On the other hand, there is a sequence $(\dot x_\xi : \xi < \lambda_5 )$ of $\PP^{pre}$-names for distinct reals belonging to $N^*$.
Hence $(\dot x_\xi : \xi \in \lambda_5 \cap N^*)$ is a sequence of $\PP$-names for distinct reals. Since $| \lambda_5 \cap N^* | = \theta_8$,
$\PP$ forces $\cc \geq \theta_8$.
\end{proof}

%%%

\subsection{Further results}
\label{submodel-further}

Using the preparatory forcing from~\cite{KST19} (cf Theorem~\ref{CM-compact-a}) together with the submodel technique, one obtains:

\begin{thm}[Goldstern, Kellner, Mej\'ia, and Shelah~\cite{GKMSta1}] 
Assume GCH and $(\lambda_i : 1 \leq i \leq 9)$ is a $\leq$-increasing sequence of uncountable cardinals with $\lambda_i$ regular
for $i \leq 8$ and $\lambda_9$ of uncountable cofinality. Then there is a ccc partial order forcing that
\[ \begin{array}{c} \aleph_1 \leq \add (\N) = \lambda_1 \leq  \add (\M) = \bb = \lambda_2 \leq \cov (\N) = \lambda_3 \leq \non (\M) = \lambda_4 \leq  \hskip 2truecm \\  \\
\hskip 2truecm \leq \cov (\M) = \lambda_5 \leq  \non (\N) = \lambda_6 \leq \dd = \cof (\M) = \lambda_7 \leq \cof (\N) = \lambda_8 \leq \cc = \lambda_9 \end{array} \]
\end{thm}

The only difference between this result and Theorem~\ref{CM-ZFC} is the order relationship of $\add (\M) = \bb$ and $\cov (\N)$ and, dually, of 
$\cof (\M) = \dd$ and $\non (\N)$ (cf the difference between Theorems~\ref{CM-compact} and~\ref{CM-compact-a}). 

Further cardinal invariants can be included in the picture (cf Theorem~\ref{CM-compact-more}):

\begin{thm}[Goldstern, Kellner, Mej\'ia, and Shelah~\cite{GKMSta2}] 
Assume GCH and $(\lambda_i : 1 \leq i \leq 12)$ is a $\leq$-increasing sequence of uncountable cardinals with $\lambda_i$ regular
for $i \leq 11$ and $\lambda_{12}$ of uncountable cofinality. Then there is a cofinality-preserving partial order forcing that
\[ \begin{array}{c} \aleph_1 \leq \mm = \lambda_1 \leq \pp = \lambda_2 \leq \hh = \gg = \lambda_3 \leq \add (\N) = \lambda_4 \leq  \cov (\N) = \lambda_5 
\leq \add (\M) = \bb = \lambda_6 \leq  \non (\M) = \lambda_7 \leq  \hskip 2truecm \\  \\
\hskip 2truecm \leq \cov (\M) = \lambda_8 \leq   \dd = \cof (\M) = \lambda_9 \leq \non (\N) = \lambda_{10} \leq \cof (\N) = \lambda_{11} \leq \cc = \lambda_{12} \end{array} \]
\end{thm}

\begin{thm}[Goldstern, Kellner, Mej\'ia, and Shelah~\cite{GKMSta4}] 
In the previous result, letting $\lambda_1 = \aleph_1$ and $\lambda_\ss$ and $\lambda_\rr$ two regular cardinals with
$\lambda_3 \leq \lambda_\ss \leq \lambda_6$, $\lambda_9 \leq \lambda_\rr \leq \lambda_{12}$, and, for $i \leq 2$, $\lambda_\ss \in [ \lambda_{3 + i} , \lambda_{4+i} ]$
iff $\lambda_\rr \in [ \lambda_{11 - i }, \lambda_{ 12 - i }]$, there is a cofinality-preserving partial order forcing the given distribution together
with $\ss = \lambda_\ss$ and $\rr = \lambda_\rr$.
\end{thm}

While it is well-known that $\bb < \ss$ and $\rr < \dd$ are consistent~\cite{Sh84,BS89}, it is open whether one can force $\bb < \ss \leq \non (\M) \leq \cov (\M) < \rr < \dd$ in this theorem.

%%%%%%%%%%%%%%%%%%%%%%%

%

%

%

%

%

%

%%%%%%%%%%%%%%%%%%%%%%%%%%%%%%%%%%%%%%%%%%%%%%%%%%%%%%%%%%%%%%%%%%%%%%%

\bigskip
\bigskip
\bigskip

\noi{\bf Acknowledgments.} We thank Martin Goldstern, Jakob Kellner, and Diego Mej\'ia for many comments and for correcting several errors in an earlier version of this work.
We are also grateful to the referee for a very careful reading of our paper and many suggestions which considerably improved the presentation.

%%%%%%%%%%%%%%%%%%%%%%%%%%%%%%%%%%%%%%%%%%%%%%%%%%%%%%%%%%%%%%%%%%%%%%%

\end{document}